\documentclass[a4paper,11pt,preprint]{amsart}

\usepackage[a4paper, centering]{geometry}
\usepackage{pifont}

\usepackage{fancyvrb}
\usepackage{enumitem,ulem}
\usepackage{amsmath}
\usepackage{amsfonts}
\usepackage{amsthm}
\usepackage{amssymb}
\usepackage{upref}
\usepackage{color}
\usepackage[table,xcdraw,svgnames,dvipsnames]{xcolor}
\usepackage{graphicx}

\usepackage{listings}

\usepackage[T1]{fontenc}
\usepackage{bigfoot} 
\usepackage[framed]{matlab-prettifier}
\usepackage[colorlinks]{hyperref}
\hypersetup{
	allbordercolors=.,
	citecolor=NavyBlue,
	linkcolor=DarkRed,
	urlcolor=NavyBlue
}
\usepackage{array}
\usepackage{multicol}

\graphicspath{{./pics/}}

\theoremstyle{plain}
\newtheorem{thm}{Theorem}[section]
\newtheorem{lemma}[thm]{Lemma}

\newtheorem{prop}[thm]{Proposition}

\newtheorem{rem}[thm]{Remark}
\theoremstyle{definition}

\newcommand{\Pol}[1]{\mathcal P_{#1}}
\newcommand{\bo}[1]{{\bf#1}}
\newcommand{\Id}{\operatorname{\bo{Id}}}

\newcommand{\aaa}{a}
\newcommand{\bbb}{b}
\newcommand{\Usz}{{U}_{S_0}}
\newcommand{\Uszh}{\overline{{U}_{S_0}}}
\newcommand{\Usp}{{U}_{S_+}}
\newcommand{\Usm}{{U}_{S_-}}
\newcommand{\reg}{\text{reg}}
\newcommand{\sing}{\text{sing}}

\newcommand{\Ppp}{\Bbb{P}_+}

\usepackage{indentfirst}

\newcommand{\ra}{\rightarrow}
\def\R{\mathbb{R}}

\def\Pp{{\mathbb{P}_n}}
\newcommand{\vps}{\varepsilon}

\newcommand{\Om}{\Omega}

\newcommand{\lb}{\lambda}

\newcommand{\sm}{\setminus}

\newcommand{\sq}{\subseteq}
\newcommand{\ov}{\overline}

\definecolor{darkred}{rgb}{0.55, 0.0, 0.0}

\def\G{\Gamma}

\title{Polygonal Faber-Krahn inequality: Local minimality via validated computing.}
\author{Beniamin Bogosel, Dorin Bucur}

\address [Beniamin Bogosel]{CMAP, CNRS, \'Ecole polytechnique, 
	Institut Polytechnique de Paris, 91120 Palaiseau, France}
\email {beniamin.bogosel@polytechnique.edu}

\address[Dorin Bucur]{
	Laboratoire de Math\'ematiques UMR CNRS 5127 \\
	Universit\'e de Savoie,  Campus Scientifique \\
	73376 Le-Bourget-Du-Lac (France)
}
\email{dorin.bucur@univ-savoie.fr}

\begin{document}
	
	\maketitle
	
	\begin{abstract}
		The main result of the paper shows that the regular $n$-gon is a local minimizer for the first Dirichlet-Laplace eigenvalue among $n$-gons having fixed area for $n \in \{5,6\}$. The eigenvalue is seen as a function of the coordinates of the vertices in $\R^{2n}$. Relying on fine regularity results of the first eigenfunction in a convex polygon, an explicit \emph{a priori} estimate is given for the eigenvalues of the Hessian matrix associated to the discrete problem, whose coefficients involve the solutions of some  Poisson equations with singular right hand sides. The \emph{a priori} estimates, in conjunction with certified finite element approximations of these singular PDEs imply the local minimality for $n \in \{5,6\}$. All computations, including the finite element computations, are realized using interval arithmetic.
	\end{abstract}

	\section{Introduction}
	In the previous work \cite{Bogosel_Bucur_Polya} it is shown that a local minimality result for a classical discrete shape optimization problem due to Poly\`a and Szeg\"o, still unsolved today, depends on the positivity of a number of constants. Theory alone is not able to answer this question, for now. Using explicit finite element estimates for solutions of some singular partial differential equations together with interval arithmetic in Intlab \cite{intlab}, controlling all roundoff errors in the computations, allows to establish the local minimality results. 
	
	More precisely, for a bounded open set $\Omega \subset \Bbb{R}^2$ let us consider the eigenvalue problem for the Laplace operator with Dirichlet boundary conditions
	\begin{equation}\left\{ \begin{array}{rcll}
	-\Delta u & =& \lambda u & \text{ in }\Omega,\\
	u &= &0 & \text{ on }\partial \Omega.  
	\end{array} \right.
	\label{eq:dir-eigenvalue}
	\end{equation}
	The spectrum consists only on eigenvalues, counted with corresponding multiplicity,
	$$0< \lb_1(\Om) \le \lb_2(\Om) \le \dots \le \lb_k(\Om) \dots \ra +\infty.$$
	Lord Rayleigh conjectured in 1877 that the first eigenvalue is minimal on the disc, among all other planar domains of the same area. Faber and Krahn gave rigorous proofs in 1923 (see \cite{Da11} for a description of the history of the problem and \cite{he17b,henroteigs} for a survey of the topic).  
	
	Drawing a parallel with the isoperimetric inequality, it is natural to assert that the polygonal shape with $n$ sides and fixed area which minimizes the fundamental Dirichlet Laplace eigenvalue is the regular polygon. This is also in accord with the intuition that a "rounder" shape takes longer to cool down from any initial state, keeping in mind that $\lambda_1(\Omega)$ is the dominant rate of the vanishing exponential involved. 
	
	In their book of 1951, P\' olya  and Szeg\"o conjecture a polygonal version of this inequality (see \cite[page 158]{PoSz51}). Denoting by $\Pol{n}$ the family of  simple polygons with $n$ sides in $\Bbb{R}^2$ and for every $n \geq 3$ consider the problem 
	\begin{equation}
	\min_{P \in \mathcal P_n} \lambda_1(P) |P|.
	\label{eq:polya-conj}
	\end{equation}
	
	Formulation \eqref{eq:polya-conj} is one of the multiple equivalent statements of this problem and has the advantage of being scale invariant: $\lambda_1(tP) = \lambda_1(P)/t^2$ for any $t>0$.
	
	\medskip
	\noindent {\bf P\' olya-Szeg\"o Conjecture (1951).} {\it The unique solution to problem \eqref{eq:polya-conj} is the regular polygon with $n$ sides.
	}
	\medskip
	
	The conjecture holds true for $n=3$ and $n=4$: see \cite[Chapter 3]{henroteigs} where the Steiner symmetrization principle is used. Steiner symmetrization does not work for the case $n \geq 5$ since, the number of vertices could possibly increase after symmetrization. A new approach, which applies only to triangles, was proposed by Fragal\`{a} and Velichkov in \cite{FragalaVelichkov19},  establishing that equilateral triangles are the only critical points for the first eigenvalue. Indrei finds in \cite{Indrei2024} a manifold in which the regular polygon is optimal. Numerical simulations suggesting that the Poly\`a-Szeg\"o conjecture is true are given in \cite{nigam_numerics}, \cite{AFpolygons}, \cite{beniphd}.
	
	In \cite{Bogosel_Bucur_Polya} the authors show that the proof of optimality of the regular $n$-gon for \eqref{eq:polya-conj} can be reduced to a finite number of validated numerical computations. This is achieved by the following results: 
	
	\begin{itemize}[noitemsep]
		\item computation of  second order shape derivatives of simple Dirichlet-Laplacian eigenvalues on Lipschitz domains (including polygons)
		\item computation of the Hessian matrix of $\lambda_1(P)$ in terms of the coordinates of the vertices of $P$. Taking into account the symmetry, the eigenvalues of the Hessian matrix for the regular $n$-gon are characterized using solutions of three PDEs.
		\item Explicit a priori estimates for piece-wise affine finite elements approximation for these three PDEs are given (eigenvalue equation and material derivatives for the eigenfunction with respect to one vertex). Numerical computations justify, {\it up to machine errors}, that for $5\le n \leq 8$ the regular $n$-gon is a local minimizer in \eqref{eq:polya-conj}. 
		\item Qualitative results are given concerning the stability of the Hessian eigenvalues under vertex perturbations and bounds on geometric quantities of the optimal polygon are given.
	\end{itemize}
    In \cite{Clutterbuck2024} the local maximality of the Robin eigenvalue problem for the Laplacian is established for quadrilaterals, using similar lines and exploiting the fact that in the case of the square, the PDEs involved in the Hessian computation become more explicit. 
    
    In order to conclude with local minimality, the missing step in  \cite{Bogosel_Bucur_Polya} was precisely the certification of the machine errors, certification which could not be performed because of the complexity and size of the required computations.
    
    Proof strategies based on certified computations have already been used successfully in articles related to spectral theory. In the paper \cite{schiffer}, Schiffer's conjecture is proved for the regular pentagon: a Neumann eigenfunction that is positive on the boundary and not identically constant is found numerically through certified computations. In \cite{dahne-payne} a counter example to a famous conjecture of Payne is given through validated numerics: an example of two dimensional domain with $6$ holes is given for which the nodal line of the second Dirichlet eigenfunction is closed and does not touch the boundary of the domain. In \cite{three_eigs_triangle} a conjecture stated in \cite{AFpolygons} is settled, proving that there exist triangles which are not isometric and for which the first, second and fourth eigenvalues coincide.
    
   The numerical framework chosen for performing the certified computations in this article uses finite elements. Let us briefly justify this choice. There exist precise numerical methods that compute eigenvalues and eigenfunctions for the Laplace operator on smooth domains or on polygons. For example \cite{dahne-salvy} uses particular solutions coupled with interval arithmetic to obtain guaranteed enclosures for the Laplace eigenvalues which is then used in \cite{dahne-payne} to answer a question of Payne regarding the nodal line of the first eigenfunction on a domain with holes. In \cite{barnett-hassell} the authors give tight inclusions for the Neumann eigenvalues and in \cite{trefethen} a fast numerical method is used for solving Laplace problems. While the results in \cite{trefethen} are promising in singular contexts, there are no explicit quantifiable error bounds that can be applied directly for some of the problems of interest in this work. We also point out \cite{BaNi06} and \cite{BaNi08} for a numerical analysis  of general Poisson equations with distributional data.
    	
   Several non-trivial differences, with respect to the case of the Laplace operator treated in the references above, occur in our case: the operator appearing in our PDEs  is $-\Delta -\lambda_1(\Pp)I$, it acts on a strict subspace contained in  $H^1_0$ and the right-hand side is singular, so that  solutions are not globally in $H^2$.  Here is the first comment: the general approximation results for singular data available in the literature do apply, but the (optimal) approximation order is too small and, in view of the requirements of the \emph{a priori} estimates, the linear systems are too large to be solved  with interval arithmetic software like Intlab. A second comment concerns the use of spectral methods.  
   Because of the singularity of the right hand side, the solutions we deal with have discontinuous normal derivatives across segments inside the regular polygon $\Pp$. Adapting the spectral methods enumerated above to this case is not straightforward. Piecewise affine finite element methods, on the other hand, can capture singular behavior if the singularities themselves are captured exactly by the underlying mesh. One main contribution of the paper is to quantify precisely the corresponding error estimates in such singular contexts.
    
    Numerical computations performed in \cite{Bogosel_Bucur_Polya} correspond to meshes having triangle sizes small enough such that the \emph{a priori} error estimates guarantee that all the non-zero Hessian eigenvalues are positive. In order to turn a numerical computation in a mathematical proof, all computations need to be certified, taking into account even the errors coming from floating point computations. Many such proof techniques emerged recently, using specialized software like Intlab \cite{intlab}, which replaces floating point variables with intervals. Operations performed with intervals guarantee that the output interval associated to a particular mathematical operation contains the desired result. All errors are taken into account, including roundoff errors present in floating point arithmetic.
    
 Certified computation of eigenvalues problems in the $\bo P_1$ finite elements context is described in \cite{LiOi13}. The same paper allows to handle the context corresponding to Homogeneous Dirichlet and Neumann boundary conditions. For non-homogeneous Neumann boundary conditions the paper \cite{xuefeng-neumann} applies. In Section \ref{sec:material} an alternative formulation is proposed for a particular problem involving non-homogeneous Neumann boundary conditions, allowing the validation of transmission problems, in particular. This framework should extend to other situations where the boundary data implies the $H^2$ regularity of solutions, following classical results in \cite{grisvard}.
    
	The results of this paper finalize the proof of the local minimality of the regular $n$-gon for problem \eqref{eq:polya-conj} for $n \in \{5,6\}$ by performing all computations (including the finite element computations) using interval arithmetic in Intlab \cite{intlab}. To achieve this goal, in order to make the computation effort acessible, significant improvements over the \emph{a priori} error estimates given in \cite{Bogosel_Bucur_Polya} are done. They rely on the following theoretical issues:	
	\begin{itemize}
		\item The regularity of the material derivatives of the first eigenfunction of the Dirichlet-Laplace operator on the regular $n$-gon (see the next section), is  characterized (and improved) locally, on a suitable partition of the polygon. The piece-wise $H^2$ regularity is established and explicit error estimates for the piecewise affine finite elements are given. Because of the localization procedure which isolates the singularities, the order of the estimates and the constants are significantly improved compared with \cite{Bogosel_Bucur_Polya}.
		\item The new error estimate allows meshes which have significantly fewer elements to be considered for approximating the eigenvalues \eqref{eq:dir-eigenvalue} and the material derivatives. The size of the linear generalized eigenvalue problem and linear systems involved becomes small enough such that Intlab \cite{intlab} can be used to certify the corresponding computations.
	\end{itemize} 
	
\smallskip	
	\noindent \bo{Structure of the paper.} Section \ref{sec:prelim} introduces preliminary results from \cite{Bogosel_Bucur_Polya} regarding shape derivatives, polygonal perturbations and the Hessian matrix of the first eigenvalue on the regular $n$-gon. Explicit \emph{a priori} estimates regarding finite element approximations of the eigenvalues are recalled. 
	
	Section \ref{sec:material} deals with explicit estimates for solutions to PDEs involving the Laplace operator with a singular source term, namely a measure with density in $H^{1/2}$ supported on a segment,  which gives rise to transmission conditions. These PDEs appear in the expression of the eigenvalues of the Hessian matrix described above. If the singularities are meshed exactly, the optimal order of convergence $O(h)$ can still be reached in the approximation process. The main tool is proving that the solution is piece-wise $H^2$, adapting results from \cite[Chapter 4]{grisvard}. The method proposed can extend to other cases where solutions are piecewise $H^2$ and the the mesh captures the singularities exactly.
	
	Section \ref{sec:num-validation} gives various results regarding error estimates related to linear algebra, which justify the choices made in the validation process. Section \ref{sec:validation} presents the validation strategy and the corresponding results concerning the local minimality of the regular $n$-gon in the Poly\`a-Szeg\"o conjecture for $n \in \{5,6\}$. The validation strategy illustrates the potential sources of errors: error between continuous problems and exact finite element solutions, errors in the discrete eigenvalue problems and linear systems, roundoff errors appearing when working in floating point precision. 
	
	For the sake of completeness, developments related to particular aspects in the code, possibly relevant for other applications, are detailed in the Appendices. Appendix \ref{app:morley} describes a possible implementation for the assembly of Morley finite elements which helps computing tight certified upper bounds for the interpolation constant for $\bf P_1$ finite elements, based on the method proposed in \cite{kobayashi}. Appendix \ref{app:code} describes the main functions in the code allowing to verify the local minimality. The corresponding code is provided at the repository:
	
	\begin{center}
		\href{https://github.com/bbogo/PolyaSzego}{\nolinkurl{https://github.com/bbogo/PolyaSzego}}
	\end{center}
	allowing to verify and reproduce the results proved in this paper. 
	Appendix \ref{app:intlab} details a slight modification in the function \texttt{verifyeig} from Intlab (related to certification of discrete eigenvalue problems), replacing a matrix inversion with validated linear systems. For sparse matrices of large size this improves significantly the performance of the process of certifying solutions to generalized eigenvalue problems.
	
	\section{Preliminaries}
	\label{sec:prelim}
	
	\subsection{Dirichlet-Laplacian eigenvalues on polygons: first and second derivatives} 
	Shape derivatives are introduced to study the behavior of shape functionals with respect to boundary variations. Given $\xi \in W^{1,\infty}(\Bbb{R}^2,\Bbb{R}^2)$  Lipschitz vector fields the  shape derivative verifies
	\[ J((I+\xi)(\Omega)) = J(\Omega)+J'(\Omega)(\xi)+o(\|\xi\|_{W^{1,\infty}}),\]
	where $\xi \mapsto J'(\Omega)(\xi)$ is a linear form in $\xi$. 
	
	Shape derivatives for a simple Dirichlet-Laplace eigenvalue are well known in the literature when the domain $\Omega$ is smooth (\cite{henrot-pierre-english}, \cite{Sokolowski1992}, \cite{Delfour2011}). For a simple eigenvalue one has 
	\begin{equation}\label{eq:sh-deriv-boundary}
	 \lambda'(\Omega)(\xi) =- \int_{\partial \Omega} (\partial_n u)^2 \xi \cdot \bo n, 
	 \end{equation}
	where $\partial_n$ denotes the normal derivative, $\bo n$ is the outer normal to the domain $\Omega$ and $u$ is an $L^2$ normalized eigenfunction associated to $\lambda$. The previous formula extends to all cases where $u \in H^2(\Omega)$. It can be noted that in this case the shape derivative depends only on $\xi\cdot \bo n$ on the boundary $\partial \Omega$. This structure is verified in most cases where the functional and the domain are smooth enough \cite[Chapter 5]{henrot-pierre-english}.
	
	In situations where the domain is less smooth, known results do not apply directly. In \cite{Bogosel_Bucur_Polya}, based on ideas developed in \cite{laurain2ndDeriv}, the following volume form of the shape derivative is derived
	\begin{equation}\label{eq:sh-deriv-vol}
	 \lambda'(\Omega)(\xi) = \int_\Omega \bo S_1 : D\xi,
	\end{equation}
	where 
	\begin{equation}\label{eq:S1}
	\bo S_1 = (|\nabla u|^2-\lambda u^2)\Id -2(\nabla u \otimes \nabla u).
	\end{equation}
	 The matrix scalar product is denoted by $:$, the identity matrix by $\Id$ and $D\xi$ is the Jacobian of $\xi$. This formula applies whenever $u \in H^1(\Omega)$, thus including less regular contexts. By analyzing singularities of the eigenfunction $u$ on polygons in \cite{Bogosel_Bucur_Polya} it is shown that \eqref{eq:sh-deriv-boundary} remains valid on polygons. 
	
	An advantage of the volumic formulation \eqref{eq:sh-deriv-vol}  is that it can be further differentiated without additional regularity assumptions. Thus, second order shape derivatives for simple Dirichlet-Laplace eigenvalues are computed in \cite{Bogosel_Bucur_Polya} for all domains where $u \in H^1(\Omega)$ and for all Lipschitz perturbation fields. The polygonal framework is thus included. For brevity, the formula for the second shape derivative is ommitted here, since all details can be found in \cite{Bogosel_Bucur_Polya}. In the following, the polygonal case which is of interest in this work is detailed.
	
\begin{figure}
		\centering 
		\includegraphics[height=0.3\textwidth]{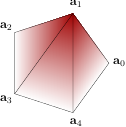}\quad
		\includegraphics[height=0.3\textwidth]{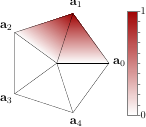}\quad
				\includegraphics[height=0.3\textwidth]{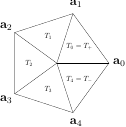}
		\caption{Examples of admissible triangulations used for defining perturbations on a polygon and graphical view of the function $\varphi_1$.}
		\label{fig:simple-triangulations}
	\end{figure}

	{\bf Polygons and polygonal perturbations.} Suppose $\Omega$ is a $n$-gon and follow the strategy described in \cite{Bogosel_Bucur_Polya}, \cite{laurain2ndDeriv}. Starting from a perturbation of the vertices, the perturbation field $\xi \in W^{1,\infty}(\R^2)$ will be built as follows. Denote the vertices of the polygon by $\bo a_i \in \Bbb{R}^2,\ i=0,...,n-1$ and for each vertex consider the vector perturbation $\xi_i \in \Bbb{R}^2,\ i=0,...,n-1$. Whenever necessary, we suppose that the indices are considered modulo $n$. Consider a triangulation $\mathcal T$ of $\Omega$ such that the edges of the polygon are complete edges of some triangles in this triangulation.  Moreover, consider the following globally Lipschitz functions  $\varphi_i$ for $0\leq i\leq n-1$ that are piecewise affine on each triangle of $\mathcal T$ and satisfy
	\begin{equation}\label{eq:phi}
	\varphi_i(\bo a_j) = \delta_{ij} = \begin{cases}
	1 & \text{ if } i=j \\
	0 & \text{ if } i\neq j
	\end{cases}.
	\end{equation}
	Several choices are possible, as the two examples  of Figure \ref{fig:simple-triangulations} show, their extension outside the polygon being irrelevant. 
	Then, we build a global perturbation in $\R^2$ given by
	\begin{equation}\label{eq:phi2}\xi = \sum_{i=0}^{n-1} \zeta_i \xi_i \in W^{1,\infty}(\R^2).
	\end{equation}
	
	Plugging the polygonal perturbation field $\xi$ given in \eqref{eq:phi2} into the shape derivative \eqref{eq:sh-deriv-vol} shows that the gradient of a simple eigenvalue with respect to the vertex coordinates $\bo x = (\bo a_0, \bo a_1, ... \bo a_{n-1})$ is given by
	\[ \nabla \lambda(\bo x) = \left( \int_\Omega \bo S_1 \nabla \varphi_i\right)_{i=0}^{n-1},\]
	where $\varphi_i$ are the piecewise affine functions on the triangulation $\mathcal T$ verifying \eqref{eq:phi} and $\bo S_1$ is the matrix defined in \eqref{eq:S1}. For simplicity, $\varphi_0$ which will appear often in the following is simply denoted by $\varphi$. 
	
	The area of a polygon is explicit in terms of vertex coordinates and, assuming the vertices are ordered in the anti-clockwise direction, it is given by:
	\[ \mathcal A(\bo x) = \frac{1}{2} \sum_{i=0}^{n-1} (x_iy_{i+1}-x_{i+1}y_i).\]
	
	Consider the case of the regular $n$-gon $\Pp$. In the following we use the notation $\theta = 2\pi/n$. Assume that $\Pp$ is inscribed in the unit disk with $\bo a_0 = (1,0)$. Consider the symmetric triangulation $\mathcal T$ illustrated in Figure \ref{fig:simple-triangulations} with an extra vertex at the center of $\Pp$. Denote by $T_i$ the triangles in the triangulation $\mathcal T$ like in center image from Figure \ref{fig:simple-triangulations}. The triangles $T_0$ and $T_{n-1}$ will sometimes be denoted by $T_+, T_-$, respectively and correspond to the support of the $\bo P_1$ function $\varphi$. 
	
	The gradient of $\lambda(\bo x)$ with respect to perturbations in the first vertex $\bo a_0$ of the regular polygon $\Pp$ is 
	\[ \int_\Pp \bo S_1\nabla \varphi_0 = \begin{pmatrix}
	-2\lambda_1/n \\ 0
	\end{pmatrix}.\]
	
	Plugging polygonal perturbation fields into the second shape derivative formula obtained in \cite{Bogosel_Bucur_Polya} gives an explicit formula for the Hessian matrix $D^2 \lambda_1(\bo x)$ for $\lambda_1(\bo x)$ on the regular $n$-gon. In \cite[Section 4]{Bogosel_Bucur_Polya} explicit formulas for the eigenvalues of the Hessian are given for the regular $n$-gon.
	
	Denote by $a(u,v)$ the bilinear form
	\[ a(u,v) = \int_\Pp \nabla u \cdot \nabla v-\lambda_1(\Pp) uv,\]
	which is positive semi-definite on $H_0^1(\Pp)$. The only functions which cancel $a(\cdot,\cdot)$ are multiples of the first eigenfunction $u_1$.
	
	First, the material derivatives corresponding to perturbations of the first vertex $\bo a_0$ of the regular $n$-gon $\Pp$ are defined by $\bo U = (U^1,U^2) \in H_0^1(\Pp)^2$ verifying: 
	\begin{align}
	 &a(U^1,v) =  \int_\Pp (\nabla u\cdot \nabla \varphi)\partial_x v+\int_\Pp (\nabla v \cdot \nabla \varphi) \partial_x u-\frac{2\lambda_1}{n}\int_\Pp uv, \forall v \in H_0^1(\Pp)\notag \\
	 &a(U^2,v) =  \int_\Pp (\nabla u\cdot \nabla \varphi)\partial_y v+\int_\Pp (\nabla v \cdot \nabla \varphi) \partial_y u, \forall v \in H_0^1(\Pp) \label{eq:material-U0}
	\end{align}
	Since the equations for $U^1,U^2$ only characterize solutions up to the addition of a multiple of $u_1$, the orthogonality conditions
	\begin{equation}\label{eq:orthogonalityU}
	\int_\Pp U^i u_1 = 0
	\end{equation}
	are considered, to fix unique solutions in \eqref{eq:material-U0}. In \cite{Bogosel_Bucur_Polya} it is shown that the Hessian matrix $D^2\lambda_1(\bo x)$ does not depend on the chosen normalization. Consider the following notations from \cite[Prop. 4.10]{Bogosel_Bucur_Polya} for $0 \leq k \leq n-1$, recalling that $\theta=2\pi/n$:
	
	\begin{align*} 
			A_k
			&= \sum_{j=0}^{n-1} (\cos (j+1)k\theta+\cos jk\theta) \int_{T_j}\nabla u_1 \cdot \nabla U^1\\ 
			&+ \sum_{j=0}^{n-1} \frac{\cos(j+1)k\theta-\cos jk\theta}{\sin \theta} \int_{T_j}\begin{pmatrix}
			-\sin (2j+1)\theta & \cos (2j+1) \theta \\
			\cos(2j+1)\theta & \sin (2j+1)\theta
			\end{pmatrix}\nabla u_1 \cdot \nabla U^1 \\
			B_k &
			=\frac{\cos \theta}{\sin \theta}\sum_{j=0}^{n-1} (\cos(j+1)k\theta-\cos jk\theta) \int_{T_j} \nabla u_1 \cdot \nabla U^2\\
			&+ \sum_{j=0}^{n-1} \frac{\cos(j+1)k\theta-\cos jk\theta}{\sin \theta} \int_{T_j}
			\begin{pmatrix}
			-\cos(2j+1)\theta & -\sin(2j+1)\theta \\ 
			-\sin(2j+1)\theta & \cos(2j+1)\theta 
			\end{pmatrix}\nabla u_1 \cdot \nabla U^2\\
			C_k & 
			= \frac{\cos \theta}{\sin \theta}\sum_{j=0}^{n-1} (\sin(j+1)k\theta-\sin jk\theta)  \int_{T_j} \nabla u_1\cdot  \nabla U^1\\
			&+ \sum_{j=0}^{n-1}\dfrac{\sin(j+1)k\theta-\sin jk\theta}{\sin \theta} \int _{T_j} \begin{pmatrix}
			-\cos (2j+1)\theta & -\sin (2j+1) \theta \\
			-\sin (2j+1)\theta & \cos (2j+1)\theta
			\end{pmatrix}\nabla u_1 \cdot \nabla U^1\\	
			D_k &  = \sum_{j=0}^{n-1}(\sin (j+1)k\theta  + \sin jk\theta) \int_{T_j} \nabla u_1 \cdot \nabla U^2 \\
			&+ \sum_{j=0}^{n-1}\dfrac{\sin(j+1)k\theta-\sin jk\theta}{\sin \theta} \int _{T_j} \begin{pmatrix}
			-\sin (2j+1)\theta & \cos (2j+1) \theta \\
			\cos (2j+1)\theta & \sin (2j+1)\theta
			\end{pmatrix}\nabla u_1 \cdot \nabla U^2 
	\end{align*}
	
	Then, according to \cite[Section 4]{Bogosel_Bucur_Polya}, the eigenvalues of the Hessian matrix $D^2(\lambda_1(\bo x)\mathcal A(\bo x))$ associated to the scale invariant formulation \eqref{eq:polya-conj} are given by the following result.  
	
	\begin{thm}\label{thm:eig-hessian}
		For $0\leq k \leq n-1$, $\theta=2\pi/n$ denote
		\begin{align*}
		\alpha_k & = \frac{2n(1-\cos(k\theta))}{\sin \theta} \int_{T_0}(\partial_x u_1)^2-2|\Pp|A_k \\
		\beta_k & = \frac{2n(1-\cos(k\theta))}{\sin \theta}\int_{T_0} (\partial_y u_1)^2 -2|\Pp| B_k\\
		\gamma_k & = -2|\Pp|C_k  = 2|\Pp|D_k
		\end{align*}
		and
		\[  \mu_{2k} = 0.5(\alpha_k+\beta_k- \sqrt{(\alpha_k-\beta_k)^2+4\gamma_k^2}), \ \mu_{2k+1} = 0.5(\alpha_k+\beta_k+ \sqrt{(\alpha_k-\beta_k)^2+4\gamma_k^2}).\]
		The eigenvalues of the Hessian matrix of $\lambda_1(\bo x)\mathcal A(\bo x)$ given in are given by $\mu_j$, $j=0,...,2n-1$. 
		\label{thm:explicit-eigenvalues}
	\end{thm}

	Note that $k=0$ trivially gives two zero eigenvalues in the formula above. Moreover, in \cite[Proposition 4.12]{Bogosel_Bucur_Polya} it is shown that the Hessian of $\lambda(\bo x)\mathcal A(\bo x)$ has two additional zero eigenvalues. This is in agreement with the fact that $\lambda(\Omega)|\Omega|$ is scale invariant: translations, rotations and homotheties do not change the objective function. Proving that the regular $n$-gon is a local minimum amounts to the following (see \cite[Section 4]{Bogosel_Bucur_Polya}).  
	
	\medskip
	
	\noindent {\bf Objective.} Given $n\geq 5$ prove that $2n-4$ of the eigenvalues of the Hessian matrix of $\bo x \mapsto \lambda_1(\bo x)\mathcal A(\bo x)$ are strictly positive. 
	
	\medskip
	
	Theorem \ref{thm:explicit-eigenvalues} shows that eigenvalues of the Hessian matrix are explicit in terms of the following quantities:
	\[ \int_{T_0} (\partial_x u_1)^2, \int_{T_0} (\partial_y u_1)^2,\]
	\[ \int_{T_j}\partial_x u_1 \partial_x U^1, \int_{T_j}\partial_x u_1 \partial_y U^1, \int_{T_j}\partial_y u_1 \partial_x U^1, \int_{T_j}\partial_y u_1 \partial_y U^1, 0 \leq j \leq n-1, \]
	\[ \int_{T_j}\partial_x u_1 \partial_x U^2, \int_{T_j}\partial_x u_1 \partial_y U^2, \int_{T_j}\partial_y u_1 \partial_x U^2, \int_{T_j}\partial_y u_1 \partial_y U^2, 0 \leq j \leq n-1.\]
	Since purely theoretical results allowing to prove the positivity of eigenvalues in Theorem \ref{thm:explicit-eigenvalues} are not available, numerical computations are used in the proof for $n \in \{5,6\}$. Explicitly quantified error estimates for fintie element computations together with interval arithmetic for controlling floating point errors are used. 
	
	\subsection{General approximation framework} 
	
	{\bf a) Interpolation errors.} Consider $\Omega$ a polygon meshed exactly with a triangulation $\mathcal T_h$. Consider the $\bf P_1$ finite element space on $\mathcal T_h$ consisting of piece-wise affine functions on triangles in $\mathcal T_h$. The norms $L^2, H^1, H^2$ are considered on the full domain, unless specified otherwise.
	
	In the following denote $\Pi_{1,h}$ the interpolation operator associated to the $\bo P_1$ Lagrange finite element. More precisely, $\Pi_{1,h}w$ is the piecewise affine function in $\mathcal V_h$ such that $w$ and $\Pi_{1,h}w$ have the same values at the nodes of the mesh. Results of \cite{LiOi13} show that whenever $w\in H^2(\Omega)$ we have
	\begin{equation}\label{eq:interpolation} \|\nabla w-\nabla \Pi_{1,h} (w)\|_{L^2} \leq C_1h \|D^2w\|_{L^2}.
	\end{equation}
	The above inequality becomes explicit when an upper bound for $\|D^2w\|_{L^2}$ is known. The constant $C_1$ is explicit, depending on the triangles of the mesh $\mathcal T_h$.  We shortly recall of the results  of  \cite[Theorem 4.3]{LiOi13}. 
	
	In each triangle $T_i\in {\mathcal T}_h$, the ratio between the smallest edge and the middle one $L_i$ is denoted $\alpha_i$ and the angle between these two edges is $\tau_i$. Then, we denote
	$$C(T_i):= 0.493  L_i\frac{1+ \alpha_i^2+\sqrt{1+2 \alpha_i^2 \cos(2\tau_i) + \alpha_i^4}}{\sqrt{2\big (1+ \alpha_i^2-\sqrt{1+2 \alpha_i^2 \cos(2\tau_i) + \alpha_i^4}\Big)}}.$$
	Following \cite[Section 2]{LiOi13}, we introduce the constant
	{\begin{equation}\label{eq:C1}
		C_1= \sup_h \frac{C(T_i)}{h},
		\end{equation}}
	where the parameter $h$ dictating the size of the mesh is the size of the median edge. In the applications concerning regular polygons it is useful to consider $\mathcal T_h$ consisting of congruent triangles each one similar to $T_0$ in Figure \ref{fig:simple-triangulations}. This renders the constant \eqref{eq:C1} completely explicit in terms of $\theta$ and $h$.
	
	Since the mesh used of the regular polygon used in the computations is constructed using congruent triangles with angles $(2\pi/n, \pi/2-\pi/n, \pi/2-\pi/n)$ more precise estimates for the interpolation constant can be found through direct estimation. We use the technique described in \cite{kobayashi}, which shows that discrete problems related to the Morley finite element gives explicit bounds for the constant $C_1$ in \eqref{eq:interpolation}. A description of the procedure for obtaining certified upper bounds for constant $C_1$ defined in \eqref{eq:interpolation}, giving better bounds than \eqref{eq:C1} for triangles of interest in our case, is described in Appendix \ref{app:morley}. The constants used in the validation computations are shown in Table \ref{tab:interp-const}.
	
\medskip	
	\noindent {\bf b) General Poisson problems.} 	
	For $\Omega$ a polygon in the plane, consider $\Gamma \subset \partial \Omega$ a subset of its boundary and $\Gamma_N = \partial \Omega \setminus \Gamma$. Consider the problem 
	\[ -\Delta w = f,\ \  \partial_n w = 0 \text{ on } \Gamma_N, w=0 \text{ on } \Omega \setminus \Gamma_N\]
    where  $f \in H^{-1}(\Omega)$. The associated variational formulation reads: $w \in H_\Gamma(\Omega): = \{ v \in H^1(\Omega) : v=0\text{ on } \Gamma \subset \partial \Omega\}$,
	\begin{equation}\label{eq:cont_f}  \int_\Omega \nabla w \cdot \nabla v = (f,v)_{H^{-1},H^1} , \text{ for every } v \in H_\Gamma(\Omega).
	\end{equation}
	Assume that the solution is $H^2$ regular, for example if $\Omega$ is convex.
	
	Suppose $\Omega$ is polygonal and triangulated, obtaining a mesh $\mathcal T_h$ and an associated space of piecewise affine functions $v_h \in \mathcal V_h$ such that $v_h$ vanishes on nodes in $\Gamma$. Assume $\mathcal V_h \subset H_\Gamma(\Omega)$ and consider $p_h$ the solution to the problem
	\begin{equation}\label{eq:disc_f} p_h \in \mathcal V_h :  \int_\Omega \nabla p_h \cdot \nabla v_h = (f,v_h)_{H^{-1},H^1} , \text{ for every } v_h \in \mathcal V_h.
	\end{equation}
	Then $p_h$ is the projection of $u$, solution of \eqref{eq:cont_f}, on $\mathcal V_h$. In particular, for every $v_h \in \mathcal V_h$ we have
	\[ \|\nabla w-\nabla v_h\|_{L^2}^2 = \|\nabla w-\nabla p_h\|_{L^2} ^2+\|\nabla p_h-\nabla v_h\|_{L^2} ^2\geq \|\nabla w-\nabla p_h\|_{L^2}^2 \]
	since $\int_\Omega \nabla (w-p_h)\cdot \nabla(p_h-v_h)=0$.

	Taking $v_h = \Pi_{1,h}(w)$ shows that if $w \in H^2(\Omega)$ then 
	\[ \|\nabla w-\nabla p_h\|\leq C_1h \|D^2w \|_{L^2}.\]
	
	If one discretizes also the distribution $f$, i.e. consider the problem 
	 \begin{equation}\label{eq:disc_fh} w_h \in \mathcal V_h :  \int_\Omega \nabla w_h \cdot \nabla v_h = (f_h,v_h)_{H^{-1},H^1} , \text{ for every } v_h \in \mathcal V_h.
	 \end{equation}
	 then $w_h$ and $v_h$ verify the same equation with different right hand sides. It is immediate to observe that
	 \[ \|\nabla w_h-\nabla v_h\|_{L^2}\leq \|f-f_h\|_{H^{-1}}.\]
	 
	 Thus, we obtain
	 \[ \|\nabla w - \nabla w_h\|_{L^2} \leq \|f-f_h\|_{H^{-1}}+C_1h \|D^2 w\|_{L^2}.\]
	
	\medskip
	\noindent
	{\bf c) Dirichlet Laplacian eigenvalues and eigenfunctions.} Assume the Dirhchlet boundary conditions are applied on the whole boundary $\Gamma= \partial \Omega$ in the following.  
	Given a triangulation $\mathcal T_h$ of the polygon $\Omega$ and $\mathcal V_h \subset H_0^1(\Omega)$ the $\bo P_1$ finite elements space on $\mathcal T_h$, denote by $\lb_{k,h}, u_{k,h}$ the $k$-th eigenvalue of $\Omega$ and its associated eigenfunction approximated in $\mathcal V_h$, solving
	\begin{equation}
	u_{k,h} \in \mathcal V_h, \int_\Omega \nabla u_{k,h}\cdot \nabla v_h = \lambda_{k,h}\int_\Omega u_{k,h} v_h, \ \ \ \forall v_h \in \mathcal V_h.
	\label{eq:eigenvalue-finite-elements}
	\end{equation}
	This problem is a generalized eigenvalue problem.
	
	Results of \cite[Theorem 4.3]{LiOi13} show that for $\Omega = \Pp$
	$$ \forall k \ge 1, \;\; \lb_{k,h} \geq \lb_k \geq \frac{\lb_{k,h}}{1+ C_1^2 h^2 \lb_{k,h}^2},$$
	since the associated eigenfunctions are in $H^2(\Omega)$. Therefore the following completely explicit estimate holds
	\begin{equation}
	|\lambda_k-\lambda_{k,h}| \leq \lambda_{k,h}^3 C_1^2/(1+C_1^2h^2\lambda_{k,h}^2)\ h^2,
	\label{eq:difflam}
	\end{equation}
	with $C_1$ verifying \eqref{eq:interpolation}.
	
	For bounds concerning the approximation of the eigenfunction $u_1$ associated to the simple eigenvalue $\lambda_1$ on the regular polygon $\Pp$, the results in \cite[Section 5]{Bogosel_Bucur_Polya} are used. The key ingredients are the following:
	\begin{itemize}
		\item Projection of $u_1$ on $\mathcal V_h$, defining 
		\[ p_h \in \mathcal V_h, \int_\Pp \nabla p_h \cdot \nabla v = \lambda_1\int_\Pp u_1 v, \forall v \in \mathcal V_h.\]
		Like in the estimate for the general Laplace problem it follows that $\|\nabla u_1-\nabla p_h\|_{L^2} \leq C_1h\|D^2 u_1\|_{L^2(\Pp)} = C_1\lambda_1 h$, since $\|D^2 u_1\|_{L^2} = \|\Delta u_1 \|_{L^2} = \lambda_1$ \cite[Chapter 4]{grisvard}. The convexity of $\Pp$ and the Aubin-Nitsche lemma implies also that $\|u_1-p_h\|\leq C_1^2 \lambda_1 h^2$. 
		\item Project $p_h$ on the orthogonal of $u_{1,h}$ in $\mathcal V_h$: $p_h = \alpha u_{1,h} +\overline p_h$ with $\int_\Omega \overline p_h u_{1,h} =0$. Without loss of generality assume $\alpha>0$. This shows that
		\[ \|\nabla p_h-\nabla u_{1,h}\|_{L^2}=\|(\alpha-1)\nabla u_h  +\nabla \overline p_h\|_{L^2}\leq |\alpha-1|\lambda_{1,h}^{1/2}+\|\nabla \overline p_h\|_{L^2}.  \]
		It follows that
		\[ \int_\Pp (\nabla \overline p_h \cdot \nabla v_h -\lambda_{1,h}\overline p_h v_h) = \int_\Pp (\lambda_1 u_1-\lambda_{1,h}p)v_h,\ \forall v_h \in \mathcal V_h.\]
		Poincar\'e's inequality in the orthogonal of $u_{1,h}$ in $\mathcal V_h$ gives
		\begin{equation}
		\lambda_{2,h}^{1/2}\|\overline p_h\|_{L^2} \leq \|\nabla \overline p_h \|_{L^2} \le \frac{\lb_{2, h}^\frac 12}{(\lb_{2, h}-\lb_{1, h})}\Big (|\lb_1-\lb_{1, h}|+ \lb_{1, h}\|u_1- p_h\|_{L^2}\Big).
		\label{eq:estimate-p-bar}
		\end{equation}
		Moreover, since $\|p_h\|_{L^2}^2= \alpha^2 + \|\overline p_h\|_{L^2}^2$ it follows that 
		\[ |1-\alpha| \leq |1-\alpha^2| \leq \int_\Pp \overline p_h^2 +\int_\Pp (u_1^2-p_h^2) \leq \|\overline p_h\|_{L^2}^2 + \|u_1-p_h\|_{L^2}(2+\|u_1-p_h\|_{L^2}).\]
	\end{itemize}
    Combining these estimates leads to
    \begin{equation}\label{eq:grad-u-estimate}
    \|\nabla  u_1-\nabla u_{1,h}\|_{L^2} \leq C_1\lambda_1h  + O(h^2),
    \end{equation}
    where the term $O(h^2)$ is explicit using the inequalities proved above.
    A similar estimate using $\|u_1-u_h\|_{L^2}\leq \|u_1-p_h\|_{L^2}+\|p_h-u_{1,h}\|_{L^2}$ 
    gives
    \begin{equation}\label{eq:L2-u-estimate}
    \|  u_1-u_{1,h}\|_{L^2} \leq C_1^2\lambda_1h^2   + O(h^2),
    \end{equation}
    where all terms are explicit. 
    
 The three estimates $|\lambda_1-\lambda_{1,h}|$, $\|\nabla u_1-\nabla u_{1,h}\|_{L^2}$, $\|u-u_h\|_{L^2}$ are linked through the following relation
    \begin{align}
    \|\nabla u_1-\nabla u_{1,h}\|_{L^2}^2-\lambda_1\|u_1-u_{1,h}\|_{L^2}^2 & = \int_\Pp |\nabla u_1|^2-2\int_\Pp \nabla u _1\cdot \nabla u_{1,h} +\int_\Pp |\nabla u_{1,h}|^2 \notag \\
    & -\lambda_1 \int_\Pp u_1^2+2\lambda_1\int_\Pp u_1 u_{1,h}-\lambda_1 \int_\Pp u_{1,h}^2 \notag \\
    & = \lambda_{1,h}-\lambda_1,   \label{eq:relation-errors}
    \end{align}
    where the weak formulations for $u_1, u_{1,h}$ are used. Relation \eqref{eq:relation-errors} is also stated in \cite[Equation (3.14)]{Vohralik2017}. One may observe that the three error estimates are related: knowing two of them gives the third one. Since the error estimates for $\|\nabla u_1-\nabla u_{1,h}\|_{L^2}$ and $\|u_1-u_{1,h}\|_{L^2}$ already contain $|\lambda_{1,h}-\lambda_1|$, the estimation procedure can be iterated, replacing a weaker estimate with a stronger one.
    
    It can be observed that even though the relation \eqref{eq:relation-errors} is stated for the regular polygon $\Pp$, it remains valid for every polygon $\Omega$, when meshed exactly. It is also verified for higher eigenpairs since the only result used is the variational formulation for $u_1, u_{1,h}$.
	
	\section{Material derivative decomposition: piecewise $H^2$ regularity and estimates}
	\label{sec:material}
	
	Similar estimates are needed for solutions of \eqref{eq:material-U0}, but it turns out that $U^1, U^2$ do not belong to $H^2(\Pp)$. Indeed, the equations solved by  $U^1, U^2$ have singular right hand sides which can be written as distributions of the form
	\[ (f^{1,2},v)_{H^{-1},H^1} = \int_\Pp f^{1,2}_\reg v +\sum_{i=0}^{N-1}\int_{S_i} g_iv, \]
	where $f^{1,2}_\reg \in L^2(\Pp)$, $S_i$ represent the rays connecting the origin to the vertices $\bo a_i$ of $\Pp$, $g_i \in H^{1/2}(S_i)$. This motivated the estimates in interpolated spaces from \cite{Bogosel_Bucur_Polya}, based on the fact that the eigenfunction $u_1$ belongs to $ H^{2+s}(\Pp)$ for some $s>0$, small. 
	
	The main observation is that the solutions of \eqref{eq:U0} are piece-wise $H^2$ on every triangle $T_i$, $i=0,...,n-1$ defined in Figure \ref{fig:simple-triangulations} and the $H^2$ norms on these triangles can be controlled. This significantly improves the convergence rate and the constants involved, as shown in the remaining of this section.
	
	Consider the context in \cite[Section 5]{Bogosel_Bucur_Polya}. We propose to give a notably better estimate than the one in Lemma 5.1 of  \cite{Bogosel_Bucur_Polya}, taking into account the particular structure of the problem.
	
	Using the notations from \cite{Bogosel_Bucur_Polya}, we have for $\bo U = \bo U_0$
	\begin{equation} 
	a(\bo U,v) = (\bo f_{\reg},v) + (\bo f_{\sing},v), \text{ for every } v \in H_0^1(\Bbb{P}_n)
	\label{eq:U0}
	\end{equation}
	with $a(u,v) = \int_{\Bbb P_n} \left(\nabla u \cdot \nabla v - \lambda_1 uv\right)$, the equality in \eqref{eq:U0} being understood in a vectorial sense. Integrating by parts in the right hand side of \eqref{eq:U0}, the regular part has the form
	\begin{equation}
	(\bo f_{\reg},v) = \begin{pmatrix}
	-2\int_{T_+\cup T_-} (\nabla \varphi\cdot \nabla (\partial_x u_1))v -\frac{2\lambda_1}{n} \int_{\Bbb P_n} u_1v \\
	-2\int_{T_+\cup T_-} (\nabla \varphi\cdot \nabla (\partial_y u_1))v
	\end{pmatrix}.
	\label{eq:f-reg}
	\end{equation}
	Denoting by $\partial_r u_1$ the derivative of $u_1$ in the radial direction. The singular part is given by
	\begin{equation}
	(\bo f_{\sing},v) = \begin{pmatrix}
	-\int_{S_+}\frac{1}{\tan \theta} \partial_r u_1 v -\int_{S_-}\frac{1}{\tan \theta} \partial_r u_1 v +\int_{S_0}\frac{2}{\tan \theta} \partial_r u_1 v\\
	-\int_{S_+}\partial_r u_1 v + \int_{S_-} \partial_r u_1 v
	\end{pmatrix}, 
	\label{eq:f-sing}
	\end{equation}
	where $\theta = 2\pi/n$. 
	It is straightforward to see that since $u_1 \in H^2(\Bbb{P}_n)$, the distribution $\bo f_{\reg}$ can be identified with an element in $L^2(\Pp)$. Moreover, a straightforward computation using integration by parts shows that $(\bo f_{\reg},u_1) = 0$. 

\begin{figure}
		\begin{tabular}{cc}
		\includegraphics[height=0.27\linewidth]{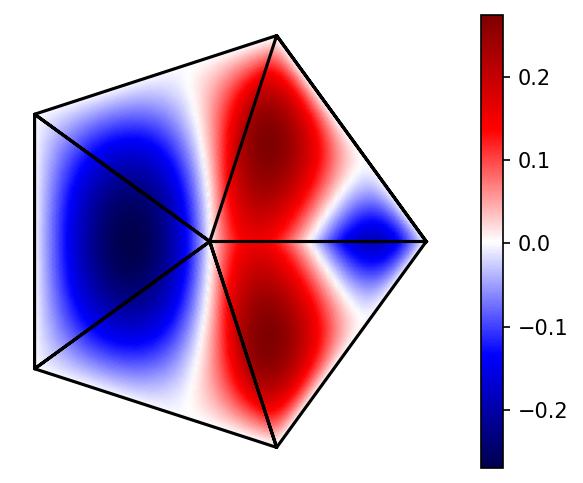} &
		\includegraphics[height=0.27\linewidth]{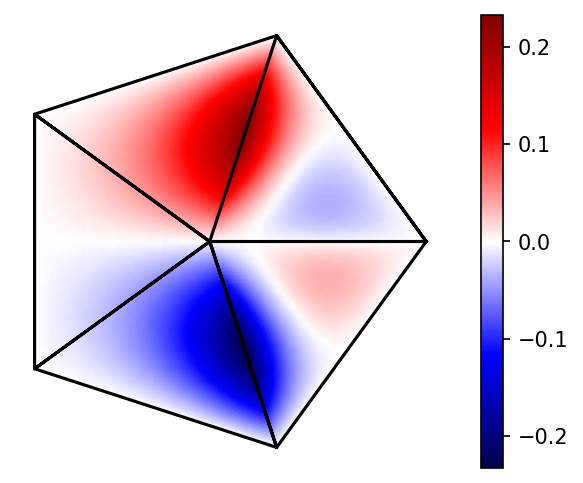}\\
		$U^1$ & $U^2$ 
		\end{tabular}
	\caption{Graphical representation of finite element solutions for problems \eqref{eq:U0}. Discontinuities of normal derivatives are visible across certains rays connecting the center of the regular $n$-gon to the vertices.}
	\label{fig:material}
	\end{figure}

	The regular part of the distribution, being orthogonal on $u_1$, there exists a unique solution for
	\begin{equation}\bo U_{\reg} \in H_0^1(\Pp)^2,\ \ 
	a(\bo U_{\reg},v) = (\bo f_{\reg},v), \forall v \in H_0^1(\Bbb P_n), \int_{\Bbb P_n}\bo U_{\reg}u_1 = 0.
	\label{eq:Ureg}
	\end{equation}
	Consider now the solution of 
	\begin{equation}
	a(\Usz,v) = \int_{S_0} \partial_r u_1 v -2c_0\int_{\Bbb P_n}u_1 v, \forall v \in H_0^1(\Bbb P_n), \int_{\Bbb P_n} \Usz u_1 = 0.
	\label{eq:U0sing}
	\end{equation}
	where $2c_0 = \int_{S_0} \partial_r u_1 u_1$ ensures that the right hand side in \eqref{eq:U0sing} is orthogonal on $u_1$. A simple integration by parts shows that $c_0 = u_1(0)^2/4$. In a similar way one can define $\Usm$ and $\Usp$ with the right hand side as integrals on the segments $S_-$ and $S_+$. Then by linearity we have
	\[ \bo U = \bo U_{reg} + \begin{pmatrix}
	-\frac{1}{\tan \theta} \Usp-\frac{1}{\tan \theta} \Usm +\frac{2}{\tan \theta} \Usz \\
	-\Usp + \Usm
	\end{pmatrix}.\]
	The singular behavior of the functions $(U^1,U^2)$ is underlined in the illustrations given in Figure \ref{fig:material}. Normal derivatives of $U^{1,2}$ across certain segments inside $\Pp$ are discontinuous. Nevertheless, the numerical computation suggests that piecewise $H^2$ regularity holds on triangles defined by the symmetric triangulation shown in Figure \ref{fig:simple-triangulations}.
	
	Using results in \cite[Chapter 2]{grisvard} it is classical that since $\Bbb P_n$ is convex, $\bo U_{\reg} \in H^2(\Bbb{P}_n)$. Following \cite[Lemma 4.3.1.3]{grisvard} we know that 
	\[ \|D^2 \bo U_{\reg}\|_{L^2} = \|\Delta \bo U_{\reg}\|_{L^2},\]
	allowing us to obtain an \emph{a priori} error bound for the associated finite element problem. The key idea is that for functions $v,w \in H^1(\Pp)$ such that $(v,w)$ is orthogonal to a fixed vector for every side of $\Pp$, it is possible to interchange derivatives in the integration: 
	\[ \int_{\Pp} \partial_x v\partial_y w = \int_{\Pp} \partial_y v \partial_x w.\]
	Since $\bo U_\reg$ is zero on the boundary of $\Pp$, its tangential gradient vanishes and applying the previous equality to $v = \partial_y \bo U_\reg$, $w = \partial_x \bo U_\reg$ gives the $H^2(\Pp)$ identity stated above.  
	
	For the singular part the estimate is more involved. The functions $\Usz, \Usm, \Usp$ defined by \eqref{eq:U0sing} do not belong to $H^2(\Bbb{P}_n)$ as their gradients are discontinuous on the segments $S_0, S_-, S_+$, respectively. More precisely, {the tangential gradients are continuous accross these segments, while the normal components have jumps}. One can notice that $\Usp$ and $\Usm$ are rotations of $\Usz$ with an angle equal to $\pm 2\pi/n$. Therefore, initially, it is enough to concentrate the analysis on $\Usz$. First, note that $\Usz$ is even with respect to the $y$ variable, in view of the symmetry of the domain $\Bbb{P}_n$, the orthogonality constraint and the fact that $S_0$ lies on the horizontal axis, which is a symmetry line of $\Bbb P_n$. 
	
	Even though $\Usz$ is not in $H^2(\Bbb P_n)$, it is $H^2$ regular on both sides of the symmetry axis. Let us denote $\Bbb P_+$ the part of $\Bbb P_n$ contained in the upper half plane $\{y > 0\}$. Following results in \cite{grisvard} summarized in \cite[Theorem 2.1]{mghazli} it is immediate that $\Usz|_{\Bbb P_+}$ is indeed $H^2(\Bbb P_+)$, when restricted to $ \Bbb P_+$. 
	
	In order to be more precise, consider the problem on half the polygon. The restriction of $\Usz$ on $\Ppp$, denoted by $\Uszh$ verifies the equation
	\begin{align}
	 \int _{\Ppp}\nabla \Uszh \nabla v -\lb_1 \int _{\Ppp}  \Uszh   v 
	  = \frac{1}{2}\int_{S_0} \partial_r u_1v-c_0\int_{\Ppp} u_1 v, \\ \forall v \in H^1(\Bbb P_n),   v = 0 \text{ on } \partial \Bbb P_n\cap \Bbb P_+, \int_{\Ppp}\Uszh u_1 = 0.
	\label{eq:U0plus}
	\end{align}
	Note that on $\Ppp$, considering Neumann boundary condition on the symmetry axis, $\lambda_1$ is still an eigenvalue with the same eigenfunction $u_1$ (keeping the same notations for restrictions). Moreover, $\lambda_2(\Pp)$ is also the second eigenvalue of $\Ppp$ with the Dirichlet-Neumann conditions described above, corresponding to a symmetric eigenfunction for $\lambda_2(\Pp)$ with respect to $y$. 
	
	Indeed, this can be justified as follows: take a second eigenfunction for the mixed Dirichlet-Neumann problem on $\Ppp$. By reflection  across the horizontal axis, we build an eigenfunction on the full polygon $\Pp$, changing sign. This means that $ \lambda_2(\Ppp)$ is not smaller than the second eigenvalue of the polygon $\Pp$. On the other hand, take a second eigenfunction $u_2$ for the Dirichlet problem of the polygon $\Pp$.  For convenience, we see $u_2$ in polar coordinates $(\rho, \theta)$. Then $\tilde u_2(\rho, \theta) := u_2(\rho, -\theta)$ is also a second eigenfunction of the polygon. If  $u_2+\tilde u_2\not \equiv 0$, we have build a new second eigenfunction having $0$ normal derivative on the horizontal line. So the second eigenvalue of the polygon is also an eigenvalue for the mixed problem on $\Ppp$, of index higher than $1$ since it has to change sign, hence we conclude. If  $u_2+\tilde u_2 \equiv 0$, then   $u_2$ is a second eigenfunction of the polygon vanishing on the horizontal line.  Let us consider the function $\overline u_2(\rho, \theta) = u_2(\rho, \theta+ \frac{2\pi}{n})$. This is another second eigenfunction which cannot vanish on the horizontal line, otherwise it would have more than two nodal domains. So we can reproduce the first argument. 	
	
	Notice the following:
	\begin{itemize}
		\item $\Uszh=0$ on $\partial \Pp \cap \{y\geq 0\}$
		\item $\partial_n \Uszh = 0$ on $\partial \Ppp \cap \{(x,0): x\leq 0\}$
		\item $\partial_n \Uszh = \partial_y \Uszh = \frac{1}{2} \partial_x u_1$ on $\partial \Ppp \cap \{(x,0):x\geq 0\}$. 
	\end{itemize}
    Since $\partial_x u_1(0,0) = \partial_x u_1(1,0)=0$ and $u_1\in H^2(\Ppp)$ it follows that $\partial_x u_1 \in H^{1/2}(\partial \Ppp \cap \{y=0\})$. It follows that $\partial_n \Uszh \in H^{1/2}(\Ppp\cap \{y=0\})$. Therefore $\Uszh \in H^2(\Ppp)$. As a direct consequence, the solution $\bo U=(U^1,U^2)$ of \eqref{eq:U0} is locally in $H^2$ when restricted to $T_+, T_-$ or $\Bbb{P}_n \setminus (T_+\cup T_-)$.

	As a first step, let us get an \emph{a priori} estimate on $\|D^2 \Uszh\|_{L^2(\Bbb P_+)}$. The more problematic term is the mixed term involving mixed derivatives of the form $\partial_{xy}$. Below, an approach similar to \cite[Lemma 4.3.1.3]{grisvard} is used. 
	
	{	
	\begin{lemma}\label{lem:grisvard} We have:
		\[ \int_{\Ppp} \partial_{xx}\Uszh \partial_{yy} \Uszh - \int_{\Ppp} \partial_{xy} \Uszh \partial_{xy} \Uszh = \frac{1}{2}\langle  \partial_\tau \tilde w, \tilde v\rangle_{H^{-1/2}(\R), H^{1/2}(\R)}
		\]
		where $\tilde w$ is the extension by $0$ of $\partial_x u_1 |_{S_0}$ on $\R\sm S_0$ and $\tilde v$ is any extension in $H^{\frac 12+\delta} (\R)$ of $\partial _x \ov {U_{S_0}} |_{S_0}$.	
	\end{lemma}

\emph{Proof:} Following \cite[Section 4.3]{grisvard}, when $\Omega$ is a convex polygon $\bo a_0,...,\bo a_{m-1}$ with edges $\Gamma_i=[\bo a_i,\bo a_{i+1}]$, $0\leq i \leq m-1$ (indices considered modulo $n$ when necessary), considering an orientation of $\partial \Omega$ in the trigonometric sense, for $\sigma>1$ define the space
\[ G^{\sigma}(\Omega) = \{(v,w) \in H^{\sigma}(\Omega)\times H^{\sigma}(\Omega) : (\eta_i,\mu_i) \cdot (v,w) = 0 \text{ on }\Gamma_j, 0\leq j \leq m-2, w(\bo a_0)=w(\bo a_{m-1})=0\},\]
where $\eta_i^2+\mu_i^2\neq 0$, $0\leq i \leq m-2$. Assume furthermore that $\mu_0\mu_{m-2}\neq 0$. Note that compared with \cite[Section 4.3]{grisvard}, one side is omitted in the definition above and pointwise values are imposed for $w$ at $\bo a_0$ and $\bo a_{m-1}$. Note that the definition of $G^{\sigma}$ is meaningful, since pointwise values can be defined for $H^\sigma$, when $\sigma>1$.

Consider $(v,w) \in G^2(\Omega)$ and apply Green's formula:
\begin{equation}\label{eq:bb201}
\int_\Omega \partial_x v\partial_y w -\int_\Omega \partial_y v\partial_x w = \int_{\partial \Omega} v \partial_\tau w,
\end{equation}
where $\partial_\tau$ indicates the tangential derivative. On sides $\Gamma_i$ where $\mu_i\neq 0$ we have $w = -\frac{\eta_i}{\mu_i} v$
which implies
\begin{equation}\label{eq:grisvard-gamma} \int_{\Gamma_i} v\partial_\tau w = -\frac{\eta_i}{\mu_i} \int_{\Gamma_i} v\partial_\tau v = -\frac{\eta_i}{2\mu_i} (v^2(\bo a_{i+1})-v^2(\bo a_{i})).
\end{equation}
Note that this identity is meaningful for $v\in H^2(\Omega)$, since $v \in C(\overline \Omega)$. If $\mu_i=0$ then by definition we have $v = 0$ on $\Gamma_i$, and by continuity $v(\bo a_i)=v(\bo a_{i+1})=0$. 
Therefore when summing contributions \eqref{eq:grisvard-gamma} only vertices $\bo a_0, \bo a_{m-1}$ and those for which the neighboring segments $\Gamma_{i-1},\Gamma_i$ have the coefficients $\mu_{i-1}, \mu_i$ different from zero remain. 
It is immediate to notice that 
\begin{align*}  \int_{\partial \Omega} v \partial_\tau w &= \int_{\Gamma_{m-1}}v\partial_\tau w + \frac{1}{2}\sum_{\mu_i\mu_{i-1}\neq 0} \left( \frac{\eta_i}{\mu_i}-\frac{\eta_{i-1}}{\mu_{i-1}}\right) v^2(\bo a_i)+\frac{\eta_0}{\mu_0} v^2(\bo a_0)-\frac{\eta_{m-1}}{\mu_{m-1}}v^2(\bo a_{m-1}).
\end{align*}
Note that if $\left( \frac{\eta_i}{\mu_i}-\frac{\eta_{i-1}}{\mu_{i-1}}\right)\neq 0$ then $(v,w)$ is orthogonal to two non-colinear vectors at $\bo a_i$, giving $v(\bo a_i)=0$. Therefore
\begin{equation}\label{eq:grisvard-seg}
  \int_\Omega \partial_x v\partial_y w -\int_\Omega \partial_y v\partial_x w = \int_{\Gamma_{m-1}} v \partial_\tau w +\frac{\eta_0}{\mu_0} v^2(\bo a_0)-\frac{\eta_{m-1}}{\mu_{m-1}}v^2(\bo a_{m-1})
\end{equation}

Let now $\delta >0$ (and small) and consider $(v,w) \in G^{1+\delta}(\Omega)$. Since $H^{1+\delta}(\Omega) \subset C(\overline \Omega)$ continuously, a variant of \eqref{eq:grisvard-seg} can be written, as follows.   From the proof in \cite[Lemma 4.3.1.3]{grisvard} it follows that $G^2(\Omega)$ is dense in $G^{1+\delta}(\Omega)$. Consequently, we can approach $(v,w)$ strongly in $H^{1+\delta} (\Om)$ by  sequences of functions $(v_\vps,w_\vps)_\vps \sq G^{2}(\Omega)$ and hence write \eqref{eq:grisvard-seg} for $(v_\vps, w_\vps)$. We also have
$$ \int_\Omega \partial_x v_\vps\partial_y w_\vps -\int_\Omega \partial_y v_\vps\partial_x w_\vps \to \int_\Omega \partial_x v\partial_y w -\int_\Omega \partial_y v\partial_x w,$$
$$\frac{\eta_0}{\mu_0} v_\vps^2(\bo a_0)-\frac{\eta_{m-1}}{\mu_{m-1}}v_\vps^2(\bo a_{m-1})\to \frac{\eta_0}{\mu_0} v^2(\bo a_0)-\frac{\eta_{m-1}}{\mu_{m-1}}v^2(\bo a_{m-1}).$$
In order to pass to the limit the term $\int_{\Gamma_{m-1}} v_\vps \partial_\tau w_\vps$, some observations are in order. Since $w_\vps$ vanishes at $\bo a_0,\bo a_{m-1}$, integration by parts gives $\int_{\Gamma_{m-1}} v_\vps \partial_\tau w_\vps = - \int_{\Gamma_{m-1}}  \partial_\tau v_\vps w_\vps$. Moreover, denoting $\tilde w_\vps\in H^{ 1}(\R)$ the extension by $0$ of $w_\vps$ on $\R \sm \G_{m-1}$ (here we assume $\G_{m-1} \sq \R \times \{0\}$ and identify $\R\times \{0\}$ with $\R$) and $\tilde v_\vps$ any extension in $H^{1}(\R)$ of $v_\vps$, we get
$$ - \int_{\Gamma_{m-1}}  \partial_\tau v_\vps w_\vps = - \int_{\R}  \partial_\tau \tilde v_\vps \tilde w_\vps=\int_{\R}   \tilde v_\vps \tilde \partial_\tau w_\vps= \langle  \partial_\tau \tilde w_\vps, \tilde v_\vps\rangle_{H^{-1/2}(\R), H^{1/2}(\R)}.$$
For instance, we shall take for $\tilde v_\vps$ any $H^2$ extension of $v_\vps$ outside $\Om$ which is uniformly bounded in $H^{1+\delta} (\R^2)$. Then we have that $\tilde v_\vps \to \tilde v$ and $\tilde w_{\vps} \to \tilde {w}$ strongly in $H^{\frac 12+\delta}(\R)$, where $\tilde v$ is an extension on $\R\sm \G_{m-1}$ of $v$ and $\tilde w$ is the extension by $0$ of $w$ on $\R\sm \G_{m-1}$. Consequently, 
$$\langle  \partial_\tau \tilde w_\vps, \tilde v_\vps\rangle_{H^{-1/2}(\R), H^{1/2}(\R)} \to \langle  \partial_\tau \tilde w, \tilde v \rangle_{H^{-1/2}(\R), H^{1/2}(\R)},$$
hence
 
\begin{equation}\label{eq:grisvard-seg2}
\int_\Omega \partial_x v\partial_y w -\int_\Omega \partial_y v\partial_x w = \langle  \partial_\tau \tilde w, \tilde v\rangle_{H^{-1/2}(\R), H^{1/2}(\R)} +\frac{\eta_0}{\mu_0} v^2(\bo a_0)-\frac{\eta_{m-1}}{\mu_{m-1}}v^2(\bo a_{m-1})
\end{equation}

	Apply \eqref{eq:grisvard-seg2} for $\Omega = \Ppp$. Let $k = \lfloor n/2\rfloor$ and denote $\bo b = \bo a_k$ if $n$ is even or the midpoint of $[\bo a_{k-1},\bo a_{k}]$ if $n$ is odd. Then the vertices of $\Ppp$ are $\bo a_0,...,\bo a_{k-1}$, $\bo b$ and the origin $\bo o$.
	
	Furthermore, $\Gamma_{m-1} = [\bo o \bo a_0]$ and taking 
	$v = \partial_x \Uszh \in H^{1+\delta}(\Ppp), w = \partial_y \Uszh\in H^{1+\delta}(\Ppp)$. On the segment $\Gamma_{m-1}$ we have $w = \frac{1}{2} \partial_x u_1$, therefore $w(\bo o) = w(\bo a_0) = 0$. Note that $w = \partial_y \Uszh$ verifies $ w = 0$ on $\Gamma_{m-2}$, i.e. $(0,1)\cdot (v,w) = 0$. On $\Gamma_0$ the function $\Uszh$ verifies Dirichlet boundary conditions $\Uszh = 0$ thus $\nabla \Uszh \cdot \overrightarrow{\bo a_0 \bo a_1} = 0$. Moreover, we have $\partial_y \Uszh(\bo a_0) = \partial_x u_1(\bo a_0)=0$. Therefore 
	\eqref{eq:grisvard-seg2} holds with $v(\bo a_0)=0$ and $\eta_{m-1} = 0$ and the result follows.	
	\hfill $\square$
}

	\medskip 
	
	\medskip
	\medskip
	
	\begin{lemma}[Extension operator with computed norm]\label{lem:extension}
		Let $\Pp$ be the regular polygon with $n$ sides and circumcircle of unit radius. There exists a linear extension operator 
		\[ \mathcal E : H^1(\Pp) \to  H^1(\Bbb{R}^2)\]
		such that $\|\mathcal E\|\leq C_n$, where $C_5=4$ and for $n\geq 6$
		\[ C_n = \left[4+24\cos^2 \frac{2\pi}{n}\right]^{1/2}.\]
	\end{lemma}

	\begin{rem}
		Note that if such an extension operator exists then the same upper bound $C_n$ can be used for an extension operator on the half plane $\Bbb{R}_+$.
	\end{rem}
	
	\emph{Proof:} Given a function in $H^1(\Pp)$, an extension to $H^1(\Bbb R^2)$ having controlled norm is constructed using reflections and a cutoff function. 

\smallskip
\noindent
	\bo{Reflection.} A distinction needs to be made between the case of the pentagon and $n\geq 6$.
	
	Consider the regular pentagon inscribed in the unit circle with center $o$, denoted $\Bbb P_5 = \bo a_0\bo a_1\bo a_2\bo a_3\bo a_4$. Consider $\bo b_1$ the reflection of $\bo o$ across the line $\bo a_0 \bo a_1$ and $\bo b_{n-1}$ the reflection of $\bo o$ across the line $\bo a_0\bo a_{n-1}$. The region in the angle $\angle \bo b_1\bo a_0 \bo b_{n-1}$ needs to be covered with additional copies of the current figure. 
	
\begin{figure}
		\centering 
		\includegraphics[height=0.35\textwidth]{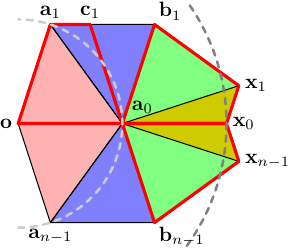}\qquad 
		\includegraphics[height=0.35\textwidth]{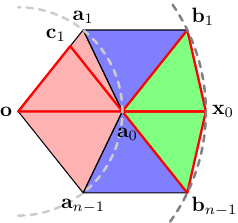}
		\caption{Reflection procedure: for the pentagon at least four reflections are needed (left), for $n \geq 6$ three reflections suffice (right). The thresholds for the application of the cutoff function are also illustrated with dashed lines. The cutoff function equals $1$ inside the small disk, zero outside the big disk and is affine in the radial direction in-between.}
		\label{fig:extension}
	\end{figure}

	To this end construct the isosceles trapeze $\bo T = \bo o \bo a_0 \bo c \bo a_1$ with $\bo c_1 \subset [ \bo a_1 \bo b_1]$. Then place two copies of $\bo T$, $\bo b_1 \bo a_0 \bo x_0 \bo x_1$ and $\bo b_{n-1} \bo a_0 \bo x_0 \bo x_1$, in the angle $\bo b_1 \bo o \bo b_{n-1}$ such that the side congruent to $\bo a_0\bo o$ in $\bo T$ aligns with $\bo a_0 \bo b_1$ and $\bo a_0 \bo b_{n-1}$. See Figure \ref{fig:extension} for the illustration. This reflection is repeated for every vertex of the pentagon. The first stage of the reflection produces a copy of any point in $\Pp$ and the second stage produces at most $4$ extra copies. Thus, globally, the reflection process produces at most $6$ extra copies of any point in $\Bbb P_5$.
	
    For $n \geq 6$ a general extension procedure producing at most $4$ extra copies of any point in $\Pp$ is available. Consider $\bo a_0,...,\bo a_{n-1}$ the vertices of $\Pp$ and construct $\bo b_1,\bo b_{n-1}$ reflections of the origin $\bo o$ across $\bo a_0\bo a_1$ and $\bo a_0\bo a_{n-1}$, respectively.
    
    Since $\bo o \bo a_1 \geq \bo a_0 \bo a_1$ there exists $\bo c_1$ on $\bo o \bo a_1$ such that the triangle $\bo T=\bo o \bo a_0 \bo c_1 $ is isosceles. Consider two copies of $\bo T$, denoted $\bo b_1 \bo a_0 \bo x_0$, $\bo b_{n-1}\bo a_0 \bo x_0$, covering the angle $\angle \bo b_1 \bo a_0 \bo b_{n-1}$ such that the side $\bo a_0 \bo o$ of $\bo T$ overlaps with $\bo a_0 \bo b_1$, $\bo a_0 \bo b_{n-1}$, respectively. This reflection contains at most four copies of any point in $\Pp$. See Figure \ref{fig:extension} for an illustration.
    
  \smallskip
\noindent  \bo{Cutoff.} Consider now for any function $u \in H^1(\Pp)$ the reflected function $\overline u$ following the procedure above and extended arbitrarily afterwards in $H^1(\Bbb{R}^2)$. Consider a cutoff function $\Phi_n$ verifying:
    \[\Phi_n = \begin{cases}
    1 & |x|\leq 1 \\
    0 & |x|\geq 1+\ell_n\\
    (\ell_n+1-|x|)/l_n & |x|\in [1,1+\ell_n],
    \end{cases}\]
    where $\ell_n$ is chosen maximal such that the ball of radius $1+\ell_n$ is included in the polygon obtained from $\Pp$ after the reflection process described in the first step.

    With these considerations, an explicit extension operator is defined for any $\Pp$ by 
    \[ \mathcal E(u) = \Phi_n \overline u.\]
    Let us evaluate the $H^1(\Bbb R^2)$ norm of $\mathcal E(u)$ and compare it with the norm of $u$ in $H^1(\Pp)$. Let $k_n$ be the minimum number of copies for points in $\Pp$ involved in the construction of the reflected polygon. Then 
    \[ \int_{\Bbb{R}^2} \mathcal E(u)^2 = \int_{\Bbb{R}^2} (\Phi_n \overline u)^2 \leq \int_{B(0,1+\ell_n)} \overline u^2 \leq k_n \int_\Pp u^2.\]
    The similar estimation for the gradient gives
    \begin{align*}
     \int_{\Bbb{R}^2} |\nabla \mathcal E(u)|^2 & = \int_{\Bbb{R}^2} |\nabla \Phi_n \overline u+\nabla \overline u \Phi_n|^2 \\
     & \leq 2\int_{\Bbb R^2} |\nabla \varphi|^2 \overline u^2+2\int_{\Bbb{R}^2} |\nabla \overline u|^2 \Phi_n^2 \\
     & \leq \frac{2}{\ell_n^2}\int_{B(0,1+\ell_n)\setminus B(0,1)} \overline u^2 +2\int_{B(0,1+\ell_n)}|\nabla \overline u|^2 \\
    & \leq 2k_n \int_\Pp |\nabla u|^2+\frac{2(k_n-1)}{\ell_n^2}\int_\Pp \overline u^2.
    \end{align*}
    Combining the two estimates gives
    \[ \|\mathcal E(u)\|_{H^1(\Bbb{R}^2)}^2\leq \max\left\{2k_n, \frac{2(k_n-1)}{\ell_n^2}+k_n\right \}.\]

	For the pentagon we have $k_n=6$, $\ell_n = 1$ giving
	$C_5 = 4$. For $n \geq 6$ we have $k_n=4$, $\ell_n = \frac{1}{2\cos \frac{2\pi}{n}}$, giving
		\[ C_n = \sqrt{4+24\cos^2 \frac{2\pi}{n}}.\]
		\hfill $\square$

	Combining the previous results we may now find an explicit estimate for $\|D^2 \Uszh\|_{L^2(\Ppp)}$.

    \begin{thm}\label{thm:estimate-H2}
    	There exists an explicit upper bound for $\|D^2 \Uszh\|_{L^2}$ in terms of $\lambda_1$, $u_1$ and the geometry of $\Ppp$.
    	
    	The explicit finite element error estimate $\|\nabla \Usz-\nabla \Pi_{1,h} \Usz\| \leq C_1h \sqrt{2} \|D^2 \Uszh\|$ holds for the symmetrized $\Usz$, provided the segment $S_0$ is meshed exactly.
    \end{thm}

\emph{Proof:} Consider $\delta >0$, $S = [0,1] \times \{0\}$ and $f,g \in H^{1/2+\delta}(S)$ such that $f ( 0)=f(1)=0$. Consequently, we can extend $f$ by $0$ on $\R \sm(0,1)$ to get that   $f\in  H^{1/2+\delta}(\R)$ and so $f' \in H^{-\frac 12+\delta} (\R)$. To simplify notations, assume   $ g$ is an arbitrary extension of $g$ to $H^{1/2+\delta}(\Bbb R)$ we have
\begin{equation}\label{eq:estimation-Fourier}
 |\langle f',  g \rangle_{H^{-1/2}(\R),H^{1/2}(\R)}|\leq \|f\|_{H^{1/2}(\Bbb{R})} \|  g\|_{H^{1/2}(\Bbb{R})}.
 \end{equation}
 Indeed, taking the Fourier transform and noting that $f$ extends with zero from $S$ to $\Bbb R$ we have
\[  \langle f',  g \rangle_{H^{-1/2}(\Bbb R),H^{1/2}(\Bbb R)}=\int_{\Bbb R} \hat f \hat g' = \int_{\Bbb R} ix \hat f \hat g.  \]
Therefore
\begin{align*}  |\langle f',\tilde g \rangle_{H^{-1/2}(\Bbb R),H^{1/2}(\Bbb R)}| &\leq \int_{\Bbb R}(1+|x|^2)^{1/4} |\hat f| (1+|x|^2)^{1/4}|\hat g| \\
& \leq \left[\int_{\Bbb R}(1+|x|^2)^{1/2} |\hat f|^2  \right]^{1/2} \left[ \int_{\Bbb R}(1+|x|^2)^{1/2}|\hat g|^2 \right]^{1/2}\\
& = \|f\|_{H^{1/2}(\Bbb R)} \|g\|_{H^{1/2}(\Bbb R)}.
\end{align*}

Next, according to \cite[Equation (2.3)]{PaPa11} we have
\begin{equation}\label{eq:trace-H12} \|u\|^2_{H^{1/2}(\Bbb R)} \leq C_{1,1}\|u\|^2_{H^1(\Bbb R^2)}
\end{equation}
with $C_{1,1}= \frac{\Gamma(1/2)}{(4\pi)^{1/2}\Gamma(1)} =\frac{1}{2}$.

Turning to the estimation regarding $D^2 \Uszh$, since $-\Delta \Uszh = -\lambda_1 \Uszh -c_0 u_1$ we have, using Lemma \ref{lem:grisvard} (and its notations for the extensions of $\partial _xu_1, {\partial_x \Uszh}$ ) that 
	\begin{align}
	\|D^2 \Uszh\|_{L^2(\Ppp)}^2 & = \int_{\Ppp} ((\partial_{xx}\Uszh)^2+(\partial_{yy} \Uszh)^2 +2 \partial_{xy} \Uszh\partial_{xy}\Uszh ) \notag \\
	& = \int_{\Ppp} ((\partial_{xx}\Uszh)^2+(\partial_{yy})^2 \Uszh+2\partial_{xx}\Uszh \partial_{yy}\Uszh)-  \langle  \partial_\tau \tilde w, \tilde v\rangle_{H^{-1/2}(\R), H^{1/2}(\R)}\notag \\
	& = \|\Delta \Uszh\|_{L^2(\Ppp)}^2-  \langle  \partial_\tau \tilde w, \tilde v\rangle_{H^{-1/2}(\R), H^{1/2}(\R)} \notag \\
	& =\|\lambda_1 \Uszh+c_0 u_1\|_{L^2(\Ppp)}^2-  \langle  \partial_\tau \tilde w, \tilde v\rangle_{H^{-1/2}(\R), H^{1/2}(\R)}. 
	\label{eq:estimateUS0}
	\end{align}

To estimate the first term above start from
\[
\|\lambda_1 \Uszh+c_0 u_1\|_{L^2(\Ppp)}^2=\lambda_1^2 \int_{\Ppp} \Uszh^2+\frac{c_0^2}{2}.
\]
Using $\Uszh$ as test function in the variational formulation for $\Uszh$ gives
\[ \int_{\Ppp} |\nabla \Uszh|^2-\lambda_1 \int_\Omega \Uszh^2 = \int_{S_0} \partial_r u_1 \Uszh.\]
and the trace theorem implies
\[ \int_S \Uszh^2 \leq \|\nabla \Uszh\|_{L^2(\Ppp)}^2.\]
The inequality above is valid for any function $U$ which is zero on $\partial \Pp$ through direct application of the mean value theorem on vertical slices. The constant in the trace theorem can be taken equal to $1$ since the polygon $\Pp$ is included in the band $\{|y|\leq 1\}$. 

Therefore 
\[
\int_{S_0} \partial_r u_1 \Uszh \leq \|\nabla \Uszh\|_{L^2(\Ppp)} \|\partial_r u_1\|_{L^2(S_0)}.
\]
Poincar\'e's inequality on the orthogonal of $u_1$ in $H^1(\Ppp)\cap H_0^1(\Pp)$ gives
\[ \left(1-\frac{\lambda_1}{\lambda_2}\right) \|\nabla \Uszh\|_{L^2(\Ppp)}^2 \leq  \|\nabla \Uszh\|_{L^2(\Ppp)} \|\partial_r u_1\|_{L^2(S)}.\]	
Therefore, again using Poincar\'e's inequality
\[ \sqrt{\lambda_2}\|\Uszh\|_{L^2(\Ppp)}\leq \|\nabla \Uszh\|_{L^2(\Ppp)} \leq \frac{\lambda_2}{\lambda_2-\lambda_1} \|\partial_r u_1\|_{L^2(S)},\]
giving the \emph{a proiri} bound
\[ \|\Uszh\|_{L^2(\Ppp)}^2 \leq \frac{\lambda_2}{(\lambda_2-\lambda_1)^2} \|\partial_r u_1\|_{L^2(S)}^2.\]

 For the second term in \eqref{eq:estimateUS0}, using estimate \eqref{eq:estimation-Fourier} gives	
\[  | \langle  \partial_\tau \tilde w, \tilde v\rangle_{H^{-1/2}(\R), H^{1/2}(\R)} |  \leq \|\partial_x u_1\|_{H^{1/2}(\Bbb R)}\|\partial_x \Uszh\|_{H^{1/2}(\Bbb{R})},\]
where $\partial_x u_1$ is extended with zero and $\partial_x \Uszh$ is extended arbitrarily from $\Ppp \cap \{y=0\}$ to $\Bbb{R}$. Using the same extension operator defined in Lemma \ref{lem:extension} and \eqref{eq:trace-H12} gives
\[ | \langle  \partial_\tau \tilde w, \tilde v\rangle_{H^{-1/2}(\R), H^{1/2}(\R)}| \leq \frac{C_n^2}{2} \|\partial_x u\|_{H^1(\Bbb P_n^+)} \|\partial_x \Uszh \|_{H^1(\Bbb P_n^+)}\leq \]
\[ \leq  \frac{C_n^2}{2} \sqrt{(\lambda_1^2+\lambda_1)/2} \sqrt{\|D^2 \Uszh\|_{L^2(\Ppp)}^2+\|\nabla \Uszh\|_{L^2(\Ppp)}^2},\]
where the equality $\|D^2 u_1\|^2_{L^2(\Pp)} = \|\Delta u\|^2_{L^2(\Pp)} = \lambda_1^2$ was used.

Combining the previous estimates we obtain
\begin{align*}
\|D^2 \Uszh\|_{L^2(\Ppp)}^2\leq  & \lambda_1^2 \|\Uszh\|_{L^2(\Ppp)}^2+\frac{c_0^2}{2}  +\frac{C_n^2}{2} \sqrt{\frac{\lambda_1^2+\lambda_1}{2}} \sqrt{\|D^2 \Uszh\|_{L^2(\Ppp)}^2+\|\nabla \Uszh\|_{L^2(\Ppp)}^2}
\end{align*}

Denoting $X = \|D^2 \Uszh\|_{L^2(\Ppp)}$ and adding an upper bound for $\|\nabla \Uszh\|_{L^2(\Ppp)}^2$ to the previous inequality leads to a quadratic inequality of the form
\[ X^2+C \leq B+A\sqrt{X^2+C},\]
which leads to $X^2 \leq \left( \frac{A+\sqrt{A^2+4B}}{2}\right)^2-C$, where:
\begin{itemize}[noitemsep]
	\item $A$ is an upper bound for $\frac{C_n^2}{2} \sqrt{\frac{\lambda_1^2+\lambda_1}{2}}$
	\item $B$ is an upper bound for $\lambda_1^2 \|\Uszh\|_{L^2(\Ppp)}^2+\frac{c_0^2}{2}+\|\nabla \Uszh\|_{L^2(\Ppp)}^2$, $c_0 = u_1(0)^2/4= \|u_1\|_{L^\infty}^2/4 <\lambda_1/4$ (see  Grebenkov   \cite[Formula (6.22)]{GrNg13})
	\item $C$ is an upper bound for $\|\nabla \Uszh\|_{L^2(\Ppp)}^2$
\end{itemize}
This leads to an explicit upper bound for $\|D^2 \Uszh\|_{L^2(\Ppp)}$. The fact that such an upper bound can be found also shows implicitly that $\Uszh$ is indeed in $H^2(\Ppp)$. 

The second part of the statement of the theorem follows at once observing that $\|D^2 \Usz\|_{L^2(\Pp)} = \sqrt{2} \|D^2\Uszh\|_{L^2(\Ppp)}$ and applying the interpolation estimate \eqref{eq:interpolation}. The $H^2$ norm of the symmetrized $\Usz$ is understood in a piece-wise sense, considering a crack on the segment $S_0$.
\hfill $\square$

\smallskip
\noindent	
{\bf General singular distribution on rays.} Let the singular distribution be given by
\begin{equation}\label{eq:fsing-general} (f,v) = \sum_{i=0}^{n-1} q_i\left(\int_{S_i} \partial_r u_1 v\right)
\end{equation}
with $\sum_{i=0}^{n-1} q_i = 0$. Then any solution of
\begin{equation}\label{eq:full-rays}
U \in H_0^1(\Pp), \ \ a(U,v) = (f,v)_{H^{-1},H^1}, \int_\Pp U u_1 = 0
\end{equation}
is a linear combination of rotations of the singular solution $\Usz$. Considering the analogue equation for the Laplace operator, an improved estimate compared to 
Lemma 5.1 from \cite{Bogosel_Bucur_Polya} is found.
\begin{lemma}\label{lem:general-estimate-rays}
	Let $U\in H_0^1(\Pp)$ be the solution of 
	\[-\Delta U = f \text{ in }H_0^1(\Pp)\]
	with $f$ of the form \eqref{eq:fsing-general}. Consider $P_h(U) \in \mathcal V_h\subset H_0^1(\Pp)$ the discrete solution corresponding to the same right hand side
	\[
	 -\Delta P_h(U) = f \text{ in }\mathcal V_h.
	\] 
	Then denoting $C(\bo q) = \sum_{i=0}^{n-1} |q_i| C_1 \sqrt{2}\|D^2 \Uszh\|_{L^2(\Ppp)}$ gives
	\begin{equation}\label{eq:estimate-sing-grad}
	\|\nabla U- \nabla P_h(U)\|_{L^2} \le C(\bo q)h \end{equation} 
	and
	\begin{equation}\label{eq:estimate-sing-L2}
	\|  U-   P_h(U)\|_{L^2} \le C_1C(\bo q)h^2. \end{equation} 
\end{lemma}
\begin{proof}
	Simply use the estimates above and the Aubin-Nitsche lemma {\cite[p.136]{CiarletBook}}.
\end{proof}
	
This gives rise to the following estimate for solutions of problems similar to those involving the material derivatives.	
\begin{thm}\label{thm:estimate-singular}
	Let $U\in H^1_0(\Pp)$  be the solution of 
	\begin{equation}\label{bobu28.1}
	\left\{ \begin{array}{rcll}
	-\Delta U - \lb_1U& =&f& \text{ in }  \Pp\\
	U&= &0 & \text{ on }\partial   \Pp\\
	\int_{\Pp} u_1U dx&= &0
	\end{array} \right.
	\end{equation}
	where $(f,u_1)_{H^{-1},H_0^1} = 0$, $f= f_{\reg}+{f_{\sing}} $ with   $f_{\reg} \in L^2(\Pp)$ and $f_{\sing}$ given by \eqref{eq:fsing-general}. Assume $f_h$ is a numerical approximation in $H^{-1}(\Pp)$ of $f$ which {satisfies} $(f_h,u_{1,h})_{H^{-1},H_0^1}=0$ and $(u_{1,h}, \lb_{1,h})$ a numerical approximation of $(u_1, \lb_1)$ in $H^1_0(\Pp)\times \R$. Denote
	$U_h $ the finite element solution in $\mathcal V_h\subset H_0^1(\Pp)$ for
	\begin{align}\forall v \in \mathcal  V_h, \quad \int_{\Pp} (\nabla U_h \cdot \nabla v-\lambda_{1,h} U_hv) \, dx& =( f_h, v)_{H^{-1}\times H^1_0} \label{eq:material-eig-decomposed.h}
	\end{align}
	together with the normalization
	\begin{equation}
	\int_{\Pp} u_{1,h}U_h  \, dx = 0,
	\label{eq:normalization.h}
	\end{equation}
	Then, {denoting by $V$ the solution of $-\Delta V = \lambda_1 U+f$ in $\mathcal V_h$, we have}
	\begin{align*}
	\|\nabla U-\nabla U_h\|_{L^2} & \le  C_1 h \| \lb_1 U+ f_{\reg}\|_{L^2}+  C(\bo q)h\\
	&+   \lambda_{1,h}^{\frac{1}{2}}\Big( (C_1 h)^2 \| \lb_1 U+ f_{\reg}\|_{L^2}
	+ C_1 C(\bo q)h^2)	+\|V\|_{L^2}\|u_{1,h}-u\|_{L^2}\Big)\\
	&+\frac{\lb_{2,h}^\frac 12}{\lb_{2,h}-\lb_{1,h}} \Big (|\lb_{1,h}-\lb_1| \|U\|_{L^2} + \lb_{1,h}C_1 C(\bo q)h^2 \hfill + (1+ \lb_{2,h})^\frac 12\|f-f_h\|_{H^{-1}}\Big).
	\end{align*}
	The analogue $L^2$ estimate follows
	\begin{align*}
	\|U-U_h\|_{L^2} & \leq C_1h  (C_1 h \| \lb_1 U+ f_{\reg}\|_{L^2}+  C(\bo q)h)\\
	& +  \Big( (C_1 h)^2 \| \lb_1 U+ f_{\reg}\|_{L^2}
	+ C_1 C(\bo q)h^2)	+\|V\|_{L^2}\|u_{1,h}-u\|_{L^2}\Big)\\
	& +\frac{1}{\lb_{2,h}-\lb_{1,h}} \Big (|\lb_{1,h}-\lb_1| \|U\|_{L^2} + \lb_{1,h}C_1 C(\bo q)h^2 \hfill + (1+ \lb_{2,h})^\frac 12\|f-f_h\|_{H^{-1}}\Big).
	\end{align*}
\end{thm}
\begin{proof} The proof follows the general estimate strategy: project onto $\mathcal V^h$ using discrete problems with continuous right hand side, then consider discrete problems with different right hand sides. A slight complication arises from the normalization conditions. 

\smallskip
\noindent	
\bo{a) Interpolation error.} We denote $U_{\reg}, U_{\sing}$ the solutions of 
	$$U_{\reg} \in H^1_0(\Pp), \;\; -\Delta U_{\reg}= \lb_1 U+ f_{\reg}, \quad U_{\sing}\in H^1_0(\Pp),\;\; -\Delta U_{\sing} =   f_{\sing}, $$
	so that $U=U_{\reg}+U_{\sing}$.	We introduce the auxiliary functions
	$V_{\reg}, V_{\sing} \in \mathcal  V^h$, $V_{\reg}, V_{\sing}$ the finite element solutions of
	$$V_{\reg}\in \mathcal  V^h, \quad -\Delta V_{\reg} = \lb_1 U +f_{\reg} \quad V_{\sing} \in \mathcal  V^h, \quad -\Delta V_{\sing}= f_{\sing}.$$
	For $V_{\sing}$, the estimate \eqref{eq:estimate-sing-grad} from Lemma \ref{lem:general-estimate-rays} holds, and gives   
	$$\|\nabla U_{\sing}-\nabla V_{\sing}\|_{L^2} \le  C(\bo q) h,$$
	while for $V_{\reg}$ the classical estimate gives
	$$\|\nabla U_{\reg}- \nabla V_{\reg}  \|_{L^2} \le C_1 h \| \lb_1 U+ f_{\reg}\|_{L^2}.$$
	
\smallskip
\noindent	\bo{b) Normalization error.} Let us denote $V=V_{\reg}+V_{\sing}$ and define $\tilde V= V-(\int_{\Pp} Vu_{1,h} dx ) u_{1,h}$. Then we have 
	\begin{equation}\label{bobu80}
	\|\nabla \tilde V-\nabla V\|_{L^2} = \lambda_{1,h}^{\frac{1}{2}} 
	\left|\int_{\Pp} (Vu_{1,h} -Uu_1) \right| \leq \lambda_{1,h}^{\frac{1}{2}}(\|U-V\|_{L^2}+\|V\|_{L^2}\|u_{1,h}-u\|_{L^2}).
	\end{equation}
	The Aubin-Nitsche lemma gives $\|U-V\|_{L^2} \leq C_1h \|\nabla U-\nabla V\|_{L^2}$.

\smallskip
\noindent	\bo{c) Discrete problems with different right hand sides.} We have that $\tilde V$ is the finite element solution of
	$$\tilde V \in \mathcal  V^h, \quad -\Delta \tilde V-\lb_{1,h} \tilde V= \lb_1 U +f- \lb_{1,h} V,$$
	which gives
	$$\int_{\Pp} |\nabla \tilde V -\nabla U_h|^2 dx - \lb_{1,h}\int_{\Pp} | \tilde V -U_h|^2 dx= ( \lb_1 U +f- \lb_{1,h} V-f_h, \tilde V-U_h)_{H^{-1}\times H^1_0}.$$
	By the Poincar\'e inequality in the orthogonal of $u_{1,h}$ in $\mathcal V_h$ we get
	$$\left(1- \frac{\lb_{1,h}}{\lb_{2,h}} \right) \int_{\Pp} |\nabla \tilde V -\nabla U_h|^2 dx\le \|\lb_1 U-\lb_{1,h} V\|_{L^2}\| \tilde V-U_h\|_{L^2}+ \|f-f_h\|_{H^{-1}} \|  \tilde V -  U_h\|_{H^1}$$
	$$\le \|\lb_1 U-\lb_{1,h} V\|_{L^2} \lambda_{2,h}^{-\frac{1}{2}} \| \nabla \tilde V-\nabla U_h\|_{L^2}+  \|f-f_h\|_{H^{-1}}\left(1+\frac{1}{\lb_{2,h}}\right)^\frac 12  \| \nabla \tilde V-\nabla U_h\|_{L^2}.$$
	Finally,
	\begin{equation}\label{bobu81}
	\| \nabla \tilde V-\nabla U_h\|_{L^2(\mathcal P_n)}\hskip 9cm
	\end{equation}
	$$\le \frac{\lb_{2,h}^\frac 12}{\lb_{2,h}-\lb_{1,h}} \Big (|\lb_{1,h}-\lb_1| \|U\|_{L^2} + \lb_{1,h}\|U-V\|_{L^2}+ (1+ \lb_{2,h})^\frac 12\|f-f_h\|_{H^{-1}}\Big).$$
	The $L^2$ estimate follows along the same lines. 
\end{proof}

\begin{rem}
	It can be observed that all estimates in the proof of Theorem \ref{thm:estimate-singular} can be given in terms of $\|\nabla U-\nabla V\|_{L^2},\|U-V\|_{L^2}, \|u-u_{1,h}\|, |\lambda_{1,h}-\lambda_1|, \|f-f_h\|_{H^{-1}}$. All these quantities can be evaluated explicitly in terms of the estimates for the eigenfunctions, eigenvalues and interpolation errors. 
\end{rem}

Using the notations of Theorem \ref{thm:estimate-singular}, consider two generic problems with solutions $U^\aaa, U^\bbb$ corresponding to the right hand sides $f^\aaa, f^\bbb$ (not necessarily those explicited in the previous section). Consider the associated notations $V^\aaa, V^\bbb$,  $\tilde V^\aaa, \tilde V^\bbb$, $U_h^\aaa, U_h^\bbb$, $f_h^\aaa, f_h^\bbb$ analogue to those introduced in the proof of Theorem \ref{thm:estimate-singular}. Consider the discrete and continuous the bilinear forms
$$a: H^1_0(\Bbb P_n) \times H^1_0(\Bbb P_n) \ra \R,\;\; a(u,v)= \int \nabla u \cdot \nabla v -\lb_1 \int uv,$$
$$a_h: {\mathcal V}^h \times{\mathcal V}^h \ra \R, \;\; a_h(u,v)= \int \nabla u \cdot \nabla v -\lb_{1,h} \int uv.$$

\subsection{Estimations for Hessian eigenvalues}
\label{sec:hess-eig-estimates}

The goal is to estimate the difference between analytical eigenvalues of the Hessian matrix given in Theorem \ref{thm:eig-hessian} and those computed using finite elements. This involves terms of the form
$$|a(U^\aaa,U^\bbb)-a_h(U_h^\aaa, U_h^\bbb)|.$$
Following \cite{Bogosel_Bucur_Polya} and the classical triangle inequality the following quantities are estimated, following the structure of the proof in Theorem \ref{thm:estimate-singular}:

\smallskip
\noindent \bo{First term (interpolation errors).} Observing that $U^\aaa, U^\bbb$ and $V^\aaa, V^\bbb$ solve problems having the right hand sides $f^\aaa, f^\bbb$ in $H_0^1(\Pp)$ and $\mathcal V^h$, respectively, it follows that $a(U^\aaa,V^\aaa) = a(V^\aaa,V^\bbb)= a(V^\aaa,U^\bbb)$. Therefore
\[ |a(U^\aaa,U^\bbb)-a(V^\aaa,V^\bbb)| = |a(U^\aaa-V^\aaa,U^\bbb-V^\bbb)|\leq \|\nabla (U^\aaa-V^\aaa)\|_{L^2} \|\nabla (U^\bbb-V^\bbb)\|_{L^2},\]
giving an estimate of order $h^2$. Note that the above follows from $a(v,v)\leq \int_\Pp |\nabla v|^2$ and the Cauchy-Schwarz inequality for the positive bilinear form $a(\cdot,\cdot)$.

\smallskip
\noindent \noindent{\bf Second term (normalization errors).}
$$|a(V^\aaa,V^\bbb)-a_h(\tilde V^\aaa, \tilde V^\bbb)| \le \|\nabla V^\aaa\|_{L^2} \|\nabla V^\bbb-\nabla \tilde V^\bbb\|_{L^2}+ \|\nabla \tilde V^\bbb\|_{L^2} \|\nabla V^\aaa-\nabla \tilde V^\aaa\|_{L^2}+  $$
$$\hskip 2.5cm |\lb_{1,h}-\lb_1| \|\tilde V^\aaa\|_{L^2}\|\tilde V^\bbb\|_{L^2}+ \lb_1 \|\tilde V^\bbb\|_{L^2} \| V^\aaa- \tilde V^\aaa\|_{L^2}+\lb_1 \|V^\aaa\|_{L^2} \| V^\bbb-\tilde  V^\bbb\|_{L^2},$$
which, in view of inequality \eqref{bobu80}, leads to an approximation of order $h^2$, since the estimate depends on $\|U-V\|_{L^2}$ and $\|u_{1,h}-u_1\|_{L^2}$. 

\smallskip
\noindent
\noindent{\bf Third term (discrete errors).}
$$|a_h(\tilde V^\aaa, \tilde V^\bbb)-a_h(U_h^\aaa, U_h^\bbb)|\le |a_h(\tilde V^\aaa, \tilde V^\bbb-U^\bbb_h)|+ |a_h(\tilde V^\aaa-U^\aaa_h, U_h^\bbb)|\le $$
$$ \|\nabla \tilde V^\aaa\|_{L^2} \|\nabla \tilde V^\bbb-\nabla U^\bbb_h\|_{L^2}+ \|\nabla  U^\bbb_h\|_{L^2} \|\nabla \tilde V^\aaa-\nabla U^\aaa_h\|_{L^2}.$$
The last inequality is a consequence of the fact that $a_h(\cdot, \cdot)$ is a scalar product on $\{u_{1,h}\}^\perp$ in ${\mathcal V}^h$ and of the Cauchy-Schwarz inequality together with the observation that $a_h(v,v) \le \int |\nabla v|^2$.
Using inequality \eqref{bobu81} we get an approximation of order $h$, where the $O(h)$ term is proportional to $\|f^{a,b}-f_h^{a,b}\|_{H^{-1}}$, which depends on $\|\nabla u_1-\nabla u_{1,h}\|_{L^2(\Pp)}$.

For the sake of completeness, the required information for estimating the eigenvalues of the Hessian matrix is provided below. Full computations are made in \cite{Bogosel_Bucur_Polya}. The notations correspond to those in Theorem \ref{thm:eig-hessian}.

\smallskip
\noindent \bo{Estimating $A_k$.} It can be seen that $A_k = a(U_0^1,W^{A_k})$ where $W^{A_k}$ solves
\[ a(W,v) = (f_{\reg}^{A_k},v)_{H^{-1},H_0^1}+(f_{\sing}^{A_k},v)_{H^{-1},H_0^1} ,\ \forall v \in H_0^1(\Pp) , \int_\Pp Wu_1=0\]
with
$$(f_{\reg}^{A_k}, v)_{H^{-1}\times H^1_0}= \sum_{j=0}^{n-1} (\cos (j+1)k\theta+\cos jk\theta) \lb_1 \int_{T_j} u_1v  $$
$$ - \sum_{j=0}^{n-1} \frac{\cos(j+1)k\theta-\cos jk\theta}{\sin \theta}    \int_{T_j} \big (-\sin (2j+1)\theta \partial^2_{xx}u_1 + 2 \cos (2j+1)\theta \partial ^2_{xy} u_1+ \sin(2j+1)\theta \partial ^2_{yy} u_1\big )v $$
$$(f_{\sing}^{A_k}, v)_{H^{-1}\times H^1_0} =  -\sum_{j=0}^{n-1}\int_{S_j} \frac{\cos \theta}{\sin \theta} 2 \cos jk\theta(1-\cos k \theta) \partial_r u_1 v.
$$
An explicit \emph{a priori} estimate for the error between $a(U_0^1,W^{A_k})$ and the finite element counter part follows from the discussion above. 

\smallskip
\noindent\bo{Estimating $B_k$.} It can be seen that $B_k = a(U_0^2,W^{B_k})$ where $W^{B_k}$ solves
\[ a(W,v) = (f_{\reg}^{B_k},v)_{H^{-1},H_0^1}+(f_{\sing}^{B_k},v)_{H^{-1},H_0^1} ,\ \forall v \in H_0^1(\Pp) , \int_\Pp Wu_1=0\]
with
$$(f_{\reg}^{B_k}, v)_{H^{-1}\times H^1_0}= \frac{\cos \theta}{\sin \theta}\sum_{j=0}^{n-1} (\cos (j+1)k\theta-\cos jk\theta) \lb_1 \int_{T_j} u_1v  $$
$$ - \sum_{j=0}^{n-1} \frac{\cos(j+1)k\theta-\cos jk\theta}{\sin \theta}    \int_{T_j} \big (-\cos (2j+1)\theta \partial^2_{xx}u_1 -2 \sin (2j+1)\theta \partial ^2_{xy} u_1+ \cos(2j+1)\theta \partial ^2_{yy} u_1\big )v $$
$$(f_{\sing}^{B_k}, v)_{H^{-1}\times H^1_0} =  \sum_{j=0}^{n-1}\int_{S_j} 2 \cos jk\theta(1-\cos k \theta) \partial_r u_1 v.
$$
An explicit \emph{a priori} estimate for the error between $a(U_0^2,W^{B_k})$ and the finite element counter part follows from the discussion above. Note that the singular part is even with respect to $y$, while the function $U_0^2$ is odd for both the continuous and discrete problems. Therefore the singular part cancels out and can be neglected in the computations in Theorem \ref{thm:estimate-singular} leading to a stronger estimate.

\smallskip
\noindent \bo{Estimating $C_k$.} It can be seen that $C_k = a(U_0^1,W^{C_k})$ where $W^{C_k}$ solves
\[ a(W,v) = (f_{\reg}^{C_k},v)_{H^{-1},H_0^1}+(f_{\sing}^{C_k},v)_{H^{-1},H_0^1} ,\ \forall v \in H_0^1(\Pp) , \int_\Pp Wu_1=0\]
with
$$(f_{\reg}^{C_k}, v)_{H^{-1}\times H^1_0}= \frac{\cos \theta}{\sin \theta}\sum_{j=0}^{n-1} (\sin (j+1)k\theta-\sin jk\theta) \lb_1 \int_{T_j} u_1v  $$
$$ - \sum_{j=0}^{n-1} \frac{\sin(j+1)k\theta-\sin jk\theta}{\sin \theta}    \int_{T_j} \big (-\cos (2j+1)\theta \partial^2_{xx}u_1 -2 \sin (2j+1)\theta \partial ^2_{xy} u_1+ \cos(2j+1)\theta \partial ^2_{yy} u_1\big )v $$
$$(f_{\sing}^{C_k}, v)_{H^{-1}\times H^1_0} = -  \sum_{j=0}^{n-1}\int_{S_j} 2 \sin jk\theta(1-\cos k \theta) \partial_r u_1 v.
$$
An explicit \emph{a priori} estimate for the error between $a(U_0^1,W^{C_k})$ and the finite element counter part follows from the discussion above. Note that the singular part is odd with respect to $y$, while the function $U_0^1$ is even for both the continuous and discrete problems. Therefore the singular part cancels out and can be neglected in the computations in Theorem \ref{thm:estimate-singular} leading to a stronger estimate.

The practical estimates for all quantities involved in the proof of Theorem \ref{thm:estimate-singular} follows from the following:
\begin{itemize}[noitemsep]
	\item Theorem \ref{thm:estimate-H2}, Lemma \ref{lem:general-estimate-rays}
	\item Effective practical estimates computed in \cite[Section 5]{Bogosel_Bucur_Polya}
\end{itemize}

\section{Numerical validation}
\label{sec:num-validation}

This section provides a series of  results which are used in the validated computing process for the  discrete eigenvalue problems and discrete linear systems. Although Intlab \cite{intlab} has functions dedicated to these types of problems, applying the corresponding functions directly leads to inefficient computations overcoming, ultimately, the available computational power. To simplify computational aspects, various theoretical results are shown in the following.

\begin{figure}
	\centering 
	\includegraphics[width=0.5\textwidth]{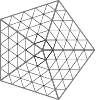}
	\caption{Symmetric mesh for a regular pentagon for $m=5$: the ray $[\bo o \bo a_0]$ is divided into $m$ equal segments.}
	\label{fig:sym-mesh}
\end{figure}

Given a mesh size $h$ and a numerical computation, all quantities in the \emph{a priori} estimates are explicit. In all our computations the mesh considered is symmetric and corresponds to a regular polygon inscribed in the unit circle with one vertex at $\bo a_0 = (1,0)$. See Figure \ref{fig:sym-mesh} for an illustration. The mesh size $h$ is given by $1/m$ where $m$ is the number of segments in which $[\bo o\bo a_0]$ is divided. For a given parameter $m$ the mesh has $N = 1+nm(m+1)/2$ nodes.

First, the eigenvalues, eigenvectors and solutions to linear systems should be enclosed with certified intervals using interval arithmetic in Intlab \cite{intlab}. Let $K,M$ be the rigidity and mass matrices for the $\bo P_1$ finite elements, defined by
\[
K = \left( \int_{\Pp} \nabla \psi_i\cdot \nabla  \psi_j\right)_{1\leq i,j \leq N},
M = \left( \int_{\Pp}  \psi_i  \psi_j\right)_{1\leq i,j \leq N}
,\]
where $\psi_i$ denotes the piece-wise affine function on the mesh of $\Pp$ taking the value $1$ at the node $1\leq i\leq N$ and zero for all other nodes. 

We are interested in solving the generalized eigenvalue problem $Ku=\lambda Mu$, imposing Dirichlet boundary conditions on nodes in $\partial \Pp$. This can be achieved as follows. Denote by $K_0,M_0 \in \Bbb{R}^{d\times d}$ the submatrices of $K$ and $M$ corresponding to the inner nodes. Then the eigenvalue equation corresponding to the Dirichlet boundary condition is simply
\begin{equation}\label{eq:gen-eig}
 K_0 u = \lambda M_0 u.
 \end{equation}
For $h$ small enough, the eigenvalues of \eqref{eq:gen-eig} having small magnitude approximate the eigenvalues of the Dirichlet Laplace operator on $\Pp$. One first objective is to solve \eqref{eq:gen-eig} rigorously, for a given mesh size, including verified enclosures for the eigenvalues $\lambda_{1,h}, \lambda_{2,h}$ and the first eigenvector $u_{1,h}$. In order to keep notations simple, since in this section we always deal with $\lambda_i$ and $u_i$ as solutions of linear systems, we drop the mention of $h$ in the previous notations. The verification process can be achieved using Intlab \cite{intlab} and every element in $K_0, M_0$ should be provided using an explicit interval enclosure by floating point numbers.

The symmetry of the mesh of $\Pp$ illustrated in Figure \ref{fig:sym-mesh} indicates that all triangles in the mesh are similar to an isosceles triangle with the top angle equal to $2\pi/n$. In particular, the matrices $K$, $M$ have precise expressions for all non-zero elements. For interior nodes, not lying on the boundary, we indicate the formulas corresponding to elements in $K_0,M_0$.

\begin{prop}\label{prop:formulas-KM}
	The elements of the rigidity matrix $K_0$ corresponding to nodes not lying on the boundary of $\Pp$ are the following:
	
	- $i$ is the central node: $K_{ii} = n \tan \frac{\theta}{2}, K_{ij} = -\tan \frac{\theta}{2}$
	
	- $i$ is not a central node: $K_{ii}=\frac{2}{\sin \theta}+2\tan \frac{\theta}{2}=\frac{2}{\tan \theta}+4\tan \frac{\theta}{2}$, $K_{ij} = -\frac{1}{\tan \theta} $ (2 times per row/column), $K_{ij}= -\tan \frac{\theta}{2} $ (4 times per row/column)
	
	If $A_h = 0.5h^2\sin \theta$ is the area of one triangle in the mesh, the elements of the mass matrix $M_0$ corresponding to nodes not lying on the boundary of $\Pp$ are
	
	- $i$ is the central node: $M_{ii} = nA_h/6, M_{ij} = A_h/6$
	
	- $i$ is not the central node: $M_{ii} = A_h, M_{ij} = A_h/6$.
\end{prop}

When dealing with equation \eqref{eq:gen-eig} multiple theoretical aspects facilitate the validation process. 

\begin{prop}\label{prop:propertiesKM}
	a) The matrices $K_0, M_0$ are non-singular and eigenpairs for \eqref{eq:gen-eig} correspond to eigenpairs for the matrix $M_0^{-1}K_0$.
	
	b) $K_0$ is an $M$-matrix, i.e. $K_0 = sI-B$ with $B\geq 0$ (all elements of $B$ are non-negative). Therefore $K_0^{-1}$ has non-negative entries. 
	
	c) Furthermore, $K_0^{-1}$ is positive (all entries are strictly positive). Therefore $K_0^{-1}M_0$ is positive and the Perron-Frobenius theorem implies its largest eigenvalue is simple and corresponds to a positive eigenvector. 
	
	d) The smallest eigenvalue of $M_0^{-1}K_0$ is simple and corresponds to a positive eigenvector. All other eigenvectors contain both positive and negative terms.
\end{prop}

\emph{Proof:} a) The result follows at once. Denote by $I$ the set of inner nodes of the mesh and $N_0$ the number of inner nodes. Consider $ x = (x_i)$ the basis coordinates of $\psi = \sum_{i=1}^{N_0} x_i \psi_i$, which is a $\bo P_1$ finite element function verifying $\psi = 0$ on $\partial \Pp$. Therefore, we have
\[  x^T K_0  x = \|\nabla \psi\|_{L^2}\geq 0,  x^T M_0 x = \|\psi\|_{L^2}\geq 0.\]
If one of the previous quantities is zero then obviously $x=0$. Positive definiteness follows.

b) The fact that $K_0$ is an $M$-matrix follows at once from Proposition \ref{prop:formulas-KM}. Furthermore, since each line or column of $K_0$ has at least one non-diagonal entry which is strictly negative, each line and column of $B$ has a non-diagonal entry which strictly positive. 

Let $s>0$ be the maximal element on the diagonal of $K_0$. Then $\frac{1}{s} K_0 = I-\frac{1}{s}B$. The spectral radius of $\frac{1}{s}B$ is smaller than $1$ (see Theorem 8.1.18 from \cite{horn-matrix}). Therefore $(\frac{1}{s} K_0)^{-1} = \sum_{k=0}^\infty (\frac{1}{s}B)^k$, which has non-negative entries since $B \geq 0$.

c) Let $e_i$ be a vector of the canonical basis of $\Bbb R^{N_0}$ and $y = K_0^{-1}e_i\geq 0$. Then
\[ (sI-B)y = e_i \Longrightarrow sy = By+e_i.\]
Since $s>0$, if $y = (y_j)_{j=1}^{N_0}$ it follows that $y_i>0$. Next, if $B_{jk}>0$ then $y_j>0$ if and only if $y_k>0$. To see this observe that
\[ sy_j \geq B_{jk}y_k \text{ and }sy_k \geq B_{jk}y_j.\]
Inductively, since $y_i>0$ all successive neighbors of the node $i$ verify $y_j>0$. Finally, since the mesh is connected, we get $y>0$. Since the result is valid for all elements of the canonical basis, we have $K_0^{-1}>0$.

The matrix $M_0$ is non negative by construction and has non-zero elements on every line and column. Therefore $K_0^{-1}M_0>0$. 

The Perron-Frobenius theorem implies that the largest eigenvalue of $K_0^{-1}M_0$ is simple and corresponds to an eigenvector whose elements are strictly positive. See \cite[Section 8.2]{horn-matrix} for more details.

d) The largest eigenvalue of $K_0^{-1}M_0$ corresponds to the inverse of the smallest eigenvalue of $M_0^{-1}K_0$. Therefore the smallest eigenvalue of \eqref{eq:gen-eig} is simple and corresponds to a positive eigenvector. All other eigenvector verifies $u_jMu_1 = 0$, for all $j>1$. Since $M$ has positive entries on every line and all elements of $u_1$ are positive, each $u_j$, $j>1$ must have both positive and negative elements.  \hfill $\square$

 The following result is useful for the numerical validation of the first three eigenvalues of the regular $n$-gons. Proposition \ref{prop:propertiesKM} already shows that the first discrete eigenvalue $\lambda_{1,h}$ is simple and corresponds to an eigenvector with constant sign. The following result shows that the second eigenvalue is multiple and its cluster is isolated from the rest by the fourth eigenvalue on the unit disk $B_1$. Therefore, for $h$ small enough, if two distinct discrete eigenvalues smaller than $\lambda_4(B_1)$ are found for \eqref{eq:eigenvalue-finite-elements} then these correspond to $\lambda_{1,h}$ and $\lambda_{2,h}$.

\begin{prop}\label{prop:three-eigs}
	a) The second eigenvalue of the regular $n$-gon is at least double. If $\mathcal T_h$ is a symmetric mesh of $\Pp$ like in Figure \ref{fig:sym-mesh} then the discrete eigenvalue $\lambda_{2,h}$ for \eqref{eq:eigenvalue-finite-elements} is multiple.
	
	b) Consider $n \geq 5$. For a mesh size $h>0$ small enough, quantifiable depending on $n$, there exist at most three discrete eigenvalues for \eqref{eq:eigenvalue-finite-elements} in the interval $[0,\lambda_4(B_1)]$, where $B_1$ is the unit disk. For $n \geq 5$, these three eigenvalues correspond to: $\lambda_{1,h}$, which is simple, and the cluster of $\lambda_{2,h}$ with multiplicity $2$. 
\end{prop}
Note that on the regular polygon the second eigenvalue is precisely double. This could be justified with some extra arguments, but   this is not necessary for our further analysis.

\emph{Proof:} a) The regular $n$-gon is connected, therefore  the first eigenvalue $\lambda_1$ is simple and corresponds to an eigenfunction $u_1$ having constant sign. Consider the second eigenvalue $\lambda_2$ and an associated eigenfunction $u_2$. It is known \cite{Al94} that the nodal line corresponding to $u_2$ touches the boundary at exactly two points. 

Assume $\lambda_2$ is simple. Then any rotation $\tilde u_2$ of $u_2$ with a multiple of $2\pi/n$ around the origin gives an eigenfunction for $\lambda_2$. Since the nodal line of $u_2$ touches the boundary at exactly two points, it is impossible that $\tilde u_2$ is a multiple of $u_2$. Therefore the multiplicity of $\lambda_2$ is at least $2$. 

Assume the regular $n$-gon is meshed with a symmetric mesh $\mathcal T_h$ like in Figure \ref{fig:sym-mesh}. Consider $\lambda_{2,h}$ a second eigenvalue for the discrete problem \eqref{eq:eigenvalue-finite-elements} together with a second eigenvector $u_{2,h}$. Since $u_{2,h}$ is orthogonal on $u_{1,h}$, which has constant sign, there exist positive and negative elements in $u_{2,h}$. Consider a permutation of the elements of $u_{2,h}$ corresponding to a rotation of angle multiple of $2\pi/n$ on the mesh. The positivity and negativity of elements of $u_{2,h}$ cannot remain invariant for all such rotations. Therefore the discrete eigenvalue $\lambda_{2,h}$ is also multiple. Passing to the limit when $h\to 0$ this also gives an alternative proof for the multiplicity of $\lambda_2(\Pp)$.

b) Since the regular $n$-gon $\Pp$ is inscribed in the unit ball, we have $\Pp \subset B_1$, therefore $\lambda_k(\Pp) \geq \lambda_k(B_1)$ for any $k \geq 1$. In \cite{Nitsch-reg-poly} it is shown that $\lambda_1(\Pp)$ is decreasing with respect to $n$. Such a result is not available for higher eigenvalues. Whenever one can show that $\lambda_2(\Pp)=\lambda_3(\Pp)<\lambda_4(B_1)$, the result follows. For $n=4$ the second eigenvalue of the square is $\lambda_2(\Bbb P_4) = 2.5\pi^2<\lambda_4(B)=j_{2,1}^2 \approx 5.13...^2 = 26.31...$, where $j_{2,1}$ is the first zero of the Bessel function of the second kind. Therefore, for all $n \geq 4$ one has $\lambda_2(\Pp)=\lambda_3(\Pp)<\lambda_4(B_1)$. For $h$ small enough the discrete eigenvalues verify the same inequality. This will be checked \emph{a posteriori} in all our computations. \hfill $\square$

In \cite{MayerEigs} multiple results regarding validated computations regarding eigenvalues and eigenvectors are presented. Theorem 16 in the previous reference shows that in the particular case of symmetric matrices it is possible to obtain bounds for all eigenvalues and eigenvectors associated to simple eigenvalues in terms of residuals. A modification of this result corresponding to generalized eigenvalues is given below, notably using the natural matrix norm associated to the mass matrix. In some situations, this allows to avoid the usage of eigenvalue and eigenvector verification routines (like \texttt{verifyeig} in Intlab \cite{intlab}) which are costly for large matrices. 

\begin{prop}\label{prop:validation-eigs}
	Let $K,M$ be symmetric matrices in $\Bbb{R}^{N\times N}$ such that $K, M$ are positive definite. Let $\tilde x \in \Bbb{R}^N$ and $\tilde \lambda \in \Bbb{R}$ be an approximate eigenpair. Denote the residual vector $r:= K\tilde x-\tilde \lambda M\tilde x$ and $|x|_M^2 = x^TMx$. Then the following assertions hold:
	
	a) The interval $\left[ \tilde \lambda - \sqrt{\frac{r^TM^{-1}r}{\tilde x^TM\tilde x}},\tilde \lambda +\sqrt{\frac{r^TM^{-1}r}{\tilde x^TM\tilde x}} \right]$ contains at least one generalized eigenvalue of $(K,M)$.
	
	b) If $a\geq r$, $|\tilde x|_M=1$ and if $[\tilde \lambda -a, \tilde \lambda +a]$ contains a simple eigenvalue $\lambda^*$ but no more eigenvalues of $(K,M)$ then there is an eigenvector $x^*$ with associated with $\lambda^*$ verifying 
	\[ |\tilde x-x^*|_M^2 \leq \frac{r^TM^{-1}r}{a^2}+\max\left\{(|\tilde x|_M-1)^2, \left(1-\sqrt{\tilde x^T Mx - \frac{r^TM^{-1}r}{a^2}} \right)^2\right\},\]
	provided the square root exists. 
\end{prop}

\emph{Proof:} The proof could follow from \cite[Theorem 16]{MayerEigs} by a change of variables, choosing $A = M^{-1/2}KM^{-1/2}$, but for the sake of completeness the following simple arguments are given below.

a) Since generalized eigenvalues of $(K,M)$ correspond to eigenvalues of $M^{-1}K$, which is symmetric, positive definite, there exists a basis of $\Bbb{R}^N$ made of generalized eigenvectors $(v_k)$ of $(K,M)$ associated to eigenvalues $\lambda_k>0$, $k=1,...,N$. The vectors $v_k$ are assumed orthonormal with respect to the matrix $M$, in other words, $v_i^TMv_j = \delta_{ij}$. Therefore, if $\tilde x = \sum_{i=1}^N \alpha_k v_k$ then $ \tilde x^T M \tilde x = \sum_{k=1}^N \alpha_k^2$. On the other hand 
\[ r = K \tilde x -\tilde \lambda M \tilde x = \sum_{k=1}^N (\lambda_k-\tilde \lambda) \alpha_k Mv_k.\]
This implies that
\begin{equation}\label{eq:discrete-eig-bound} r^TM^{-1}r = \sum_{i=1}^N (\lambda_k-\tilde \lambda)^2 \alpha_k^2 \geq \left( \min_{k=1,...,N} |\lambda_k-\tilde \lambda|^2\right) \tilde x^TM\tilde x.
\end{equation}
This proves the first assertion. 

b) Assume $\lambda^*=\lambda_k$ is simple with eigenvector $v_k$ the normalization conventions described previously. Using the decomposition $\tilde x = \sum_{i=1}^N \alpha_k v_k$ we find 
\[ \tilde x^T M \tilde x = \sum_{j=1}^N \alpha_k^2.\]
The hypothesis on the separation of $\tilde \lambda$ with respect to $\lambda_j, j\neq k$ gives
\[ r^TM^{-1}r=(\lambda_k-\tilde \lambda)^2\alpha_k^2+\sum_{j \neq k}(\lambda_j-\tilde \lambda)^2\alpha_j^2 \geq a^2 (\tilde x^T M \tilde x-\alpha_k^2).\]
On the other hand, the previous inequality gives
\[ |\tilde x-v_k|_M^2 = \tilde x^T M \tilde x-2\alpha_k+1\leq \frac{r^T M^{-1} r}{a^2}+(\alpha_k-1)^2.\]

Without loss of generality, assume $\alpha_k\geq 0$. If $\alpha_k\geq 1$ then $1\leq \alpha_k\leq \sqrt{\tilde x^T M \tilde x}$, therefore
\[ |\tilde x-v_k|_M^2 \leq \frac{r^TM^{-1}r}{a^2}+(|\tilde x|_M-1)^2 \]
For $0 \leq \alpha_k<1$ we have $\alpha_k \geq \sqrt{|\tilde x|_M^2-\frac{r^TM^{-1}r}{a^2}}$, provided $|\tilde x|_M^2-\frac{r^TM^{-1}r}{a^2}\geq 0$. In this case
\[ |\tilde x-v_k|_M^2 \leq \frac{r^TM^{-1}r}{a^2}+\left(1-\sqrt{\tilde x^T Mx - \frac{r^TM^{-1}r}{a^2}} \right)^2. \]
The conclusion follows. \hfill $\square$

\begin{rem}\label{rem:eigsM}
	It can be observed that the technique used to prove the second point in Proposition \ref{prop:validation-eigs} gives an estimate concerning the distance to the eigenspace associated to a cluster in case the eigenvalue is multiple. In this case the residual together with the minimal distance from the cluster to the closest neighboring eigenvalues is enough to give the estimate.
	
	The expressions involving the matrix norm $|\cdot|_M$ can be evaluated either explicitly by solving a validated linear system or precise properties of the matrix $M$ (which is very particular in our case: see Proposition \ref{prop:formulas-KM}) can help to turn these estimates into a more explicit form.
	
	According to \cite{cond_numbers} the eigenvalues of the mass matrix verify
	\[ \frac{1}{2}\min M_{ii} \leq \lambda_{\min}(M) \leq \lambda_{\max}(M) \leq 2\max M_{ii}.\]
	The explicit formulas for $M_{ii}$ give
	\[ \frac{1}{12}\min\{n,6\}A_h \leq \lambda_{\min}(M) \leq \lambda_{\max}(M) \leq \frac{1}{3}\max \{n,6\}A_h,\]
	where $A_h = 0.5h^2 \sin \theta$ is the area of a small triangle in the mesh shown in Figure \ref{fig:sym-mesh}. Therefore
	\[ r^TM^{-1}r \leq \frac{1}{\lambda_{\min}(M)}|r|^2\leq  \frac{12}{\min\{n,6\}A_h}|r|^2, \ \ |x|^2\leq \frac{1}{\lambda_{\min}(M)}|x|_M^2 \leq \frac{12}{\min\{n,6\}A_h}|x|_M^2.\] 
	Similar inequalities remain valid for $(K_0,M_0)$ which are submatrices of $(K,M)$. Classical interlacing inequalities show that eigenvalues of $K_0$, $M_0$ verify the same bounds.
\end{rem}
	
A second key point is the approximation of the material derivatives defined in \eqref{eq:material-U0} which imply solving the discrete systems associated to the finite element formulations. For simplicity, the index $h$ is dropped in the following. It should be noted that in the remaining of this section, only the discrete problem is considered. In this case, given the (discrete) eigenfunction $u_1$ associated to $\lambda_1$ the system
\begin{equation}\label{eq:saddle} \begin{pmatrix}
A& b \\
b^T  & 0 
\end{pmatrix} 
\begin{pmatrix}
U \\ 0
\end{pmatrix}
= \begin{pmatrix}
f\\
0
\end{pmatrix}
\end{equation}
needs to be solved for
\[ A = K_0-\lambda M_0 \in \Bbb{R}^{d\times d} , b = M_0 u \in \Bbb{R}^d.\]
The compatibility condition $f^T u_1 = 0$ needs to be verified by $f \in \Bbb{R}^{d}$, where $u_1$ is the eigenvector associated to the smallest eigenvalue $\lambda_1$. 

Recall that $T_\pm$ are the two triangles delimited by $y=0$ and $y = \pm x\tan \frac{2\pi}{n}$. Next, $\varphi$ is the $\bo P_1$ function defined on the triangulation of $\Pp$ with $n$ triangles where one node is at the center (see Figure \ref{fig:simple-triangulations}) taking value $1$ only at $(1,0)$. We have
\[ \nabla \varphi = \left(1,-\frac{1}{\tan \theta}\right)\chi_{T_+} + \left(1,\frac{1}{\tan \theta}\right)\chi_{T_-},\]
where $\theta = 2\pi/n$. 
The variational formulation for \eqref{eq:U0} reads: $\bo U \in H_0^1(\Pp)^2$ verifying
\[ a(\bo U,v) = \int_\Pp \left[ (\nabla u_1 \cdot \nabla \varphi)\nabla v +(\nabla \varphi \cdot \nabla v)\nabla u_1 - \begin{pmatrix}
2\lambda_1/n \\ 0 
\end{pmatrix}u_1v. \right] \]
The explicit equations verified $\bo U = (U_1,U_2)$ defined in \eqref{eq:material-U0} are as follows: $U_i \in H_0^1(\Pp)$

\begin{equation}\label{eq:U1}
\int_\Pp \nabla U_1 \cdot \nabla v -\lambda_1 U_1 v = \int_{T_+\cup T_-} \left( 2\partial_x u \partial_x v +\partial_y \varphi (\partial_x u \partial_y v+\partial_y u \partial_x v)\right) -2\frac{\lambda_1}{n}\int_\Pp uv  
\end{equation}

\begin{equation}\label{eq:U2}
\int_\Pp \nabla U_2 \cdot \nabla v -\lambda_1 U_1 v = \int_{T_+\cup T_-} \left( 2 \partial_y \varphi \partial_y u \partial_y v + (\partial_x u \partial_y v+\partial_y u \partial_x v) \right) 
\end{equation}
for every $v \in H_0^1(\Pp)$, with constraints $\int_\Pp U_i u_1 = 0$.

Consider the following partial rigidity matrices:
\[ K_{xx} = \left(\int_{T_+\cup T_-} \partial_x \psi_i \partial_x \psi_j \right)_{i,j=1}^{N_0}, K_{yy} = \left(\int_{T_+} \partial_y \psi_i \partial_y \psi_j - \int_{T_-} \partial_y \psi_i \partial_y \psi_j\right)_{i,j=1}^{N_0}\]
\[ K_{xy}^\pm = \left(\int_{T_\pm} \partial_x \psi_i \partial_y \psi_j+\partial_y \psi_i \partial_x \psi_j \right)_{i,j=1}^{N_0}
\]
Note that these matrices are only considered for interior vertices of the mesh, since Dirichlet boundary conditions are taken on $\partial \Pp$. For the triangulation illustrated in Figure \ref{fig:sym-mesh} the elements of the matrices above are all explicit in terms of $\theta=\frac{2\pi}{n}$ and the usual trigonometric functions. The terms can be found using the following result.
\begin{prop}\label{prop:elements-xy}
	For a triangle which is homothetic to $T_+$ with factor $h$, with vertices labeled $i,j,k$ in the order of $(0,0),(1,0),(\cos \theta, \sin \theta)$, $\theta=\frac{2\pi}{n}$, the gradients of the associated $\bo P_1$ finite element functions $\psi_i,\psi_j,\psi_k$ are given by
	\[ \begin{array}{lll}\partial_x \psi_i = -\frac{1}{h}, & \partial_x \psi_j = \frac{1}{h}, & \partial_x \psi_k = 0 \\
	\partial_y \psi_i = -\frac{1}{h}\tan \frac{\theta}{2}, & \partial_y \psi_j = -\frac{1}{h\tan \theta}, & \partial_y \psi_k = \frac{1}{h\sin \theta}.
	 \end{array}\]
	These formulas coupled with the fact that all elements in the symmetric mesh have area $\frac{1}{2}h^2\sin \theta$ show that the elements of $K_{xx}, K_{yy}, K_{xy}, K_{yx}$ are completely determined in function of $\theta$ and do not depend on $h$. 
\end{prop}

Then the discrete systems associated to \eqref{eq:U1}, \eqref{eq:U2} when using $\bo P_1$ finite elements have the form
\begin{equation}\label{eq:U1-discrete}
 (K_0-\lambda_1 M_0)U_1 = 2K_{xx} u_1 -\frac{1}{\tan \theta} (K_{xy}^+-K_{xy}^-)u_1-\frac{2\lambda_1}{n} M_0 u_1
\end{equation}
\begin{equation}\label{eq:U2-discrete}
(K_0-\lambda_1 M_0)U_1 = -\frac{2}{\tan \theta}K_{yy} u_1+(K_{xy}^+ + K_{xy}^-)u_1,
\end{equation}
verifying the following orthogonality relations 
\begin{equation}\label{eq:Ui-discrete-orthogonal}
u_1 M_0 U_1 = u_1 M_0 U_2 = 0.
\end{equation}

More precisely, the system of the form \eqref{eq:saddle} needs to be solved for
\[ A = K_0-\lambda_1 M_0, b = M_0u_1,\]
where $u_1$ is the first generalized eigenvector associated to $(K_0,M_0)$ and the right hand sides are computed using \eqref{eq:U1-discrete}, \eqref{eq:U2-discrete}. 

According to \cite{Saddle-Kimura-Chen}, in order to precondition the saddle point system it is enough to consider a perturbation of the form 
\begin{equation}\label{eq:mw} A(w) = A+wbb^T	\end{equation}
and 
\begin{equation}\label{eq:bigM}
\mathcal A(w) = \begin{pmatrix}
A(w) & 0 \\
0 & b^T A(w)^{-1} b
\end{pmatrix}.\end{equation}
Then according to \cite{Saddle-Kimura-Chen} and \cite{Golub2005} we have the following result:

\begin{lemma}\label{lem:saddle}
	a) Denoting $\mathcal H = \begin{pmatrix}
	A&b \\ b^T& 0
	\end{pmatrix}$, if $w>0$, $A$ is positive semi-definite and $\mathcal H$ is non-singular then $A(w)$ given in \eqref{eq:mw} is positive definite.
	
	b) The eigenvalues of $\mathcal A(w)^{-1} \mathcal H$ belong to 
	\[ \left[ -1 ,\frac{1-\sqrt{5}}{2}\right] \cup \left[ 1,\frac{1+\sqrt{5}}{2}\right].\]
	
	c) If $U^*$ is the exact solution of $\mathcal H U^* = \overline f=\begin{pmatrix}
	f\\0
	\end{pmatrix}$ then
	\begin{equation}\label{eq:estimate-saddle}
	\|U^*-U\|_2 \leq \frac{2}{\sqrt{5}-1} \max\{\|A(w)^{-1}\|_2, \|A(w)\|_2 \|b^Tb\|_2^{-1}\}\|\mathcal AU-\overline f\|_2.
	\end{equation} 
\end{lemma}

\emph{Proof:} a) Let $x \in \Bbb{R}^d$.  Then $x^T A(w) x= x^TAx + w(b^T x)^2\geq 0$, since $A$ is positive semi-definite. 

Moreover, $x^T A(w) x=0$ implies $x^TAx=0$ and $b^Tx=0$. Therefore If $\bar x = (x,0)$ then $\mathcal H \bar x=0$ and since $\mathcal H$ is non-singular, $x=0$. Threrfore $A(w)$ is positive definite. 

b) This is proved in \cite[Theorem 2.5]{Golub2005}.

c) This is proved in \cite[Theorem 2.2]{Saddle-Kimura-Chen}.

\hfill $\square$ 

In Lemma \ref{lem:saddle} we have are two degrees of freedom regarding the choice of parameters. The first one is the scalar parameter $w$. Secondly, notice that multiplying the vector $b$ with a scalar does not change solutions to the saddle point system \eqref{eq:saddle}. In the following $b$ is chosen of the form $b = \gamma M_0 u_1$.

In the following, the dimension $d$ of the system \eqref{eq:saddle} is $N_0$, the number of nodes in the mesh, not lying on the boundary $\partial \Pp$. Consider the set of generalized eigenpairs $(\lambda_i, u_i)$, $i=1,...,d$ for \eqref{eq:gen-eig} which forms a basis of $\Bbb{R}^d$, orthonormalized with respect to $M_0$. Assume $\lambda_1\leq ...\leq \lambda_d$. For a general $x \in \Bbb{R}^d$ let $x = \sum_{i=1}^d \alpha_i v_i$ its decomposition in this basis. 

Then, for $b = \gamma M_0u_1$ we have
\[ x^T A(w) x= x^T(K_0-\lambda_1 M_0)x +\gamma^2 w(u_1^TM_0x)^2 = \sum_{i=1}^d \alpha_i^2 (\lambda_i-\lambda_1) +\gamma^2 w\alpha_1^2,\]
since $b^T x = \gamma \sum_{i=1}^n \alpha_i u_i^T M_0 u_1 = \gamma\alpha_1$. Therefore $x^T A(w)x \geq \min\{\lambda_2-\lambda_1,\gamma^2 w\}$. Since $A(w)$ is symmetric, $\lambda_1$ is simple (Proposition \ref{prop:propertiesKM}) and $w>0$ we find
\[ \|A(w)^{-1}\|_2 \leq \max\{\frac{1}{\lambda_2-\lambda_1},\frac{1}{\gamma^2 w}\}.\]

Continuing the estimates, we have $\|A(w)\|_2\|b^Tb\|_2^{-1} \leq (\|A\|_2+w\|b\|^2)/\|b\|^2\leq \|A\|_\infty/\|b\|_2^2+w$. Replacing $b = \gamma M_0u_1$ gives
\[ \|A(w)\|_2\|b^Tb\|_2^{-1} \leq \frac{\|K_0-\lambda_1M_0\|_\infty}{\gamma^2 \|M_0u_1\|_2^2} +w.\]
The matrix norm $\|K_0-\lambda_1M_0\|_\infty$ is easily evaluated using Intlab. Thus, once $w$ and $\gamma$ are chosen, estimate \eqref{eq:estimate-saddle} is explicit once the residual $\|\mathcal AU-f\|_2$ is known.

The inequalities above show that:
\begin{itemize}
	\item if $\gamma^2w \geq \lambda_2-\lambda_1$ then $\|A(w)^{-1}\|_2 \leq 1/(\lambda_2-\lambda_1)$.
	\item choosing $w = (\lambda_2-\lambda_1)/\gamma^2$ gives 
	\[ \|A(w)\|_2\|b^Tb\|_2^{-1} \leq \frac{\|K_0-\lambda_1M_0\|_\infty}{\gamma^2 \|M_0u_1\|_2^2} +\frac{\lambda_2-\lambda_1}{\gamma^2},\]
	therefore, when $\gamma$ increases the estimate gets stronger.
	\item Nevertheless, taking $\gamma$ large will increase the interval evaluating how well the floating point approximation verifies the normalization constraint $u_1 M_0 U = 0$. Choosing $\gamma$ too large will deteriorate the global residual estimate. In practical computations, choosing $\gamma = \gamma_0/\|M_0u_1\|_2$ for $\gamma_0 \in [1,10]$ gives satisfying results. 
\end{itemize}

Given $\mathcal H$ the explicit matrix from \eqref{eq:saddle}, an interval enclosure $\overline{\mathcal H}$ is considered. If enclosures for the discrete eigenvector $u_1$ are available, then an enclosure $\overline f$ for the right hand side of \eqref{eq:U1-discrete}, \eqref{eq:U2-discrete} are computed. Given $U$ a floating point solution found for solving \eqref{eq:saddle} for the system given by the midpoint of $\overline{\mathcal H}$ (through an iterative approach, Conjugate Gradient for example) with right hand side the midpoint of $\overline f$, the residual $\|\overline{\mathcal H} U -\overline f\|_2$ is evaluated using interval arithmetic. Thus, a validated upper bound for $\|{\mathcal H} U -f\|_2$ is found and an explicit bound for the distance $\|U-U^*\|_2$ is obtained.



\section{Validation strategy, implementation and results}
\label{sec:validation}

\subsection{Detailed description of the validation procedure}
In this section, the validation strategy is presented, based on results proved in previous sections. All computations are performed using interval arithmetic in Intlab \cite{intlab}.

\smallskip
\noindent \bo{Inputs.} 

-  $n$, the number of vertices of the regular polygon 

-  $m$ the number of segments in which the segment $[0,1]$ is divided in the mesh shown in Figure \ref{fig:sym-mesh}. This gives rise to the mesh size parameter $h=1/m$.

\begin{figure}
	\centering 
	\includegraphics[height=0.35\textwidth]{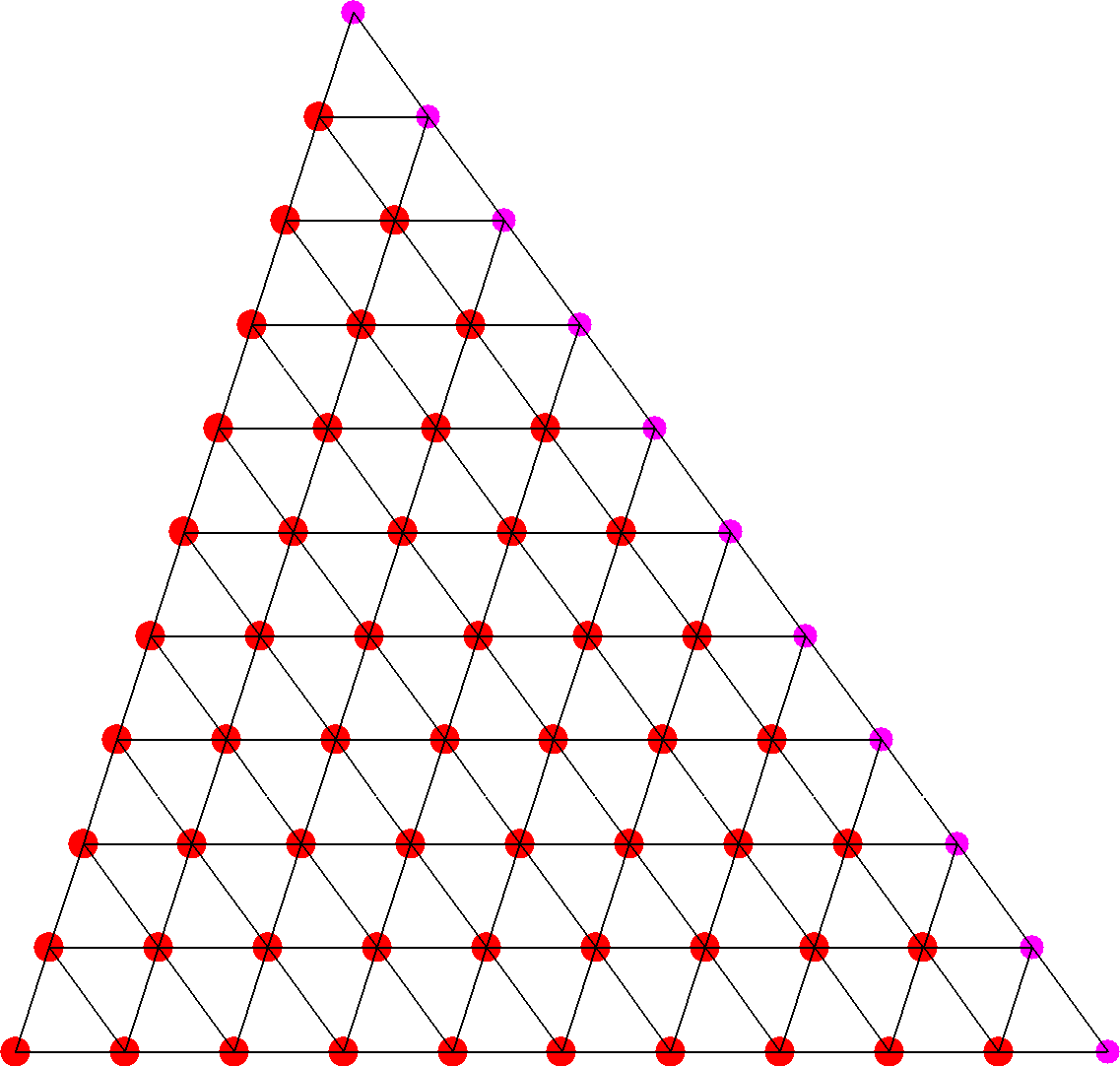}\quad
	\includegraphics[height=0.35\textwidth]{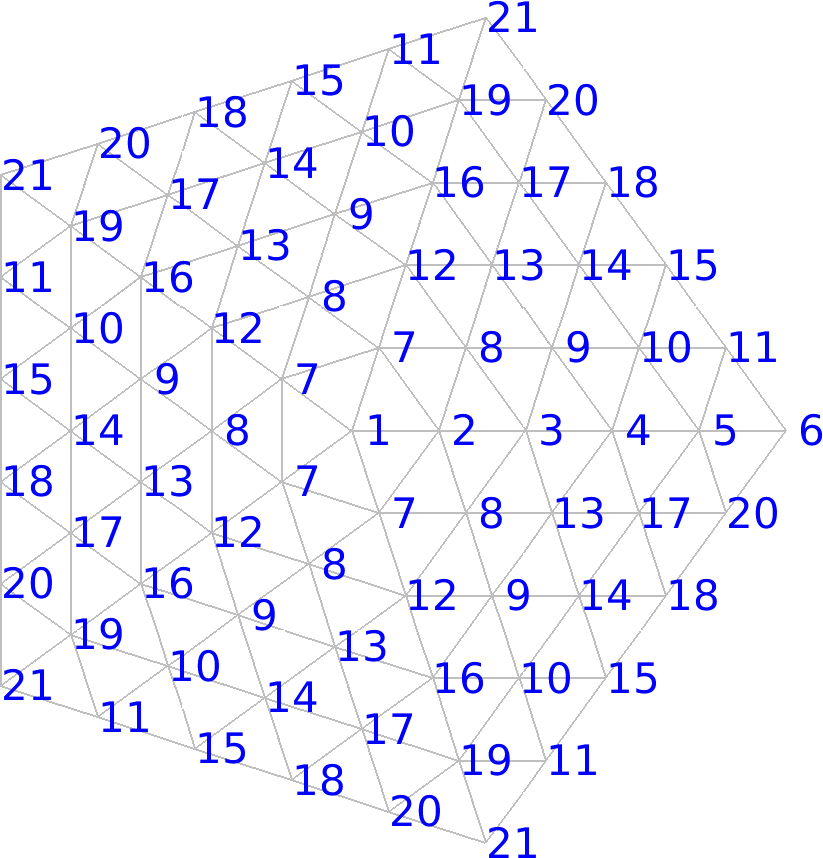}
	\caption{Mesh of a triangular slice $\mathcal T_h^0$ of the regular $n$-gon (left). Mapping from the slice $\mathcal T_h^0$ to the whole mesh used to extend a function by dihedral symmetry (right).}
	\label{fig:mesh-slice-extension}
\end{figure}

The regular polygon $\Pp$ is meshed with a mesh $\mathcal T_h$ like in Figure \ref{fig:sym-mesh}, in particular the rays $S_i$ are meshed exactly. The mesh is constructed manually and the submesh $\mathcal T_h^0$ corresponding to points in $\{x\geq 0, y \in [0,x\tan \theta ]\}$ is considered. An explicit mapping allows to extend finite element functions on $\mathcal T_h^0$ to $\mathcal T_h$ respecting the dihedral symmetry. See Figure \ref{fig:mesh-slice-extension} for an example.

\medskip
\noindent \bo{A. Finite element computations.}

\smallskip
\noindent \bo{a) Rigidity and Mass matrices.} Following Proposition \ref{prop:formulas-KM} the elements of the rigidity and mass matrices $K_0,M_0$ associated to the interior nodes have explicit formulas in terms of $\theta = 2\pi/n$. The matrices are constructed using these expressions in Intlab. Thus $K_0,M_0$ are interval matrices containing the exact rigidity and mass matrices. Since $K_0$, $M_0$ are submatrices associated to the interior nodes, this automatically imposes the Dirichlet boundary conditions. 

\smallskip
\noindent \bo{b) The first discrete eigenvalue and eigenfunction.} The first two eigenvalues for the discrete eigenproblem \eqref{eq:eigenvalue-finite-elements}, $K_0u_{i,h}=\lambda_{i,h}M_0u_{i,h}$ are computed initially using floating point arithmetic. Proposition \ref{prop:validation-eigs} gives two interval enclosures for $\lambda_{1,h}$ and $\lambda_{2,h}$ in terms of the residuals. Following Proposition \ref{prop:three-eigs} if these two intervals together with the \emph{a priori} estimates \eqref{eq:difflam} give upper bounds smaller than the fourth eigenvalue of the unit ball $\lambda_4(B_1)$ then the two enclosures correspond indeed to the first two discrete eigenvalues of $(K_0,M_0)$. 

For the first eigenvalue $\lambda_{1,h}$, which is simple (see Proposition \ref{prop:propertiesKM}) an enclosure is found for the first eigenvector $u_{1,h}$. This is achieved either using a modification of the routine \texttt{verifyeig} in Intlab \cite{intlab} (where the matrix inversion is replaced with three validated linear systems solved with \texttt{verifylss}) applied on the slice $\mathcal T_h^0$ or by using the direct bound using the residual, shown in Proposition \ref{prop:validation-eigs}. Whenever Intlab can be used for bounding the first eigenvector it gives a tighter enclosure than the residual estimate. 

\begin{figure}
	\centering
	\includegraphics[width=0.4\textwidth]{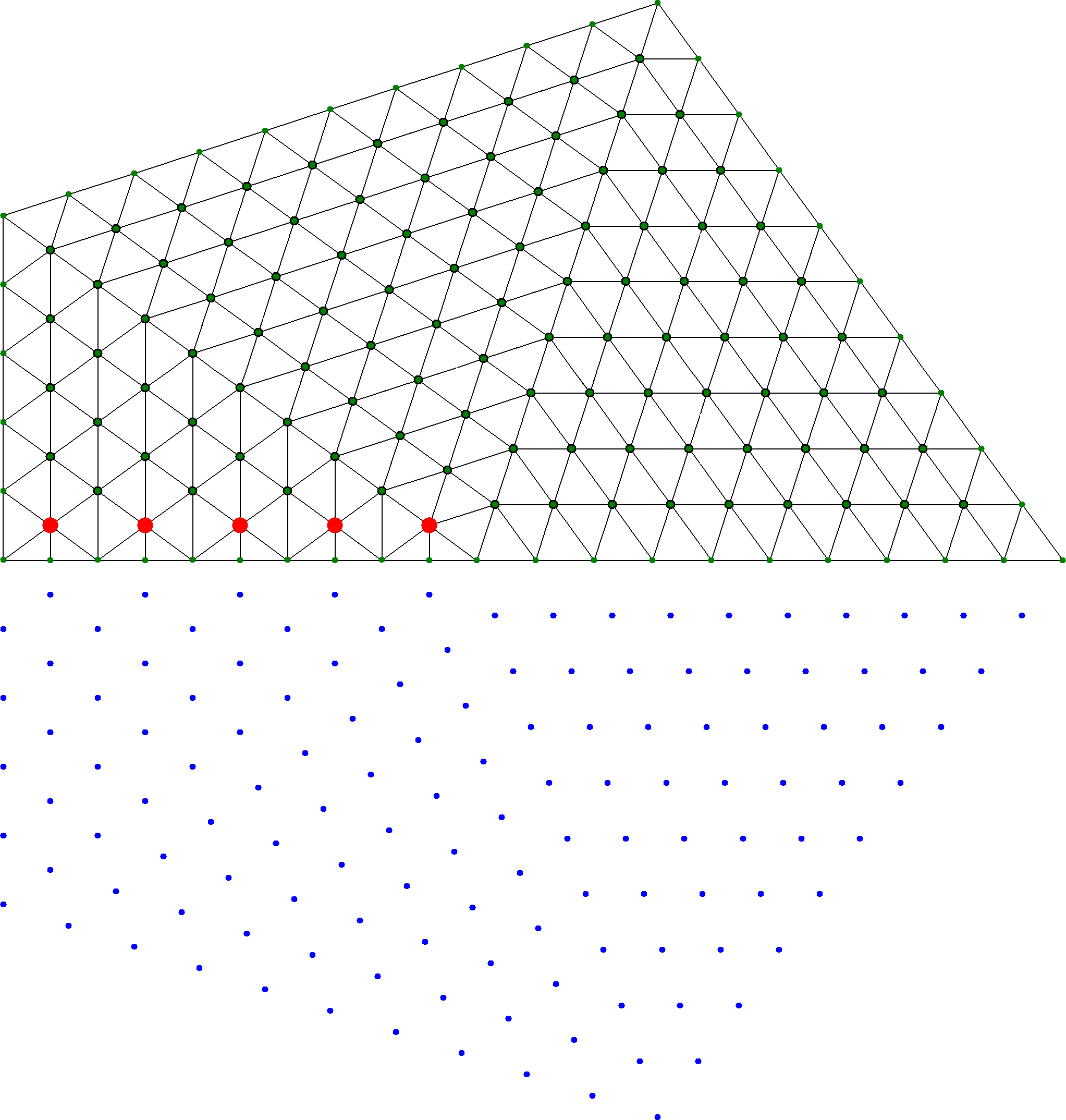}
	\caption{Using half the mesh of $\Pp$ to solve the discrete problem associated to the second eigenvalue in \eqref{eq:eigenvalue-finite-elements}. The associated mass and rigidity matrices are submatrices of those for the full polygon. The diagonal elements of the points underlined need to be modified slightly.}
	\label{fig:half-mesh}
\end{figure}

The second eigenvalue is always validated using the residual estimate given in Proposition \ref{prop:validation-eigs}. For the second eigenvalue the whole mesh of $\Pp$ is be used. For large meshes, the existence of symmetric eigenfunctions for $\lambda_{2,h}$ allows to use only half the mesh for finding an interval enclosure for the second discrete eigenvalue. The associated matrices can be obtained taking submatrices corresponding to the full mesh and replacing the diagonal values corresponding to $m/2$ nodes underlined in Figure \ref{fig:half-mesh}.

\smallskip
\noindent \bo{c) Approximation of the material derivatives.} Next, all information required for solving \eqref{eq:U1-discrete}, \eqref{eq:U2-discrete} is computed using interval arithmetic. In particular the partial rigidity matrices described in Proposition \ref{prop:elements-xy} are computed using exact formulas in Intlab.

Systems \eqref{eq:U1-discrete}, \eqref{eq:U2-discrete} are solved initially in floating point arithmetic using the Conjugate Gradient algorithm. An initialization with the zero vector guarantees that the orthogonality condition on $u_{1,h}$ is preserved since the right hand sides are orthogonal to $u_{1,h}$. 

Following Lemma \ref{lem:saddle} the residuals for the linear systems are computed using interval arithmetic for the floating point results obtained previously. All variables appearing in the estimate \eqref{eq:estimate-saddle} are computed in floating point arithmetic. Notably, the interval enclosure of $\lambda_{2,h}-\lambda_{1,h}$ obtained in the previous step is used. 

In the end, interval enclosures for solutions of \eqref{eq:U1-discrete} and \eqref{eq:U2-discrete}, verifying the orthogonality conditions \eqref{eq:Ui-discrete-orthogonal}, are obtained.

\medskip
\noindent \bo{B. Computation of the eigenvalues of the Hessian matrix.} Up to this point interval enclosures are available for $\lambda_{1,h}, \lambda_{2,h}, u_{1,h}, U_h^1, U_h^2$. These are used to compute the eigenvalues of the Hessian matrix of the first eigenvalue, given in Theorem \ref{thm:eig-hessian}. 

The computation of the integrals of partial derivatives of $u_{1,h}$ and $U_h^i$ is achieved by assembling partial rigidity matrices for all triangles $T_i$ corresponding to slices of the regular $n$-gon into $n$ parts like in Figure \ref{fig:simple-triangulations}. These matrices depend again explicitly on $\theta = 2\pi/n$.

At the end of this step, the discrete approximations $\mu_{j}^h$, $j=0,...,2n-1$ for the eigenvalues of the Hessian are available with guaranteed interval enclosures. 

\medskip
\noindent \bo{C. A priori estimates for finite elements.} The final step is controlling the distances $|\mu_j-\mu_j^h|$, $j=0,...,2n-1$ using the \emph{a priori} error estimates for finite element computations:
\begin{itemize}
	\item $|\lambda_k-\lambda_{k,h}|$ is bounded by \eqref{eq:difflam} according to \cite{LiOi13}.
	\item $\|\nabla u_1-\nabla u_{1,h}\|_{L^2}$ and $\|u_1-u_{1,h}\|_{L^2}$ are bounded by \eqref{eq:grad-u-estimate}, \eqref{eq:L2-u-estimate}.
	\item In some cases, the equality \eqref{eq:relation-errors} may be used to iteratively replace the weakest of the three estimates above with a better one, therefore improving all estimates.
\end{itemize}

Finally, results in Section \ref{sec:material} are used to quantify the finite element errors for the solutions of the material derivatives \eqref{eq:U1}, \eqref{eq:U2} and for equations described in Section \ref{sec:hess-eig-estimates} involved in the estimation of the errors for the Hessian eigenvalues. 

\medskip
\noindent \bo{D. Certification of local minimality.}
At the end of the previous step a guaranteed error estimate is found for each one of the Hessian eigenvalues in Theorem \ref{thm:eig-hessian}. According to \cite[Section 4]{Bogosel_Bucur_Polya} four of the eigenvalues of the Hessian on the regular $n$-gon are identically equal to zero. Thus, if $2n-4$ of the interval enclosures found for the Hessian eigenvalues are included in $(0,+\infty)$, then the numerical certification of local minimality succeeds: the regular $n$-gon is a local minimizer for the first Dirichlet-Laplace eigenvalue. 

\subsection{Simulations and Results} Given $n$ the number of vertices and $m$ the number of segments on the ray $[\bo o \bo a_0]$ in the mesh illustrated in Figure \ref{fig:sym-mesh}, the computation described in the previous section is performed. The validation procedure succeeds for $n\in \{5,6\}$. The case $n\geq 7$ is inaccessible for the moment. It is not a question of computational power: the \emph{a priori} error estimates are not small enough compared to the size of the intervals associated to the variables appearing in the computation. Details are shown below.

\begin{figure}
	\centering
	\includegraphics[height=0.39\textwidth]{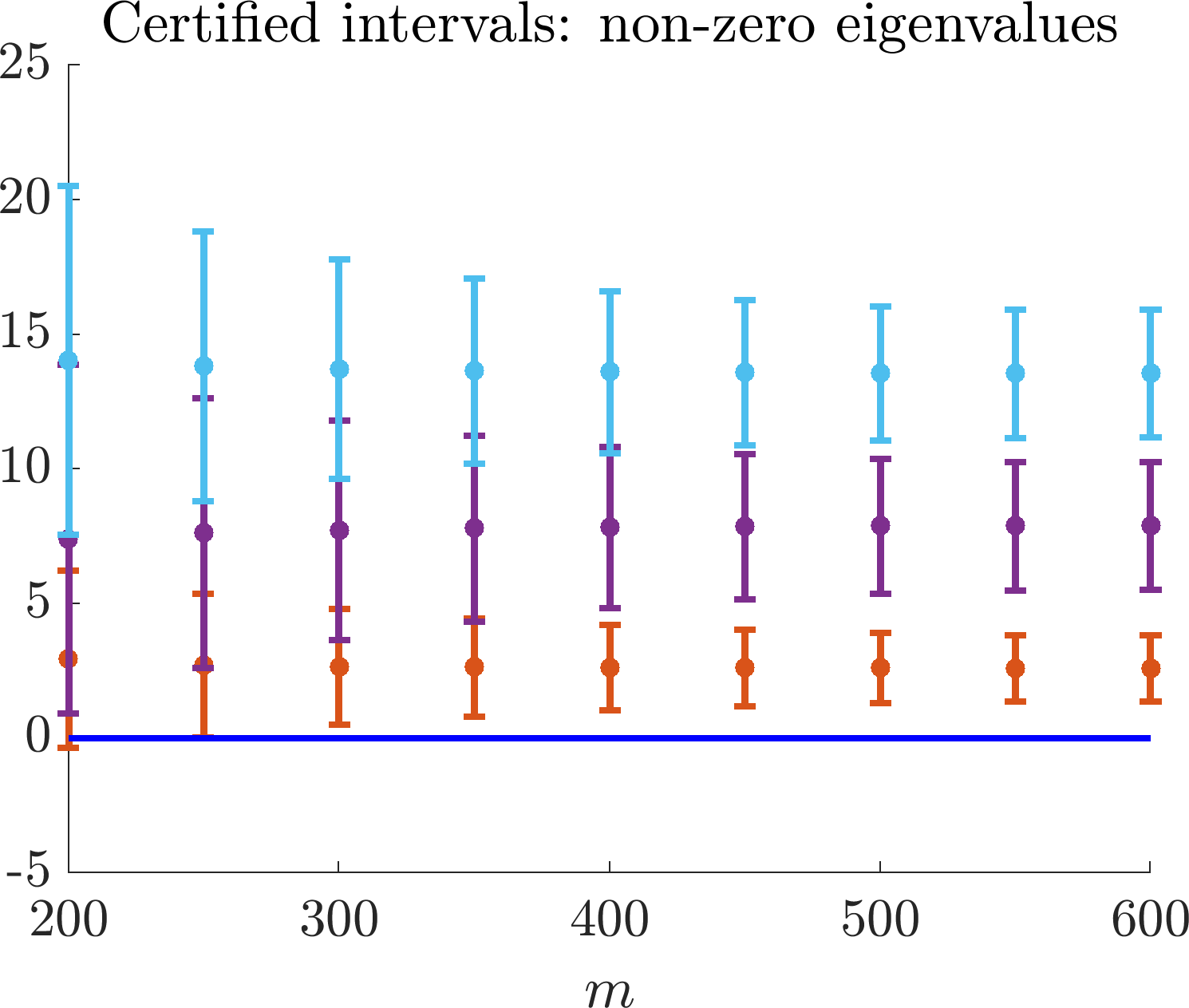}\quad
	\includegraphics[height=0.39\textwidth]{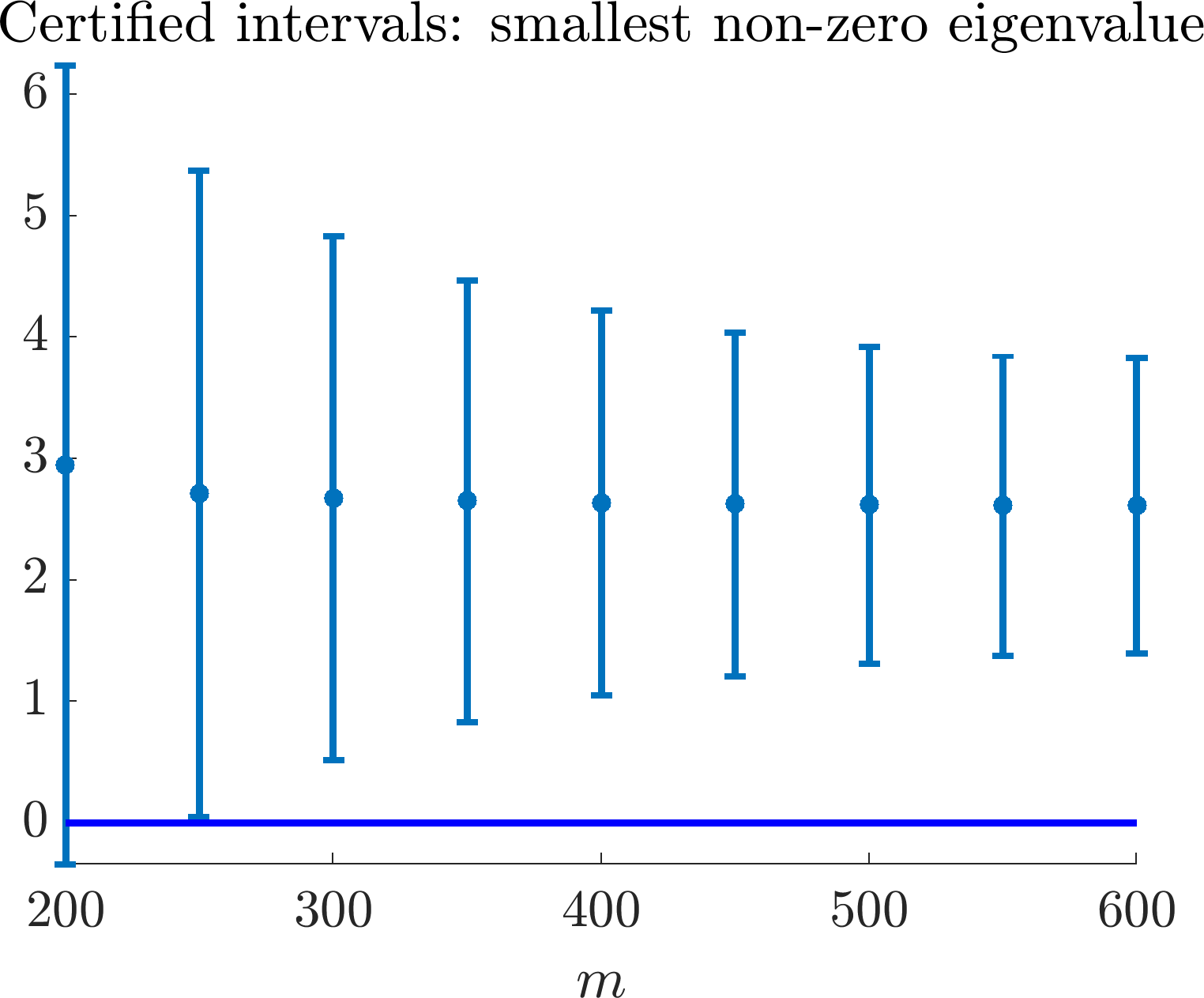}
	\caption{Computations for the \bo{regular pentagon}. Left: certified intervals containing the non-zero eigenvalues of the Hessian (three pairs of double eigenvalues). Right: certified intervals containing the smallest non-zero eigenvalue.}
	\label{fig:estimates-5} 
\end{figure}

\begin{figure}
	\centering
	\includegraphics[height=0.39\textwidth]{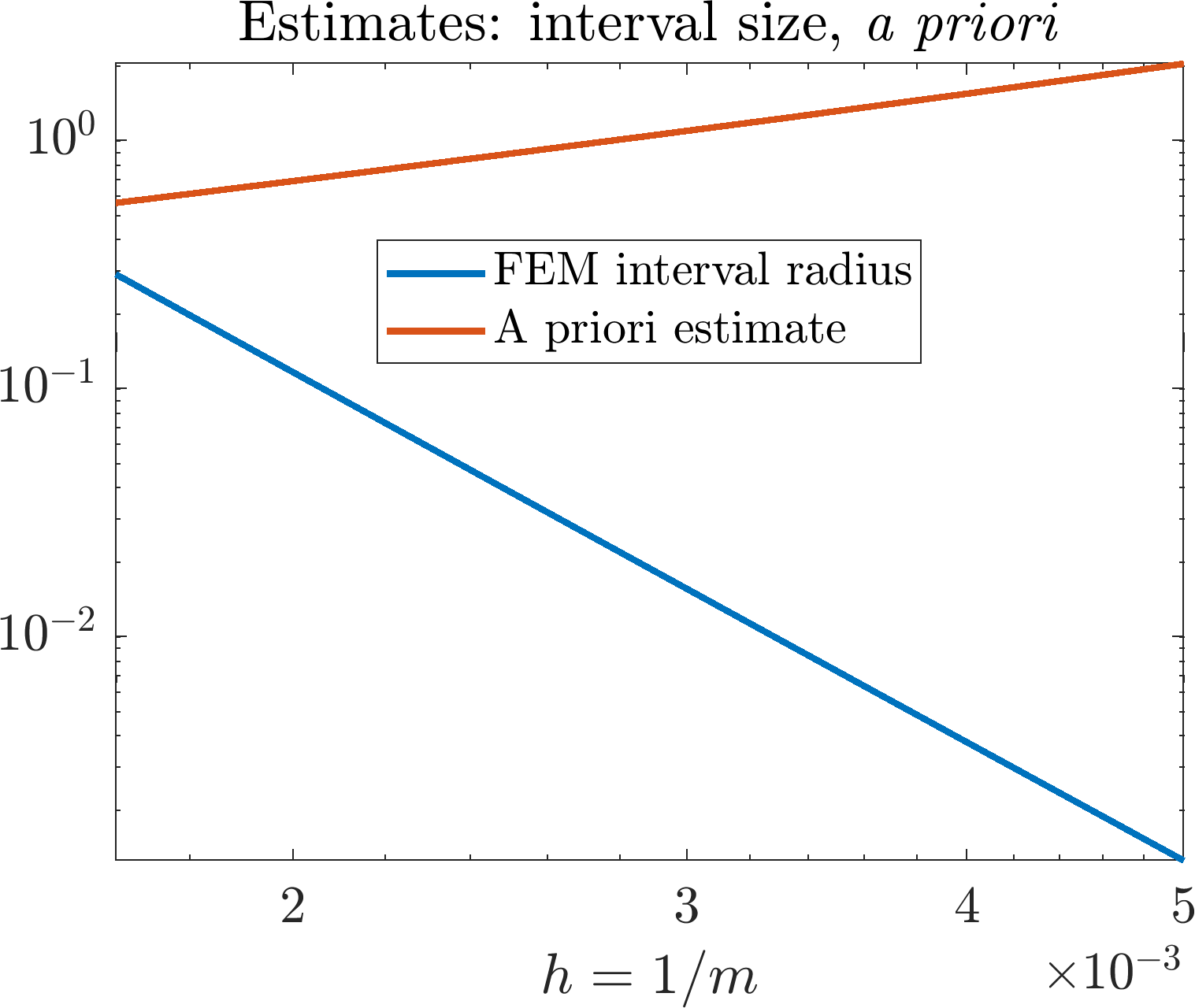}
	\caption{Computations for the \bo{regular pentagon}. The radius of the intervals for the FEM computation of the Hessian eigenvalues is shown with respect to the mesh size $h=1/m$, together with the corresponding \emph{a priori} estimate, whose leading term is $O(h)$.}
	\label{fig:FEMvsAnalytic}
\end{figure}

\bo{The regular pentagon.} The computations described in the previous section are performed for meshes corresponding to $m \in \{200,...,600\}$. It can be seen in Figure \ref{fig:estimates-5} that for $m$ large enough the $2n-4$ intervals containing the non-zero eigenvalues of the Hessian matrix (see Theorem \ref{thm:eig-hessian}) are contained in $(0,\infty)$, therefore these eigenvalues are guaranteed to be positive. This shows that the regular pentagon is a local minimizer for problem \eqref{eq:polya-conj}. The validation procedure succeeds already for $m=250$. Multiple observations are given below:

\bo{(a) Evolution of intervals enclosures with respect to $h=1/m$.} The size of the intervals stabilizes, even though the mesh becomes finer and finer. This is due to the quick growth in the interval size for the eigenvalues $\mu_j^h$ computed with finite elements (formulas from Theorem \ref{thm:eig-hessian} with finite element variables). Multiple factors contribute to this:
\begin{itemize}[noitemsep]
	\item The condition number of the rigidity matrix behaves like $O(h^{-2}) = O(m^2)$ (see \cite{cond_numbers}). Therefore, the performance of linear solvers becomes worse for finer meshes. This also affects the interval enclosures produced by Intlab \cite{intlab}.
	\item The linear systems \eqref{eq:U1-discrete}, \eqref{eq:U2-discrete} have wider right hand sides for fine meshes, therefore produce wider results. 
	\item The Hessian eigenvalues computations depend on integrals of the form $\int_{T_i} \partial_{x,y} u_{1,h} \cdot \partial_{x,y} U^{1,2}_h$. These evaluations produce intervals proportional to the number of triangles in the slice $T_i$, which is $O(m^2)=O(h^{-2})$. 
	\item The combined effect of all these aspects is illustrated in Figure \ref{fig:FEMvsAnalytic} where it can be seen that the interval size increases drastically when $h \to 0$. The \emph{a priori} estimate verify the theoretical $O(h)$ bound, while the interval sizes behave like $O(h^{-5})$.
\end{itemize}

\bo{(b) Evaluation of finite element quadratures with intervals.} Let us discuss in more detail the evaluation of finite element terms of the form $\int_{T_i} \partial_{x,y} u_{1,h} \cdot \partial_{x,y} U^{1,2}_h$ with intervals. It is apparent from Proposition \ref{prop:elements-xy} that the elements of the rigidity matrices used to compute these quantities do not depend on the scaling of the triangle. The area of the triangle is canceled when computing the integral of a term of size $1/h^2$. Therefore, assuming $u_h, U_{1,2}^h$ are computed as solutions of linear systems and eigenvalue problems, obtaining intervals of a given size $\rho$, assuming this size is not worsened when $h\to 0$ (this is often not the case), the interval enclosure for  $\int_{T_i} \partial_{x,y} u \cdot \partial_{x,y} U^{1,2}$ will be proportional to the number of triangles in the mesh of $T_i$ times $\rho$, giving $O(m^2 \rho) = O(\rho h^{-2})$ in our case. Therefore, the sole evaluation of the quadrature of a triangulation produces a drastic increase in the resulting intervals. 

Another interpretation for the increase in interval sizes when dealing with integrals of derivatives of $\bo P_1$ functions given with intervals, can be given as follows. If the discrete functions are given as intervals of size $\rho$ then an interval enclosure of the integral includes the worst case scenario when the function oscillates with amplitude $\rho$ from one vertex to the next one. Thus, errors are proportional to the number of oscillations, i.e. the number of triangles. 

For terms of the form $\int_{T_i} \nabla u_{1,h} \cdot \nabla U^{1,2}_h$, the fact that $u_{1,h}$ is a discrete eigenvalue should be used, to obtain the equivalent formula $\lambda_{1,h}\int_{T_i}  u_{1,h}  U^{1,2}_h$. This latter integral is computed using a mass matrix, whose elements are proportional to $h^2$. Therefore, there is no explosion in the interval size as the mesh increases for integrals not containing gradients. This observation motivates the use of integrals without derivatives, whenever possible, to decrease the size of the interval enclosures.

The previous observation can be seen as the interval arithmetic analogue of the \emph{a priori} estimates:
\begin{itemize}[noitemsep]
	\item The estimate $\|\nabla u-\nabla u_h\|_L^2$ is of order $h$, while $\|u-u_h\|_{L^2}=O(h^2)$, assuming $u \in H^2$.
	\item Assuming triangles of uniform area, and enclosures for finite element variables $u_h,v_h$ with uniform intervals with respect to $h$, the interval enclosure for $\int_T \nabla u_h \cdot \nabla v_h$ has a radius $O(h^{-2})$ bigger than the one for the enclosure for $\int_T u_h v_h$. 
\end{itemize}

\begin{figure}
	\centering
	\includegraphics[height=0.39\textwidth]{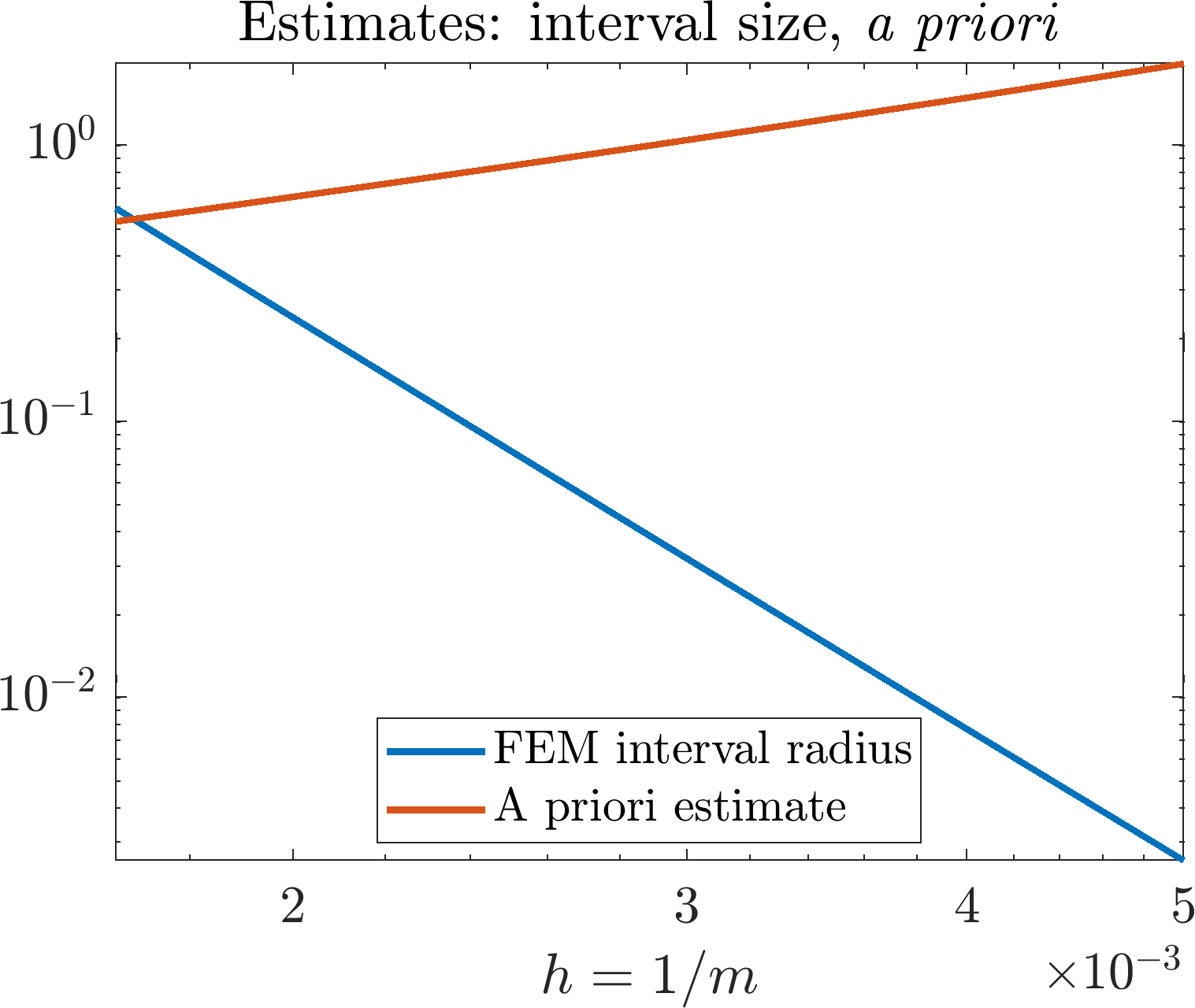}\quad
	\includegraphics[height=0.39\textwidth]{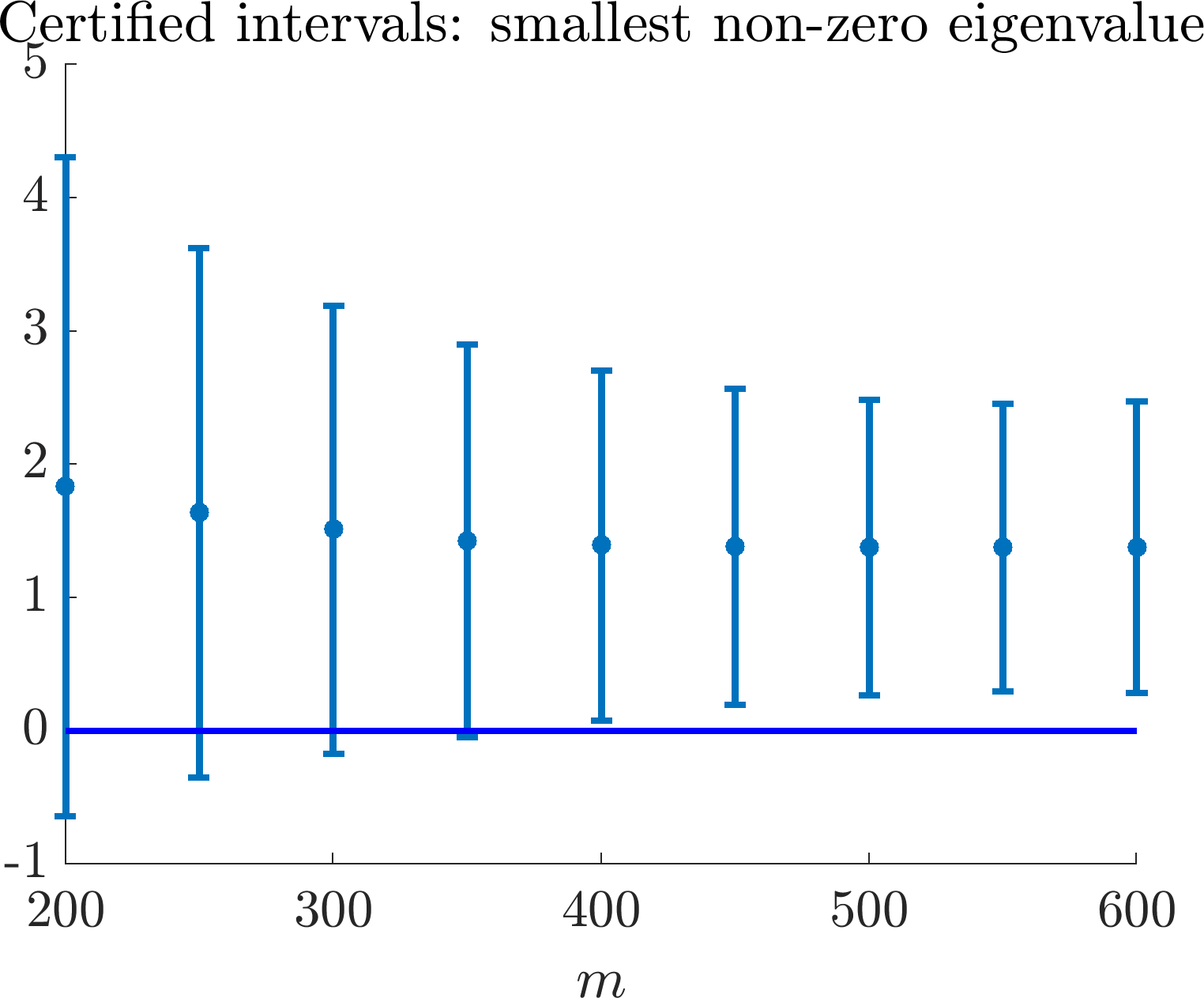}
	\caption{Computations for the \bo{regular hexagon}. Left: evolution of the \emph{a priori estimates} and the size of the intervals for $\mu_j^h$ coming from the FEM computations. Right: certified intervals containing the smallest non-zero eigenvalue. The validation of local minimality succeeds in this case.}
	\label{fig:estimates-6} 
\end{figure}

\bo{The regular hexagon.} Computations for the regular hexagon for $m \in \{200,250,...,600\}$ are shown in Figure \ref{fig:estimates-6}. The certification of local minimality using the method proposed in the previous section succeeds again for $m$ large enough: the certified lower bound for the smallest non-zero eigenvalue becomes positive. The validation procedure succeeds already for $m=380$. On the other hand, compared to the regular pentagon, for $m=600$ the radius of the intervals for the Hessian eigenvalues computed using finite elements becomes larger than the \emph{a priori} estimate. This shows that refining the mesh further will not improve the results. 

\begin{figure}
	\centering
	\includegraphics[height=0.39\textwidth]{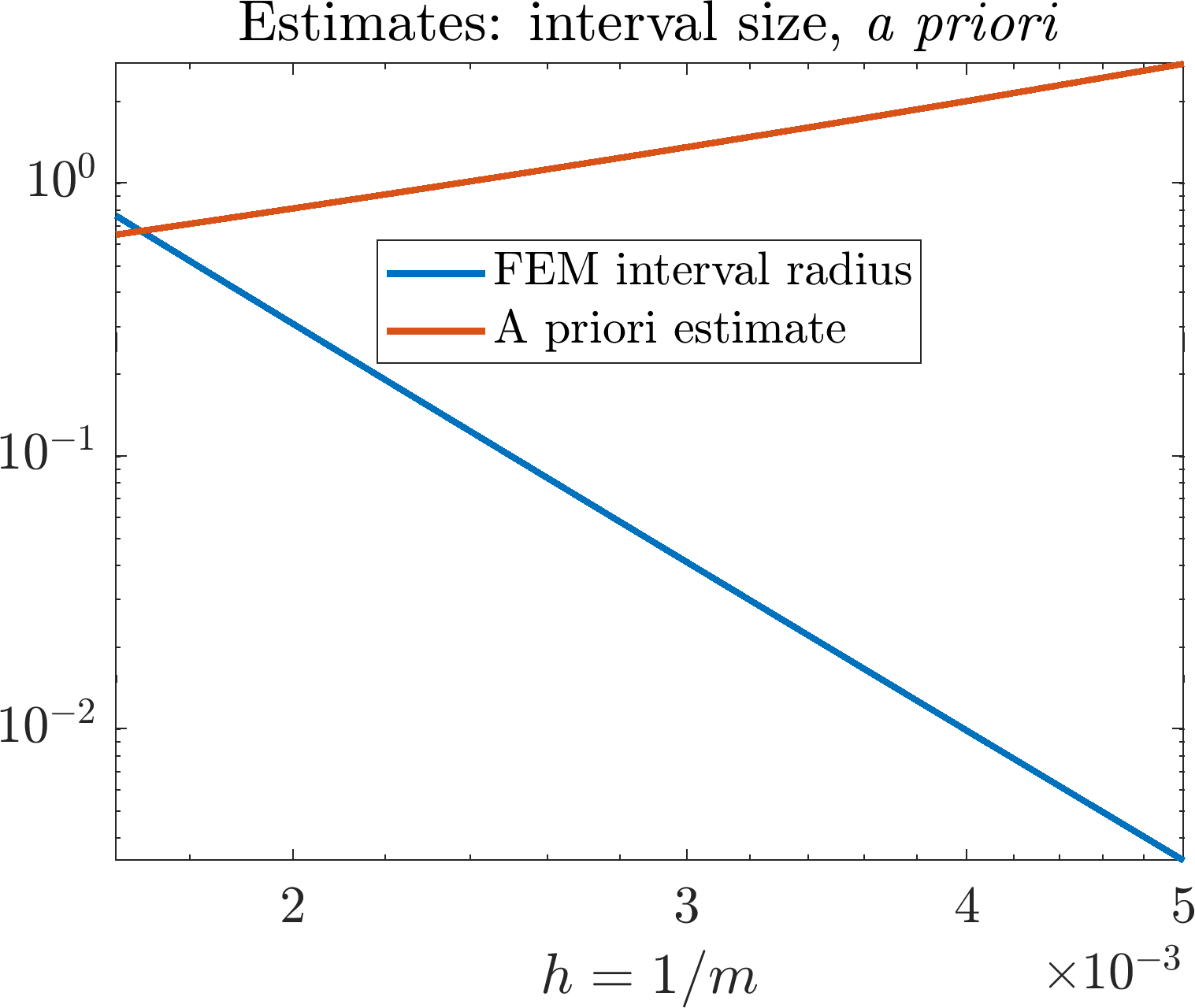}\quad
	\includegraphics[height=0.39\textwidth]{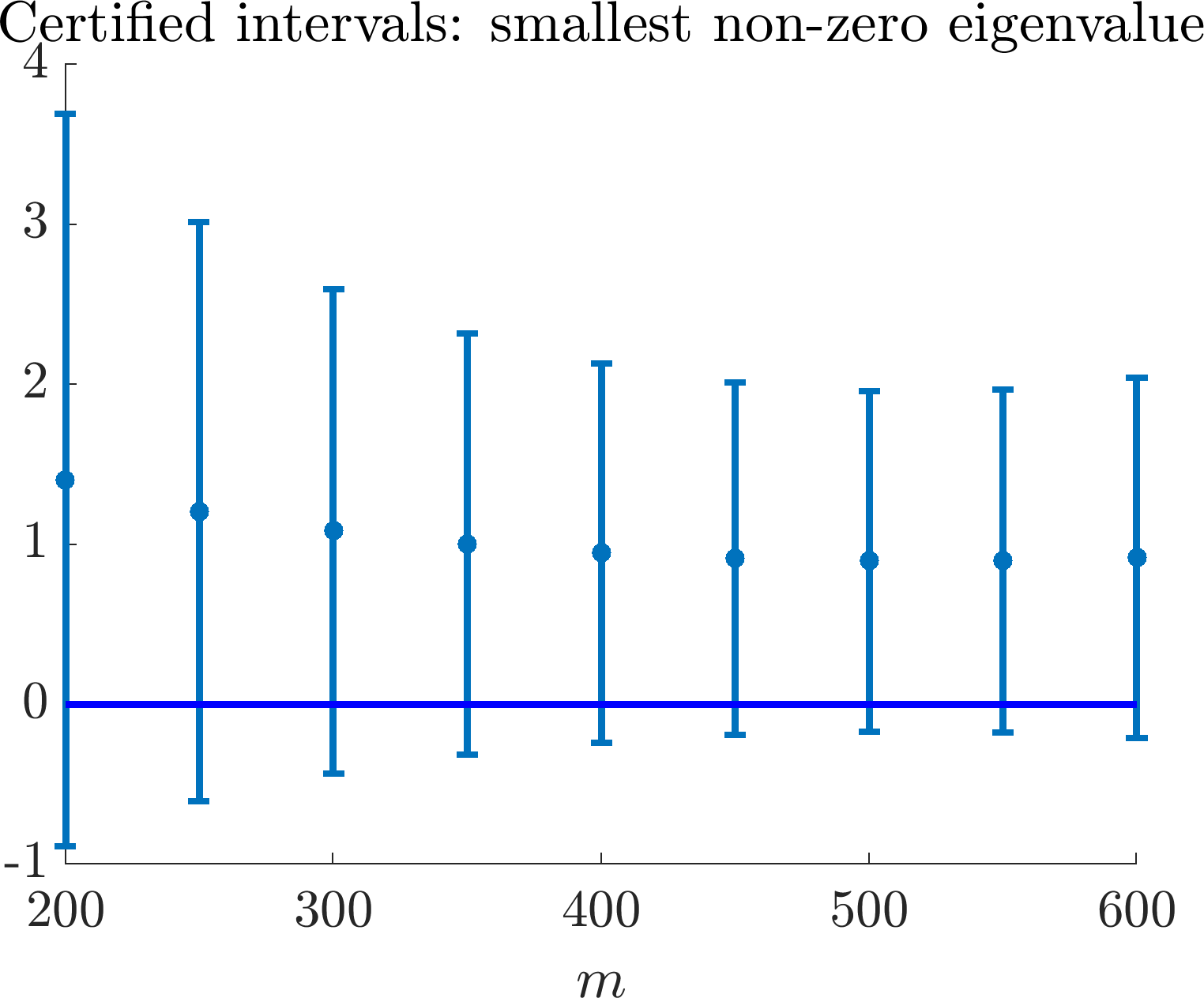}
	\caption{Computations for the \bo{regular heptagon}. Left: evolution of the \emph{a priori estimates} and the size of the intervals for $\mu_j^h$ coming from the FEM computations. Right: certified intervals containing the smallest non-zero eigenvalue. The validation of local minimality does not succeed in this case.}
	\label{fig:estimates-7} 
\end{figure}

\bo{The regular heptagon.} Computations for the regular heptagon for $m \in \{200,250,...,600\}$ are shown in Figure \ref{fig:estimates-7}. The certification of local minimality fails in this case. The radius of the intervals for the Hessian eigenvalues computed using finite elements becomes larger than the \emph{a priori} estimate for $m=600$. Moreover, it can be observed that the radius of the enclosure of the smallest eigenvalue increases for $m\geq 500$ even though the \emph{a priori} estimate continues to decrease. 

\begin{table}
	\begin{center}
		\begin{tabular}{c||c|c|c||c|c|c|c|}
			$n$ & \multicolumn{3}{c}{Results from \cite{Bogosel_Bucur_Polya}} & \multicolumn{4}{c}{Results of this paper} \\
			\hline  \hline 
			& $h$  & DoF & Intervals & $h=1/m$ & m & DoF & Intervals \\ 
			\hline
			$5$ & 9.8e-4 & 2.5 million & \ding{55} & 0.0040  & 250 & 156876 & \checkmark \\ \hline
			$6$ & 4.2e-4 & 17 million & \ding{55} & 0.0026 & 380 &  434341 & \checkmark\\ \hline
			$7$ & 1.9e-4 & 97 million & \ding{55}   & - & - & - &  \\ \hline
			$8$ & 1.35e-4 & 220 million & \ding{55} & - & - & - & \\ \hline
		\end{tabular}
	\end{center}
	\caption{Comparing the estimates in \cite{Bogosel_Bucur_Polya} with those of this paper. The improved error estimates allow us to use interval arithmetic in all computation, validating the local minimality for the regular pentagon and the regular hexagon.}
	\label{tab:comparison}
\end{table}

\bo{Comparison with previous results.} Results provided in \cite{Bogosel_Bucur_Polya} used weaker \emph{a priori} estimates. Meshes for which these estimates could validate the positivity of the non-zero eigenvalues of the Hessian matrix were too large to validate local minimality using intervals. The improvement of the constant in the interpolation constant \eqref{eq:interpolation} was also helpful to further decrease the computational cost. Table \ref{tab:comparison} shows a quick comparison of the results.  Mesh sizes are chosen in each case such that $h$ is as large as possible while \emph{a priori} estimates guarantee the positivity of the Hesian eigenvalues.

\subsection{Implementation}

The Matlab code implementing the steps described above allowing to certify the local minimality for $n \in \{5,6\}$ can be found at the following repository:

\begin{center}
\href{https://github.com/bbogo/PolyaSzego}{\nolinkurl{https://github.com/bbogo/PolyaSzego}}
\end{center}
It requires a working installation of Intlab \cite{intlab}. See Appendix \ref{app:code} for more details. 

\subsection{Conclusion, Perspectives}

The theoretical results and certified numerical simulations presented in this paper show the following result:

\begin{thm}\label{thm:loc-minimality}
	The regular $n$-gon is a local minimizer for the first Dirichlet-Laplace eiengvalue among $n$-gons with fixed area for $n\in \{5,6\}$. 
\end{thm}

Concerning perspectives for improving the results we mention the following:
\begin{itemize}[noitemsep]
	\item \bo{\emph{A posteriori} estimates.} Refining the error estimate for $\|\nabla u_1-\nabla u_{1,h}\|_{L^2(\Pp)}$ would directly influence the leading term for the error estimates discussed in Section \ref{sec:hess-eig-estimates}. Having better estimates would allow the validation procedure to succeed for meshes with fewer triangles, where using interval arithmetics for finite element computations will lead to tighter interval enclosures. The results in \cite{Vohralik2017} should be applied in this context.
	\item \bo{Higher order elements.} The context of $\bo P_2$ finite elements may be of interest since it gives improved error estimates. Results in \cite{xuefeng-quadratic} could help to adapt such techniques to our problem.
	\item \bo{Spectral methods.} The method of particular solutions was used in \cite{dahne-salvy}, \cite{dahne-payne} to compute enclosures for the Dirichlet Laplacian eigenvalues and eigenfunctions on polygons. Similar results should be applied to problems \eqref{eq:U0}, where the singular behavior should be taken into account in the choice of particular solutions. Methods described in \cite{trefethen} behave well in singular contexts. Investigating quantified error bounds for such methods is of interest concerning the results of this paper. 
\end{itemize}

\appendix

\section{Morley element and interpolation constants}
\label{app:morley}

Consider $T$ a triangle in the plane with vertices $x_1,x_2,x_3$ labeled in the positive trigonometric sense. Denote by $\gamma_1,\gamma_2,\gamma_3$ the edges opposite to $x_1,x_2,x_3$, respectively, having positive orientation. By definition, the Lagrange $\bo P_1$ interpolant $\Pi(u)$ of a function $u$ defined on a triangle $T$, is the affine function taking the same values as $u$ at the vertices of $T$. Therefore $\Pi(u)-u$ takes zero values at the vertices of $T$. Consider the space 
\[ V^2(T) = \{u \in H^2(T) : u(x_k) = 0, k=1,2,3\}.\]
Note that this space is well defined since functions in $H^2(T)$ are continuous, thus point values can be considered. In order to obtain the best constant in the inequality
\begin{equation}\label{eq:optimal-interpolation} \|\nabla (\Pi(u)-u)\|_{L^2} \leq C(T)|u|_{H^2},
\end{equation}
where $|u|_{H^2}^2 = \|\partial_{xx}u\|_{L^2(T)}^2+2\|\partial_{xy} u\|_{L^2(T)}^2+\|\partial_{yy} u\|_{L^2(T)}$ the following Rayleigh quotient is considered:
\[ C(T) = \sup_{u \in V^2(T)\setminus \{0\}} \frac{\|\nabla u\|_{L^2(T)}}{|u|_{H^2(T)}}.\]
 This constant is well defined and an upper bound can be found through classical estimates involving the eigenvalues of the Laplacian and the Poincar\'e inequality. Indeed, assume for contradiction that for some sequence $(u_n) \in H^2(T)$ we have $\|\nabla u_n\|_{L^2(T)}=1$ and $|u_n|_{H^2(T)} \to 0$. Let us denote $\overline u_n = \frac {1}{|T|} \int_T u_n$, so that we can apply the Poincar\'e-Wirtinger inequality in $H^1(T)$ to the functions $u_n -\overline u_n$, hence their $L^2$ norms are uniformly bounded. This means that the sequence $(u_n-\overline u_n)_n$ is bounded in $H^2$ and each function $v_n:= u_n-\overline u_n$ is constant on the vertices of $T$, while $\|\nabla v_n\|_{L^2(T)}=1$ and $|v_n|_{H^2(T)} \to 0$. We can assume that $v_n$ converges weakly in $H^2(T)$ to some function $v$. In particular, the convergence will also hold in $C^{0, \alpha}(T)$, for every $\alpha \in (0,1)$ and in $H^1(T)$-strongly. We get that $|v|_{H^2(T)} = 0$, so that $v$ is affine. But, $v$ is constant on the vertices of $T$, so that $v$ is a constant function, in contradiction with  the fact that $\int_T \|\nabla v\|^2 =1$.

  This constant is denoted $C_4(T)$ in \cite{kobayashi}, where it is shown that 
\begin{equation}\label{eq:kobayashi-estimate} C(T) \leq \sqrt{\frac{n^2}{n^2-1}}C^{(n)}(T),
\end{equation}
where $C^{(n)}(T)$ is defined as follows:
\begin{itemize}
	\item Given a triangle $\tau$ with vertices $p_k$ and edges $\gamma_k$, $i=1,2,3$ and $\varphi \in H^2(\tau)$ consider the interpolation operator $\Pi^{(\beta)}_\tau$ (using the same notations as in \cite{kobayashi}) such that: 
	
	a) $\Pi^{(\beta)}_\tau\varphi\in \{a_1x^2+a_2xy+a_3y^2+a_4x+a_5y+a_6 : a_i \in \Bbb{R}, i=1,...,6\}$.
	
	b) $\Pi^{(\beta)}_\tau\varphi(p_k) = \varphi(p_k),\ k=1,2,3$
	
	c) $\int_{\gamma_k} \nabla \Pi^{(\beta)}_\tau \varphi \cdot\bo  n = \int_{\Gamma_k} \nabla \varphi \cdot \bo n, \ k=1,2,3$.
	
	In other words, $\Pi^{(\beta)}_\tau \varphi$ is a second order polynomial in $x,y$ taking the same values as $\varphi$ at the vertices of $\tau$ and having the same average normal flux across the edges $\gamma_k$, $k=1,2,3$.
	
	This interpolation operator is also called the Morley interpolant. Since the gradient of a quadratic function is a linear function, the average of the normal derivative on an edge of $\tau$ is proportional to the normal derivative at the midpoint of the side. Thus, one can take as degrees of freedom for the Morley interpolant the values at the vertices and the normal derivatives at the midpoints of the sides. 
	
	\item Given $n \geq 1$ define $T'=\cup_{k=1}^{n^2}\tau_k$ the triangulation of $T$ into $n^2$ congruent triangles similar to $T$ (splitting each side of $T$ into $n$ equal segments like in Figure \ref{fig:mesh-slice-extension}). Given $u \in H^2(T)$ consider the interpolant $\Pi^{(\beta)}u$ such that
	\[ \Pi^{(\beta)}u|_{\tau_k} = \Pi_{\tau_k}^{(\beta)}u.\]
	In other words, on each sub-triangle $\tau_k$, $\Pi^{(\beta)}u$ is a quadratic polynomial, taking the same values as $u$ at the vertices of $\tau_k$ and having the same normal flux as $u$ across the sides of $\tau_k$. 
	
	\item The constant $C^{(n)}(T)$ is defined by 
	\begin{equation}\label{eq:opt-interp-discrete} C^{(n)}(T) =  \sup_{u \in V^2(T)\setminus 0} \frac{\|\nabla \Pi^{(\beta)}u\|_{L^2(T')}}{|\Pi^{(\beta)} u|_{H^2(T')}}.
	\end{equation}
	The notations $L^2(T')$ and $H^2(T')$ simply mean that functions are in the corresponding spaces for each member of the triangulation $T'$. The functions may not be globally $H^2$, for example.
\end{itemize}
The estimate \eqref{eq:kobayashi-estimate} states that the discrete eigenvalue obtained when using the Morley finite element on a triangulation of $T$ using equal triangles gives an upper bound for the constant in the interpolation estimate \eqref{eq:optimal-interpolation}. The proof of the estimate is a simple consequence of the identity
\[ |\varphi|^2_{H^2(\tau)} = |\Pi_\tau^{(\beta)}\varphi|^2_{H^2(\tau)}+|\varphi-\Pi_\tau^{(\beta)}\varphi|^2_{H^2(\tau)},\]
which follows from the definition of the Morley interpolant. In other words, the Morley interpolant decreases the $H^2$ semi-norm on a triangle. 

In order to have a tight estimate for the interpolation constant $C(T)$ it is enough to have a certified inclusion for $C^{(n)}(T)$, obtained using interval arithmetic. 

To achieve this, we show below an intuitive strategy for assembling the matrices associated to the Morley element, which can help compute $C^{(n)}(T)$. For alternative approaches one might consider \cite{kirby}. The strategy is as follows: 

\begin{itemize}[noitemsep]
	\item From the degrees of freedom for the Morley element (values at vertices, normal derivatives at midpoints of the edges) one can find the full gradient at the midpoints of the edges with an elementary observation regarding quadratic functions: if $P$ is quadratic on $T$ then $P(p_{k+1})-P(p_k)=2|p_{k+1}-p_k|\partial_\tau P\left(\frac{p_{k+1}+p_k}{2}\right)$.
	\item Gradient of a quadratic function on $T$ is affine on $T$. Knowing gradients at midpoints gives the gradients at vertices. 
	\item Knowing gradients as vertices as affine functions leads to easy assembly of the matrices involving first and second derivatives using classical ideas regarding $\bo P_1$ finite elements. 
\end{itemize}

In the following, the contribution of a triangle of the mesh to the assembly matrices is described. The procedure could be applied to a general mesh $\Omega$ where the degrees of freedom corresponding to Morley finite elements correspond to values at vertices and values on edges (corresponding, for example, to half the normal derivatives across the edges).

Supposing that $v_1,v_2,v_3$ are the values of an quadratic function $p$ at the vertex of $T$ and $\psi_{12},\psi_{31},\psi_{23}$ are half the normal derivatives at the midpoints of the sides $x_ix_j$ and $\ell_{ij} = |x_i-x_j|$ then the components of the gradient of $p$ at the vertices of $T$ can be found through the matrix product:
\[ 
\begin{pmatrix}
p_{1,x}\\
p_{1,y}\\
p_{2,x}\\
p_{2,y}\\
p_{3,x}\\
p_{3,y}
\end{pmatrix}
= \mathcal M_T
\begin{pmatrix}
v_1 \\
v_2 \\
v_3 \\
\psi_{23} \\
\psi_{31} \\
\psi_{12}
\end{pmatrix}
\]
where
\[\mathcal M_T = \begin{pmatrix}
-I_2 & I_2 & I_2 \\
I_2 & -I_2 & I_2 \\
I_2 & I_2 & -I_2
\end{pmatrix}
\begin{pmatrix}
\bo M_{23} & 0 & 0 \\
0 & \bo M_{31} & 0 \\
0 & 0 & \bo M_{12}
\end{pmatrix}
\begin{pmatrix}
0 & -\frac{1}{\ell_{23}} & \frac{1}{\ell_{23}} & 0 & 0 & 0 \\
0 & 0 & 0 & \sigma_{23} & 0 & 0 \\
\frac{1}{\ell_{31}} & 0 & -\frac{1}{\ell_{23}} & 0 & 0 & 0 \\
0 & 0 & 0 & 0 & \sigma_{31} & 0  \\
-\frac{1}{\ell_{12}} & \frac{1}{\ell_{12}} &  0 & 0 & 0 & 0 \\
0 & 0 & 0 & 0 & 0 &\sigma_{12} \\	
\end{pmatrix}.\]
The unitary matrices $\bo M_{ij} = \frac{1}{\ell_{ij}}[ (x_j-x_i) | (x_j-x_i)^{\perp}]$ rotate a coordinate system aligned with $x_j-x_i$ to make it align with the one corresponding to the canonical basis. Tangential and normal derivatives of $p$ at midpoints are transofrmed into partial derivatives $\partial_x p, \partial_y p$. The values $\sigma_{ij}\in \{-1,1\}$ are arbitrary (but fixed) sign choices for the flux on the edge $[x_i,x_j]$ giving an orientation for the flux, opposite for the two triangles adjacent to this side. 

Next, consider the permutation matrix $P$ which regroups on the first three components the values of $\partial_x p$ and on the last three the values of $\partial_y p$.
\[ \begin{pmatrix}
p_{1,x}\\
p_{2,x}\\
p_{3,x}\\
p_{1,y}\\
p_{2,y}\\
p_{3,y}
\end{pmatrix}= P\mathcal M_T \begin{pmatrix}
v_1 \\
v_2 \\
v_3 \\
\psi_{23} \\
\psi_{31} \\
\psi_{12}
\end{pmatrix}.\]
Next, consider the mass and rigidity matrices $K_T,M_T$ for $\bf P_1$ finite elements on the triangle $T$. 

Denoting by $\bo v$ a generic vector containing all degrees of freedom for the Morley elements (one per vertex and one per edge) we have the following. Assume the domain $\Omega$ is triangulated with triangles $\{T_i\}_{i=1}^N$.

A) The matrix $\mathcal K_{xx}$ which for a Morley function $\tilde u$ with degrees of freedom $\bo v$ computes $$\bo v^T \mathcal K_{xx} \bo v =  |\tilde u|^2_{H^2(T')}=\sum_{i=1}^N \int_{T_i} (\partial_{xx}\tilde u)^2+2(\partial_{xy}\tilde u)^2+(\partial_{yy}\tilde u)^2=\sum_{i=1}^N \|\nabla p_x\|_{L^2(T_i)}^2+\|\nabla p_y\|_{L^2(T_i)}^2,$$
where $p_x,p_y$ denote degree one polynomials corresponding to $\partial_x \tilde u, \partial_y \tilde u$ on the triangles of the mesh, is given as follows.  For each triangle in the mesh, the contribution corresponding to the degrees of freedom associated to the values and fluxes along the elements of the triangle is given by
\[   \mathcal M_T^T P^T \begin{pmatrix}
K_T & 0 \\
0   & K_T
\end{pmatrix} P \mathcal M_T,\]
where $K_T$ is the rigidity matrix for $\bo P_1$ elements on the triangle $T$.

B) The matrix $\mathcal M_{xx}$ which for a Morley function $\tilde u$ with degrees of freedom $\bo v$ computes $$\bo v^T \mathcal M_{xx} \bo v =  \|\nabla \tilde u\|^2_{L^2(T')}=\sum_{i=1}^N \int_{T_i} (\partial_{x}\tilde u)^2+(\partial_{y}\tilde u)^2=\sum_{i=1}^N \|p_x\|_{L^2(T_i)}^2+\|p_y\|_{L^2(T_i)}^2$$ is given as follows. Then for each triangle in the mesh, the contribution corresponding to the degrees of freedom associated to the values and fluxes along the elements of the triangle is given by
\[   \mathcal M_T^T P^T \begin{pmatrix}
M_T & 0 \\
0   & M_T
\end{pmatrix} P \mathcal M_T,\]
where $M_T$ is the mass matrix for $\bo P_1$ elements on the triangle $T$.

The assembly procedure described above may not be the most efficient, but it has the advantage of being clear and explicit. Moreover, for the triangulation of the original triangle $T$ into $n^2$ congruent triangles the assembly can be rendered explicit and evaluated with INTLAB with a tight interval enclosure. The constant $C^{(n)}(T)$ defined in \eqref{eq:opt-interp-discrete} is given by
\[ C^{(n)}(T) = \sqrt{1/\rho_1}, \text{ where } \rho_1 = \inf_{\bo v, \bo v(x_k)=0}\frac{\bo v^T \mathcal K_{xx} \bo v }{\bo v^T\mathcal M_{xx} \bo v},\]
which means that $\rho_1$ is the smallest generalized eigenvalue for the matrices $(\mathcal K_{xx},\mathcal M_{xx})$, taking into account the Dirichlet boundary conditions $\bo v(x_1)=\bo v(x_2)=\bo v(x_3)=0$. To impose these boundary condition, the corresponding lines and columns are eliminated from $(\mathcal K_{xx},\mathcal M_{xx})$ obtaining $(\mathcal K_{xx}^0,\mathcal M_{xx}^0)$.

Given $\overline \rho$, the smallest generalized eigenvalue of  $(\mathcal K_{xx}^0,\mathcal M_{xx}^0)$ approximated using a floating point computation, the routine \texttt{isspd} from INTLAB verifies if the matrix
\[ \mathcal K_{xx}^0-(\overline \rho-\varepsilon)\mathcal M_{xx}^0\]
is symmetric positive definite for some $\varepsilon>0$ (for example $\varepsilon=10^{-6}$). If the numerical validation succeeds, an upper bound for $C^{(n)}(T)$ is obtained, which gives an upper bound for $C(T)$ from \eqref{eq:kobayashi-estimate}.

Observing that the constant $C(T)$ defined in \eqref{eq:optimal-interpolation} scales linearly with the size of $T$, it is enough to compute it for a triangle with vertices $x_1=(0,0),x_2=(1,0),x_3=(a,b)$, verifying $b>0$. In this case, denoting with $\psi_1,\psi_2,\psi_3$ the $\bo P_1$ finite element functions on $T$ it follows immediately that
\[ \begin{array}{lll}
\partial_x \psi_1 = -1 & \partial_x \psi_2 = 1 & \partial_x \psi_3 = 0 \\
\partial_y \psi_1 = \frac{a-1}{b} & \partial_y \psi_2 = -\frac{a}{b} & \partial_y \psi_3 = \frac{1}{b}.
\end{array} \]
The expression of the rigidity matrix $K_T$ on the triangle $T$ follows, observing that the area of $T$ is equal to $b/2$. For the mass matrix one obtains
\[ M_T = |T| \begin{pmatrix}
1/6 & 1/12& 1/12 \\
1/12 & 1/6 & 1/12 \\
1/12 & 1/12 & 1/6
\end{pmatrix}.\]

The code which for a given pair $(a,b)$, $b>0$ gives a certified upper bound for the optimal interpolation constant $C(T)$ can be found at the following repository:

\begin{center}
	\href{https://github.com/bbogo/PolyaSzego/tree/main/Morley}{\nolinkurl{https://github.com/bbogo/PolyaSzego/tree/main/Morley}}
\end{center}

The code requires a working installation of INTLAB \cite{intlab}. It receives as inputs the vertex coordinates $(a,b)$ and the number of segments $m$ in which every side is divided for constructing the mesh. It outputs a certified upper bound for $C(T)$ if the validation method succeeds. Otherwise, it outputs a message indicating that a lower bound was not found. 

All results described in \cite{kobayashi} can be verified using the code proposed above. The estimates relevant to the contents of this paper correspond to isosceles triangles with edge length $1$ and central angle $2\pi/n$, giving $a = \cos(2\pi/n), b = \sin(2\pi/n)$, $n \in \{5,6\}$. The corresponding constants are summarized in Table \ref{tab:interp-const}, where dependence on $C(T)$ the central angle of the regular polygon is emphasized.
\begin{table}
	\begin{tabular}{|c|c|c|c|c|c|c|}
		\hline
       $\theta$ & $2\pi/5$ &  $2\pi/6$ & $2\pi/7$ & $2\pi/8$ & $2\pi/9$ & $2\pi/10$  \\ \hline 
       $C(T)$   & $0.3697$ &  $0.3200$ & $0.3146$ & $0.3107$ & $0.3104$ & $0.3128$ \\ \hline 
	\end{tabular}
\vspace{0.1cm}
\caption{Certified upper bounds for $\bo P_1$ interpolation constants on isosceles triangles corresponding to a slice of a regular $n$-gon for $n\in \{5,6,7,8,9,10\}$.}
\label{tab:interp-const}
\end{table}

\section{Description of the code}
\label{app:code}

The Matlab code proving the local minimality of the regular $n$-gon for problem \eqref{eq:polya-conj} can be found at the following repository:

\begin{center}
	\href{https://github.com/bbogo/PolyaSzego}{\nolinkurl{https://github.com/bbogo/PolyaSzego}}
\end{center}

It uses interval arithmetic and requires a working installation of INTLAB \cite{intlab}. Details regarding the implementation can be found in the code comments. For the sake of completeness, a description of the functions and of the validation process is described below. Before running the code, the subfolder \texttt{./Tools/} should be added to the Matlab path. 

\smallskip
\noindent \bo{a) Validation for the first two eigenvalues and the first eigenfunction.} The implementation is based on the discussion in Section \ref{sec:validation}. The first main function is \texttt{PolyaHessInterval.m}. For running the code corresponding to $n=5$ and $m=250$ division points on the ray $[\bo o \bo a_0]$ type: 
\begin{lstlisting}
>> PolyaHessInterval(5,250)
\end{lstlisting}
The code constructs two meshes, one for the slice $T_+$ (see Figure \ref{fig:simple-triangulations}), where the first eigenfunction is computed and validated, exploiting the symmetry. The first eigenpair is either validated using \texttt{verifyeig} in INTLAB or with a residual estimate. A second mesh of the full polygon is used for computing and validating the second eigenvalue. The validation uses the residual estimate from Proposition \ref{prop:validation-eigs}. At the end of the validation process the eigenvalue and eigenfunction enclosures are saved to a file to for a later use. 

The matrices needed in the computation are assembled exactly. Only the incidence relations in the mesh are used in the assembly, all other aspects are computed analytically. 

\smallskip
\noindent  \bo{b) Validating the local minimality of the regular $n$-gon.} The remaining of the validation procedure can be achieved using the function \texttt{PolyaHessIntervalU.m}. For running the code corresponding to $n=5$ and $m=250$ division points on the ray $[\bo o \bo a_0]$ type:
\begin{lstlisting}
>> PolyaHessIntervalU(5,250)
\end{lstlisting}

\begin{lstlisting}[caption={Output of the numerical validation code for $n=5$, $m=250$. Six out of the ten eigenvalues are positive, proving the local minimality of the regular pentagon.},captionpos=b,label=lst:output]
Number of positive eigenvalues = 6
Proof of local minimality succeeded!
Degrees of Freedom (full mesh) 156876
intval ans = 
[   -2.8083,    2.5145] 
[   -2.8083,    2.5145] 
[   -0.0001,    0.0001] 
[   -0.0001,    0.0001] 
[    0.0547,    5.3775]
[    0.0547,    5.3775] 
[    2.6266,   12.6594] 
[    2.6266,   12.6594] 
[    8.8146,   18.8474] 
[    8.8146,   18.8474] 
\end{lstlisting}

The linear systems associated to the material derivatives \eqref{eq:U1-discrete}, \eqref{eq:U2-discrete} are solved and validated and an interval enclosure for the solution is found. Conjugate gradient is used to obtain a floating point solution, which is validated using Lemma \ref{lem:saddle}. Afterwards, the eigenvalues of the Hessian matrix are evaluated from the finite element results, following Theorem \ref{thm:eig-hessian}. Next, all elements involving the \emph{a priori} estimates in Sections \ref{sec:prelim}, \ref{sec:material} quantifying the error between continuous solutions of the PDEs involved and the finite element counterparts are computed.

The code outputs a sequence of interval enclosures for the eigenvalues of the Hessian matrix. If $2n-4$ of these intervals are included in $[0,+\infty)$ then the validation procedure succeeds and local minimality is proved (see \cite[Section 4]{Bogosel_Bucur_Polya}). The local minimality validation succeeds for $n=5$ (taking $m=250)$ and $n=6$, (taking $m=380$). An example of output for $n=5$, $m=250$ is shown in Listing \ref{lst:output}.

\section{Fast validation of eigenvalues for large sparse matrices}	
\label{app:intlab}

Applying directly the routine \texttt{verifyeig} from Intlab \cite{intlab} to generalized eigenvalue problems of the form 
\[ K_0 u=\lambda M_0u\]
leads to significant computational costs when $K_0, M_0$ are of large sizes. Looking at the code provided with Intlab, it seems that the large cost comes from the inversion of a matrix. While it is true that multiple linear systems need to be soved using the matrix in question with various right hand sides, for large problems it is more efficient to consider linear systems instead. The function \texttt{verifylss} in Intlab allows to solve linear sparse systems of large sizes, therefore, this routine is used whenever a linear system is solved. The modifications are underlined in the Listings \ref{lst:1}, \ref{lst:2}. For the generalized eigenvalue problem the matrix \texttt{B} and its midpoint \texttt{midB} appear in the same places as in the original \texttt{verifyeig} function. 
\begin{center}
\begin{minipage}{0.47\textwidth}
\begin{lstlisting}[caption={The original \texttt{verifyeig} function in Intlab.},captionpos=b,label=lst:1]
.....
R = midA - lambda*speye(n);
R(:,v) = - xs;
"\colorbox{yellow}{R = inv( R );}" % matrix inversion
C = A - intval(lambda)*speye(n);
"\colorbox{yellow}{Z = - R * ( C * xs );}"
C(:,v) = - xs;
"\colorbox{yellow}{C = speye(n) - R * C;}"
Y = Z;
Eps = 0.1*mag(Y)*hull(-1,1) + midrad(0,realmin);
m = 0;
mmax = min( 15 * ( sum(sum(mag(Z(v,:))>.1)) + 1 ) , 20 );
ready = 0;
while ( ~ready ) && ( m<mmax ) && ( ~any(isnan(Y(:))) )
  m = m+1; 
  X = Y + Eps;
  XX = X;
  XX(v,:) = 0;
  "\colorbox{yellow}{Y = Z + C*X + R*(XX*X(v,:));}"
  ready = all(all(in0(Y,X)));
end
.....
\end{lstlisting}
\end{minipage}\quad
\begin{minipage}{0.47\textwidth}
\begin{lstlisting}[caption={Replacing matrix inversion with linear systems.},captionpos=b,label=lst:2]
.....
R = midA - lambda*speye(n);
R(:,v) = -xs;
% without matrix inversion                   
C = A - intval(lambda)*speye(n); 
"\colorbox{yellow}{Z = -verifylss(R,(C*xs));}" 
C(:,v) = -xs;

Y = Z;
Eps = 0.1*mag(Y)*hull(-1,1) + midrad(0,realmin);
m = 0;
mmax = min( 15 * ( sum(sum(mag(Z(v,:))>.1)) + 1 ) , 20 );
ready = 0;
while ( ~ready ) && ( m<mmax ) && ( ~any(isnan(Y(:))) )
  m = m+1; 
  X = Y + Eps;
  XX = X;
  XX(v,:) = 0;
  "\colorbox{yellow}{Y = Z+verifylss(R,(R-C)*X+XX*X(v,:));}"
  ready = all(all(in0(Y,X)));
end
.....
\end{lstlisting}
\end{minipage}
\end{center}

The modified \texttt{verifyeig} is applied to randomly generated matrices with sizes up to $1000\times 1000$ and compared with the original version. The results produced are identical (intervals have the same center and same radius). For full matrices the proposed modification is not faster than the original version for the cases tested. 

\begin{figure}
	\centering 
	\includegraphics[width=0.5\textwidth]{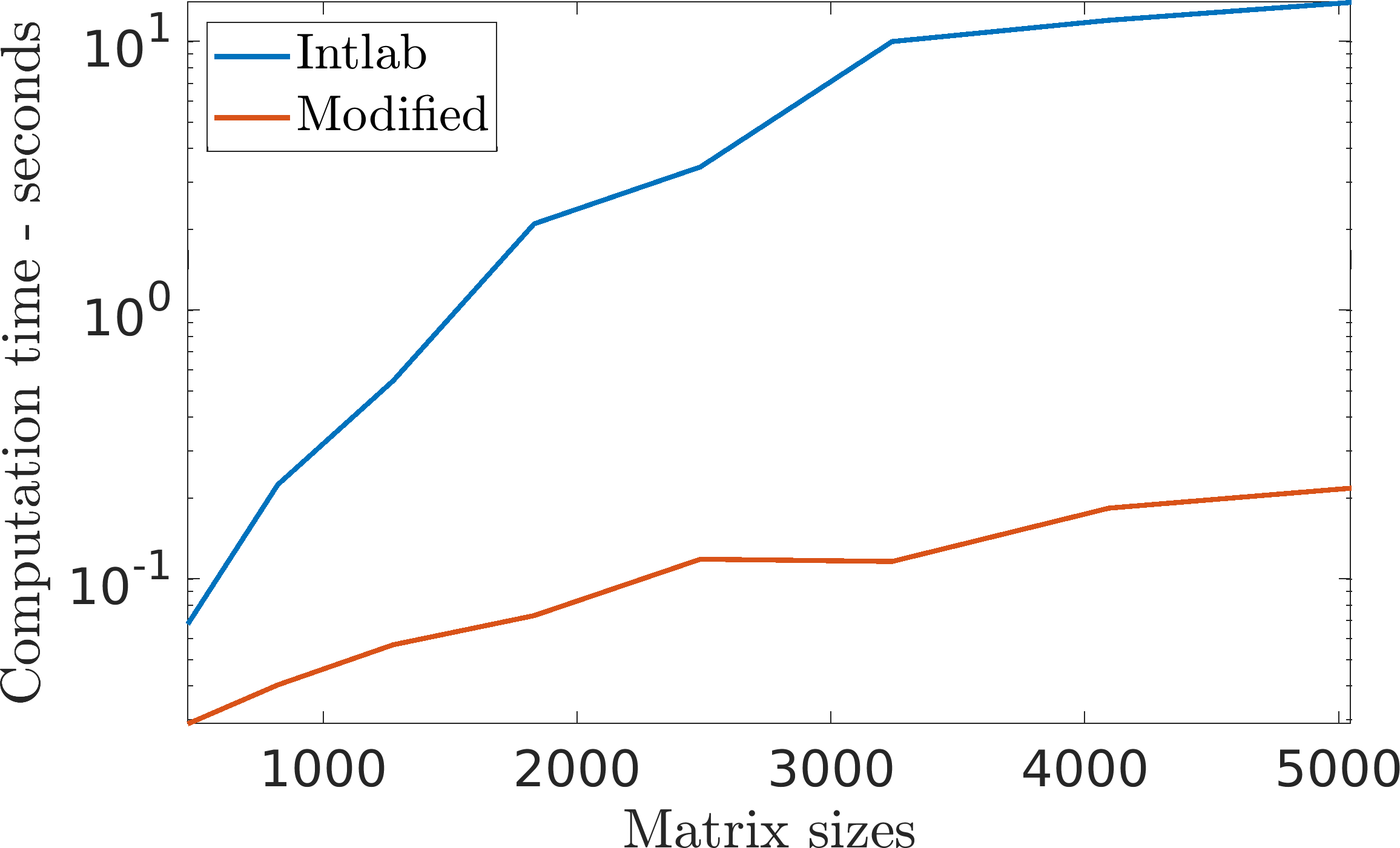}
	\caption{Comparison of computation times for the eigenvalue verification routine \texttt{verifyeig} in Intlab and the modified version, applied to sparse matrices coming from the discretization of a Laplace PDE. The significant time improvement allows to validate eigenvalues and eigenvectors for matrices of sizes up to $10^5 \times 10^5$.}
	\label{fig:intlab-comparison}
\end{figure}

Nevertheless, the situation changes when the matrices involved are sparse and come from discretization of finite element problems like \eqref{eq:eigenvalue-finite-elements}. The generalized eigenvalue problem for a triangular slice meshed with congruent triangles obtained by dividing each one of the sides into $m$ segments ($m^2$ triangles) is solved with the original version of \texttt{verifyeig} and the modified version. The computations are realized for $m\in \{30,40,...,90,100\}$ and the difference in computation times is shown in Figure \ref{fig:intlab-comparison}. The speed gain from removing the matrix inverse is significant, since Intlab uses a sparse version of the linear system verification routine \texttt{verifylss}. The verification with the modified version was tested successfully for triangle division parameter $m$ up to $450$, corresponding to matrix size of roughly $10^5 \times 10^5$. 

For even finer meshes, further modifications are provided to the \texttt{verifyeig} routine. An initial floating point system used is replaced with an iterative method. This modification was tested for the mesh parameter $m$ up to $600$ corresponding to matrices of size $180901\times 180901$. These modified versions of \texttt{verifyeig} could be of potential interest for other applications involving eigenvalues of large sparse matrices.

\smallskip
\noindent \bo{Acknowledgements:} The authors thank Maxime Breden, Xuefeng Liu and Siegfried M. Rump for valuable discussions and advice regarding the validation procedure using interval arithmetic. 

    \bibliographystyle{abbrv}
	\bibliography{./biblio.bib}
	
\end{document}